\documentclass[oneside,reqno,a4paper,11pt]{amsart}
\usepackage{amsfonts}
\usepackage{amsmath,amsfonts,amssymb,amscd}
\usepackage{color,longtable,verbatim}
\usepackage{color}
\usepackage{resizegather}
\usepackage{exscale}
\usepackage{latexsym, amsfonts,mathrsfs, amsmath, amssymb, amscd, epsfig}
\usepackage{mathrsfs}
\usepackage{amsthm}
\usepackage{etex}
\usepackage{epstopdf}
\usepackage{graphicx}
\usepackage{color}
\usepackage{ifpdf}
\usepackage{multirow}
\usepackage{float}
\usepackage{booktabs}
\usepackage{tabularx}
\usepackage{threeparttable}

\newtheorem {theorem}{Theorem}[section]

\newtheorem{theorem 1}{\textup{\textbf{Theorem}}}
\newtheorem {lemma}[theorem]{{\bf Lemma}}

\newtheorem {proposition}[theorem]{{\bf Proposition}}
\theoremstyle{remark}
\newtheorem {remark}{{\bf Remark}}[section]

\theoremstyle{problem}

\theoremstyle{definition}
\newtheorem {definition}{{\bf Definition}}[section]
\theoremstyle{plain} \numberwithin {equation}{section}
\newtheorem*{definitionA}{{\bf Definition A}}
\newtheorem*{definitionB}{{\bf Definition B}}
\newtheorem*{definitionC}{{\bf Definition C}}

\newenvironment{pf}{{\noindent \it \bf Proof.}}{{\hfill$\Box$}\\}
\newenvironment{pf2.1}{{\noindent \it \bf Proof of Theorem 2.1.}}{{\hfill$\Box$}\\}
\usepackage{fancyhdr}
\begin{document}
	\fancyhf{}
	\fancyhead[C]{FREE BOUNDARIES FOR EHD EQUATIONS WITH A GRAVITY FIELD}
	\fancyhead[R]{\thepage}
	\vspace{1cm}

	\title{Free boundary problem for two-dimensional ElectroHydroDynamic Equations with a gravity field}
	
	\author[Yuanhong Zhao,\ \ Lili Du]{Lili Du$^{\lowercase {1,2}}$, \ \ Yuanhong Zhao$^{\lowercase {2}}$}
\thanks{*This work is supported by National Nature Science Foundation of China under Grant 12125102, 12526202, Nature Science Foundation of Guangdong Province under Grant 2024A1515012794, and Shenzhen Science and Technology Program (JCYJ20241202124209011). }
\thanks{ E-mail: dulili@scu.edu.cn.  E-mail:
	zhaoyuanhong\_math@163.com. Corresponding author}\maketitle
\begin{center}
	$^1$ School of Mathematical Sciences, Shenzhen University,
	
	Shenzhen 518061, P. R. China.
	
		$^2$ Department of Mathematics, Sichuan University,
	
	Chengdu 610064, P. R. China.
	
\end{center}
\begin{abstract} 
	In this paper, we consider a two-phase free boundary problem governed by the ElectroHydroDynamic equations, which describes a perfectly conducting, incompressible, irrotational fluid with gravity, surrounded by a dielectric gas. The interface separating fluid and gas is referred to as the free boundary.
	
	It is well-established that the free surface remains smooth away from the stagnation points, where the relative velocity of the incompressible fluid vanishes. In the presence of gravity, the Stokes conjecture, proved by V\v{a}rv\v{a}ruc\v{a} and Weiss [Acta. Math. 206, 363-403, (2011)], implies that the corner type singularity will occur in the one-phase incompressible fluid. It is natural to ask whether this conjecture still holds in two-phase flow problem. As a consequence, the primary objective of this work is to characterize the possible singular profiles of the free interface near the stagnation points $x^{\text{\rm st}}$ in the presence of an electric field. 
	
	Our main result is the discovery of a critical decay rate $\left|x-x^{\text{\rm st}}\right|^{1/2}$ of the electric field near the stagnation points which indicates the classification of the singular profiles of the free surface. More precisely, we showed that when the decay rate of the electric field is faster than $\left|x-x^{\text{\rm st}}\right|^{1/2}$, its negligible effect implies that the singular profile must be the well-known Stokes corner, which is symmetric about the direction of gravity and has an opening angle of $120^\circ$. When the electric field decays as $\left|x-x^{\text{\rm st}}\right|^{1/2}$, the symmetry of the corner region may be broken, giving rise to either a Stokes corner or an asymmetric corner as the possible singular profile. If the decay rate is slower than $\left|x-x^{\text{\rm st}}\right|^{1/2}$, the electric field dominates and completely destroys the corner structure, resulting in a cusp singularity.
	
	 The analysis of these singularities relies on variational principles and geometric methods. Key technical tools include a Weiss-type monotonicity formula, a frequency formula, and a concentration-compactness argument.
\end{abstract}

\


\tableofcontents

\section{Introduction and main results}

\subsection{Introduction}
\quad

In this paper, we study a free boundary problem governed by ElectroHydroDynamic (abbreviated EHD in the following) equations in the presence of gravity. We consider a  model consisting of a stationary, irrotational, incompressible, perfectly conducting fluid under gravity surrounded by a dielectric gas (see Fig. 1). Our objective is to investigate the possible singular profiles of the interface between the perfectly conducting fluid and the dielectric gas near the stagnation points. {\textit {A stagnation point}} is defined as a point on the free boundary where the velocity field of the fluid vanishes. 
\begin{figure}[!h]
	\includegraphics[width=90mm]{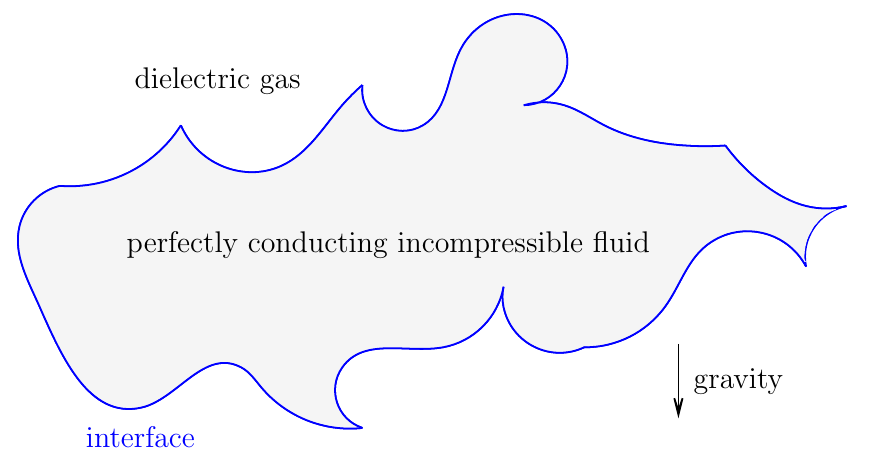}
	\caption{Two-phase EHD flow}
\end{figure} 

EHD, the study of fluid dynamics under electric fields, has wide applications in chemistry, biology and engineering (see \cite{STC}, \cite{T}, \cite{Z1} and also the references cited). Under natural atmospheric conditions, the disintegration of water droplets in a strong electric field plays a significant role in thunderstorm formation. Research indicates that when a pressure difference arises between the interior and exterior of a droplet, hydrodynamic instability causes droplet breakup, forming a conical tip with an obtuse angle at its end. This conical structure is known as the Taylor cone (see \cite{T}). A prominent application is EHD inkjet printing, a high-resolution droplet deposition system. Its operational principle relies on creating an electric potential difference between a conductive nozzle and a substrate. This electric potential gradient induces Maxwell stress tensor-driven interfacial deformation, culminating in the formation of a conical meniscus (Taylor cone) at the nozzle tip. When the voltage is increased, a jet is ejected from the apex of the cone towards the substrate (see Fig. 2). Owing to its practical significance, a deep understanding of EHD interfacial waves is of great importance. Hence, using the EHD equations as a model is of particular interest to us. 
\begin{figure}[!h]
	\includegraphics[width=115mm]{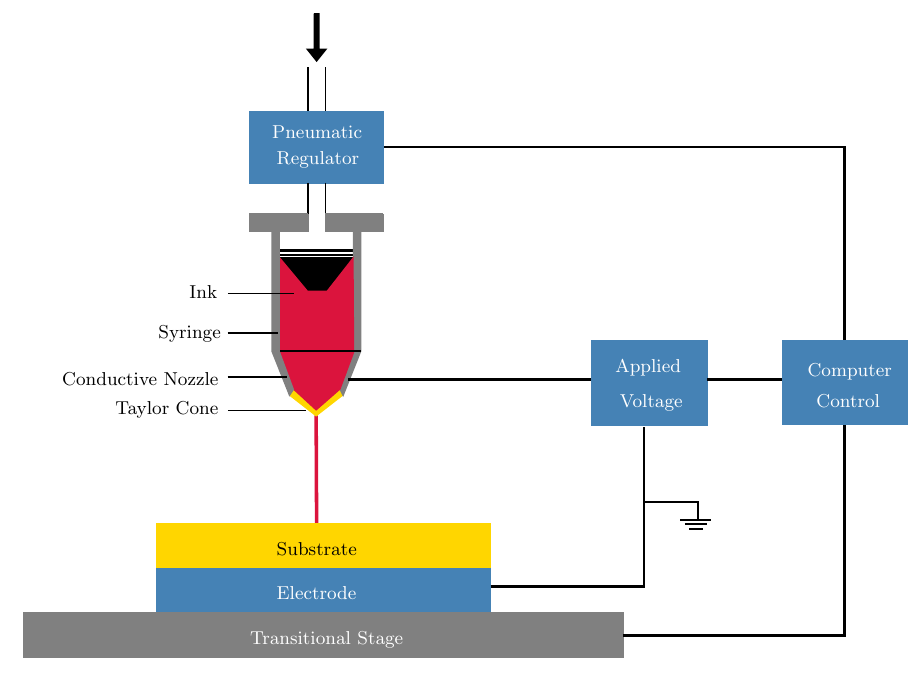}
	\caption{EHD inkjet printing apparatus, as described in \cite{JB}}
\end{figure} 

From a mathematical viewpoint, this problem can be formulated as a two-phase flow free boundary problem (the detailed derivation is provided in the next subsection). A landmark in the mathematical theory of the steady two-phase flow problems was the work of Alt-Caffarelli-Friedman in \cite{ACF1}, who established fundamental regularity properties of the free boundary to the variational problems. Shortly afterwards, Caffarelli introduced the concept of viscosity solution in a series of seminal papers \cite{C1}$-$\cite{C3}, significantly advancing the study of the optimal regularity of both solutions and free boundaries. These developments ultimately established that the free boundary is  $C^{1,\alpha}$ away from stagnation points. Furthermore, by combining the techniques from \cite{D}, De Silva, Ferrari and Salsa established the $C^{1,\alpha}$ regularity of the free boundary of two-phase flow problem governed by inhomogeneous elliptic equations in \cite{DFS1,DFS2} and \cite{DFS5} (for elliptic equations in divergence form, see \cite{DFS3,DFS4}). A significant recent breakthrough in the regularity analysis of free boundaries for two-phase flow was achieved by \cite{PSV}. That work rigorously analyzed the $C^{1,\alpha}$ regularity of the free boundaries near the branching points, where both one-phase and two-phase free boundaries coexist within any arbitrary neighborhood. Similar results were later extended to three-dimensional axisymmetric two-phase flow in \cite{DJ}.

The regularity of free boundaries near non-stagnation points for the EHD equations can be established by similar methods. In contrast, the analysis of singularities near stagnation points has rarely been addressed in the mathematical literature. Our investigation of the EHD equations with gravity is inspired by the pioneer work of Garcia, V\v{a}rv\v{a}ruc\v{a} and Weiss (see \cite{GVW}), who first studied the singularity of the free surface to the EHD equations in the absence of gravity. More precisely, they rigorously proved that at any stagnation point on the axis of symmetry, the axisymmetric free boundary forms a cusp, which is oriented either toward the fluid region or the gas region. However, the singular analysis of the interface between a perfectly conducting flow and a dielectric gas with a gravity field near stagnation points is still an open problem, even in two dimensions. 

On the other hand, the question of singular analysis of the free surface near the stagnation points in the one-phase case has been addressed in a few works, which is related to an important issue known as the  {\textit{Stokes conjecture}}. In 1880, Stokes conjectured that the free boundary of an irrotational water wave forms a symmetric corner with an angle of $\frac{2\pi}{3}$ at each stagnation point in \cite{Stokes} (see Fig. 3). An elegant and rigorous proof of this conjecture was provided by V\v{a}rv\v{a}ruc\v{a} and Weiss using a geometric approach in \cite{VW1}, without imposing any structural assumption on free surface. Specifically, for a stationary incompressible water wave under gravity with stagnation points on the free boundary, they confirmed that the asymptotics of the gravity water wave at any stagnation point is given by the ``Stokes corner flow" in \cite{VW1}. Building on this work, V\v{a}rv\v{a}ruc\v{a} and Weiss revealed in \cite{VW2}, that for the incompressible fluid with vorticity, an interesting new feature (called cusp) can form near the stagnation points, which would not occur in the irrotational case. Shortly afterwards, for a two-dimensional incompressible flow in a gravity field with a free surface subject to surface tension, the results of \cite{WZ} showed that no cusp forms and the free boundary in this water wave problem is analytic. More recently, singularities of the two-phase flow free boundary problem with a gravity field in two dimensions has been investigated in \cite{DJ2}. Regarding the three-dimensional steady axisymmetric gravity water wave problem without vorticity, a pioneering breakthrough in the singular analysis near the stagnation points and along the axis of symmetry was achieved in \cite{VW3}. More precisely, they showed that the only possible non-trivial wave profile near the stagnation points (away from the axis of symmetry) is the Stokes corner. Furthermore, they found that the only possible singular profile along the axis of symmetry (except the origin) is a cusp, while the only non-trivial possible singular profile at the origin is a new feature termed the Garabedian corner. Recently, building upon this work, the effect of vorticity was incorporated in \cite{DHP}. Specifically, the authors provided a classification of the possible singular profiles near the stagnation points and along the axis of symmetry for the steady axisymmetric flow with vorticity. Drawing on the explorations of this work, a comprehensive analysis of both regularity and singularity of the free boundary was completed in \cite{DY}.

\begin{figure}[!h]
	\includegraphics[width=85mm]{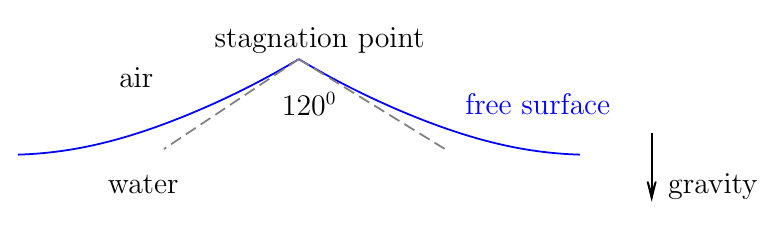}
	\caption{Stokes conjecture for 2D water wave}
\end{figure} 
As previously mentioned, to the best of our knowledge, the pioneering work \cite{GVW} remains the only study addressing singularities of the free interface for the EHD equations. In that work, the cusp singularity at stagnation points was verified through a careful observation and analysis of the transition condition on the free surface. Motivated by significant exploration of this work \cite{GVW}, we will pay attention to the singularities in free boundary for two-dimensional EHD equations under gravity in this paper.

\subsection{Description of the mathematical problem}
\quad

A brief review of Classical Electrodynamics and Fluid Mechanics (see \cite{CM}, \cite{L} and \cite{S}) provides useful background for the EHD model considered here.

For simplicity, we consider a two-dimensional model consisting of a perfectly conducting fluid under gravity, surrounded by a dielectric gas, where the interface is between fluid and gas (see Fig. 1). 

The fluid is assumed to be stationary, irrotational, incompressible, inviscid and perfectly conducting. On the other hand, the gas is assumed to be dielectric, and the electric field $\mathbf{E}$ arises from a potential difference between the perfectly conducting fluid region and the outer boundary. From the conservation principle of charge, the electric field $\mathbf{E}$ can be expressed as the gradient of a scalar electric potential $\Phi$, i.e.,  $\mathbf{E}=-\nabla\Phi$. 
Since the fluid is assumed to be perfectly conducting, the electric potential $\Phi$ within it must be constant, which implies that $\mathbf{E}\equiv0$ in the fluid region. The standard electric boundary condition enforces continuity of the tangential component of $\mathbf{E}$ across the interface. Together with the vanishing of $\mathbf{E}$ on the fluid side, yields
$$
\left(\mathbf{E},\mathbf{\nu}\right)=0 \quad\text{on the interface},
$$
where $\mathbf{\nu}$ is an arbitrary tangent vector to the interface and $\left(\,\,,\,\,\right)$ denotes the Euclidean inner product in $\mathbb{R}^2\times\mathbb{R}^2$. Hence, the interface constitutes an equipotential surface. It follows that the stationary electric potential $\Phi$ takes a single constant value throughout the fluid and on its interface. Without loss of generality, we may normalize this value to zero, namely,
$$\Phi=0,\quad\text{in the fluid and on its free surface.}$$

Let $\rho$ be the constant density of the incompressible fluid. Without losing generality, we set $\rho\equiv1$. Thus the steady incompressible inviscid flow is governed by the following two-dimensional Euler system
\begin{equation}
	\label{1.1}
	\begin{cases}
	\nabla\cdot \mathbf{u}=0,\\
	\left(\mathbf{u}\cdot\nabla\right)\mathbf{u}+\nabla p+g\mathbf{e}_2=0,\\
	\end{cases}
\end{equation}
with irrotational condition
\begin{equation}
\label{1.2}
\nabla\times \mathbf{u}=0,
\end{equation}
where $\mathbf{u}=(u_1,u_2)$ denotes the velocity field, $g$ is the gravitational constant, $p$ is the pressure, $x=(x_1,x_2)$ denotes the space variable and $\mathbf{e}_2=(0,1)$ (see Fig. 4). 
\begin{figure}[!h]
	\includegraphics[width=90mm]{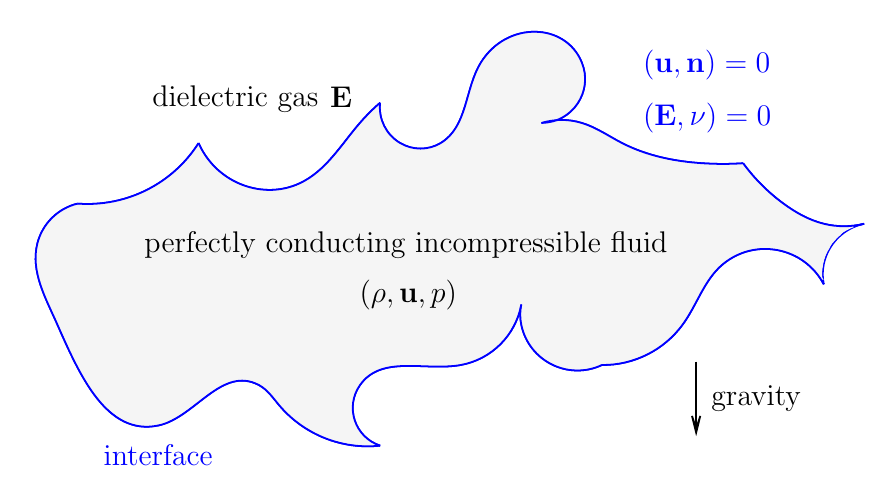}
	\caption{The physical model of two-phase EHD flow}
\end{figure} 

From the conservation of mass in \eqref{1.1}, a stream function $\Psi$ can be introduced such that
\begin{equation}
\nonumber
u_1=\frac{\partial\Psi}{\partial x_2}\quad \text{and}\quad u_2=-\frac{\partial\Psi}{\partial x_1}.
\end{equation}
It then follows from \eqref{1.2} that
\begin{equation}
\label{1.5}
\Delta\Psi=0\quad\text{in\,\,the\,\,fluid\,\,region,}
\end{equation}

On the other hand, the charge conservation implies that the electric potential $\Phi$ satisfies
\begin{equation}
\label{1.7}
\Delta\Phi=0\quad\text{in\,\,the\,\,gas\,\,region.}
\end{equation}

On the other hand, Bernoulli's law for the steady, incompressible and irrotational fluid yields
\begin{equation}
\label{B}
\frac{1}{2}\left|\mathbf{u}\right|^2+gx_2+p=\mathcal{B}\quad\text{in the incompressible fluid field},
\end{equation}
where $\mathcal{B}$ is the Bernoulli constant for the incompressible fluid.

Notice that $-p\,\mathbf{n}=p_e \,\mathbf{n}\quad\text{on the interface}$, where $\mathbf{n}$ is the unit normal vector to the interface, $p_e=\frac{\varepsilon}{2}\,\mathbf{E}^2$ is the electric energy density and $\varepsilon>0$ is the dielectric permittivity (refer to page 136 in \cite{GVPMS} and page 539 in \cite{S}). Using the relation $p=-p_e$ in Bernoulli's law \eqref{B}, the transition condition on the free boundary satisfies 
\begin{equation}
\label{New}
\frac{1}{2}\left|\mathbf{u}\right|^2+gx_2-p_e=\mathcal{B}\quad\text{on the interface}.
\end{equation}
Rewriting the velocity term as $\left|\mathbf{u}\right|^2=\left|\nabla\Psi\right|^2$ and substituting $p_e=\frac{\varepsilon}{2}\left|\nabla\Phi\right|^2$ into \eqref{New}, we obtain
\begin{equation}
\label{1.6}
\frac{1}{2}\left|\nabla\Psi\right|^2-\frac{\varepsilon}{2}\left|\nabla \Phi\right|^2=\mathcal{B}-gx_2\quad\text{on the interface}.
\end{equation}

 Furthermore, by introducing the rescaled functions
 \begin{equation}
\nonumber
 \tilde\Psi=\sqrt{\frac{1}{2g}}\Psi,\quad\tilde{\Phi}=\sqrt{\frac{\varepsilon}{2g}}\Phi,\quad x_2^0=\frac{\mathcal{B}}{g},
 \end{equation}
 we can rewrite the transition condition \eqref{1.6} in the normalized form,
 \begin{equation}
 	\label{1.8}
 	\left|\,\nabla\tilde\Psi\,\right|^2-\left|\,\nabla\tilde\Phi\,\right|^2=x_2^0-x_2,
 \end{equation}
 where $x_2^0=\mathcal{B}/g$ is called {\textit{uniform height}} of the incompressible fluid.
 
 For simplicity, we will continue to denote $\tilde{\Psi}$ and $\tilde{\Phi}$ by $\Psi$ and $\Phi$, respectively, throughout the end of this paper. Therefore, combining \eqref{1.5}, \eqref{1.7} and \eqref{1.8}, the two-phase free boundary problem for the EHD equations can be simplified as
  \begin{equation}
 \label{EHD}
 \begin{cases}
 \begin{matrix}
 \Delta\Phi=0,&	\hspace{31.7mm} \text{in the gas region},\\
 \end{matrix}\\
 \begin{matrix}
 \Delta\Psi=0,& \hspace{31.7mm} \text{in the fluid region},\\
 \end{matrix}\\
  \begin{matrix}
\left|\,\nabla\Psi\,\right|^2-\left|\,\nabla\Phi\,\right|^2=x_2^0-x_2,& \text{on the interface},\\
 \end{matrix}\\
   \begin{matrix}
\Phi=0,& \hspace{35.1mm} \text{on the interface},\\
 \end{matrix}\\
 \begin{matrix}
 \Psi=0,& \hspace{35.1mm} \text{on the interface}.\\
 \end{matrix}\\
 \end{cases}
  \end{equation}
  
Since both $\Phi$ and $\Psi$ are harmonic in their respective domains (from \eqref{EHD}),  the minimum principle implies $\Phi>0$ in the gas region (see \cite{GT}). Similarly, $\Psi>0$ in the fluid region. Furthermore, we can define the Stokes stream function $\Psi\equiv0$ in gas region, while the electric potential $\Phi\equiv0$ in the perfectly conducting fluid region. This allows us to define a continuous function $u$ across the domains as follows,
\begin{equation}
	\label{u}
	u=\begin{cases}
		\begin{matrix}
			\Phi,&		\hspace{2.75mm} \text{in the gas region,}
		\end{matrix}\\
		\begin{matrix}
			-\Psi,& \text{in the fluid region.}
		\end{matrix}\\
	\end{cases}
\end{equation}
Hence, the positive and negative parts of $u$ satisfy
$$u^+=\max\{u,0\}=\Phi~~\text{in the gas region},$$
$$ u^-=\min\{u,0\}=-\Psi~~\text{in the fluid region},$$
 and 
 $$u=0\quad\text{on the interface}.$$
Furthermore, we identify
$$\{u>0\}\,\,\text{as the gas region,}\quad \{u<0\}\,\,\text{as the fluid region,}$$ 
 and
  $$\partial\{u>0\}\cap\partial\{u<0\}\,\,\text{ as the interface between fluid and gas.}$$
   Consequently, the two-phase free boundary problem \eqref{EHD} can be reformulated in terms of $u$ as follows (see Fig. 5),
\begin{equation}
\label{1.9}
\begin{cases}
\begin{matrix}
\Delta u=0,&	\hspace{33.4mm} 	 \text{in 
}\,\,\{u\ne0\},\\
\end{matrix}\\
\begin{matrix}
\left|\nabla u^-\right|^2-\left|\nabla  u^+\right|^2=x_2^0-x_2,& \text{on}\,\,\partial \{u>0\}\cap\partial \{u<0\}.\\
\end{matrix}\\
\end{cases}
\end{equation}
\begin{figure}[!h]
	\includegraphics[width=115mm]{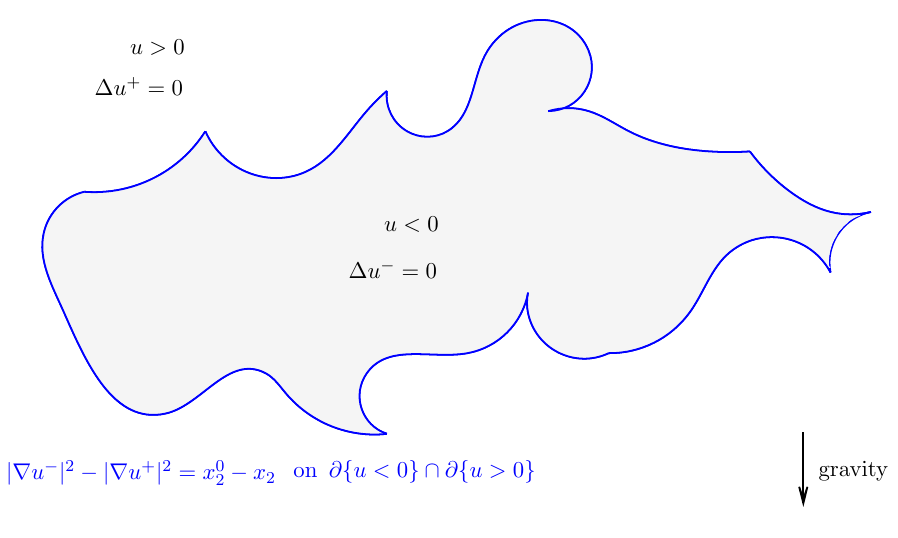}
	\caption{The two-phase free boundary problem}
\end{figure} 

In this paper, we focus on the singular profiles of the interface near the stagnation points of the incompressible fluid. A stagnation point $x^{\text{\rm st}}=\left(x^{\text{\rm st}}_1,x^{\text{\rm st}}_2\right)$ on the interface is defined by the condition $\left|\nabla u^-\left(x^{\text{\rm st}}\right)\right|=0$. From the transition condition in \eqref{1.9}, we observe that for $x_2<x_2^0=\mathcal{B}/g$, the right-hand side of the transition condition is strictly positive, which implies $\left|\nabla u^-(x)\right|>0$ as $x_2<x_2^0$. By classical regularity theory in \cite{PSV}, it follows that the free boundary $\partial\{u<0\}$ is smooth. Therefore, the main focus of this work is on the possible singularities of the interface for points with $x_2\geqslant x_2^0$. 

One of the key difficulties in this problem concerns the singular profiles near the stagnation points of the incompressible fluid, specifically regarding whether the gradient of the positive-phase $\left|\nabla u^+\right|$ is degenerate or non-degenerate. This distinction leads to fundamentally different behaviors of the positive-phase free surface, which is smooth in the non-degenerate case but singular in the degenerate case. We analyze both cases in Subsection 1.4.

The transition condition in \eqref{1.9} provides some essential and useful information, such as the possible height of the stagnation points and the decay rate of $\left|\nabla u^-\right|$ near them. Furthermore, the right-hand side $x_2^0-x_2$ of the transition condition, which originated from the gravity force, results in singular wave profiles distinct from those without gravity. For instance, in the axisymmetric EHD model without gravity studied in \cite{GVW}, the transition condition on the free surface becomes
$$
x_1\left|\nabla u^+(x_1,x_2)\right|^2-\frac{1}{x_1}\left|\nabla u^-(x_1,x_2)\right|^2=0,
$$
leading to cusp singularities, where $x_2$ is the symmetric axis. This contrast thus motivates our present study of the free boundary problem with gravity. A more detailed discussion will be provided in the next subsection.

\subsection{A variational approach}
\quad

To investigate the two-phase free boundary problem, the Alt-Caffarelli-Friedman (ACF) functional
$$
E_{\text{ACF}}(u,\Omega)=\int_{\Omega}\left(\left|\nabla u\right|^2+\lambda_+^2\chi_{\{u>0\}}+\lambda_-^2\chi_{\{u<0\}}+\lambda_0^2\chi_{\{u=0\}}\right)dx
$$
has been a central object of study since its introduction in the seminal paper \cite{ACF1}, inspiring substantial research in two-phase flows. Here $\lambda_+$, $\lambda_-$, $\lambda_0$ are non-negative constants and $\chi_D$ is the characteristic function of the set $D$. 

The associated variational problem consists of minimizing the functional $E_{\text{ACF}}$ over the convex set $\mathcal{K}$
$$
\mathcal{K}=\left\{u\in W_{\text{loc}}^{1,2}(\Omega),\,\,u=u^0\,\,\text{on}\,\,S\right\},
$$
where $S$ is a non-empty open subset of $\partial\Omega$ and $u^0\in W_{\text{loc}}^{1,2}(\Omega)$ is a given boundary function.

Minimizers are harmonic away from the zero level set $\{u=0\}$, where the jump conditions on the gradient of $u$ arise due to the discontinuous characteristic functions in the energy. The Euler-Lagrange equation is
 \begin{equation}
 	\label{ACF}
 	\Delta u=0\quad \text{in
 	}\,\,\{u\ne0\},
 \end{equation}
 and the following free boundary conditions hold
\begin{equation}
\label{ACF FB}
\begin{cases}
\begin{matrix}
\left|\nabla u^-\right|^2-\left|\nabla  u^+\right|^2=\lambda_-^2-\lambda_+^2,& \text{on 
}\,\,\Gamma_{\text{\rm tp}}=\partial\{u>0\}\cap\partial\{u<0\},\\
\end{matrix}\\
\begin{matrix}
\left|\nabla u^-\right|^2=\lambda_-^2-\lambda_0^2,&	\hspace{17.1mm} 	 \text{on 
}\,\,\Gamma_{\text{\rm op}}^-=\partial\{u<0\}\backslash\partial\{u>0\},\\
\end{matrix}\\
\begin{matrix}
\left|\nabla u^+\right|^2=\lambda_+^2-\lambda_0^2,&	\hspace{17.1mm} 	 \text{on 
}\,\,\Gamma_{\text{\rm op}}^+=\partial\{u>0\}\backslash\partial\{u<0\},\\
\end{matrix}\\
\end{cases}
\end{equation}
where $\Gamma_{\text{\rm tp}}$ denotes the set of two-phase free boundary points, while $\Gamma_{\text{op}}^\pm$ denotes the set of one-phase free boundary points. These sets are defined as follows,

(1) A point $x\in\Gamma_{\text{\rm tp}}$ if the Lebesgue measures satisfy 
$$\left|B_r(x)\cap\{u>0\}\right|>0\quad\text{and}\quad \left|B_r(x)\cap\{u<0\}\right|>0,\quad\text{for any}\,\,  r>0.$$

(2) A point $x\in\Gamma_{\text{\rm op}}^+$ if there exists a $r_0$ small enough such that $u\geqslant0$ in $B_{r_0}(x)$.

(3) A point $x\in\Gamma_{\text{\rm op}}^-$ if there exists a $r_0$ small enough such that $u\leqslant0$ in $B_{r_0}(x)$.

 Under the condition $\lambda_0=\min\left\{\lambda_+,\lambda_-\right\}$, Alt, Caffarelli and Friedman established the Lipschitz continuity of solutions and the $C^1$ regularity of free boundary for the problem \eqref{ACF} and \eqref{ACF FB} in \cite{ACF1}. Moreover, the $C^{1,\alpha}$ regularity of the free boundaries for this two-phase flow was established in \cite{C1}$-$\cite{C3} using the theory of viscosity solutions. More recently, De Silva, Ferrari and Salsa verified the $C^{1,\alpha}$ regularity of free boundary for two-phase flow problem governed by elliptic operators with forcing terms in \cite{DFS1,DFS2} and \cite{DFS5}.

An interesting feature of the two-phase free boundary $\Gamma_{\text{\rm tp}}$ is the presence of so-called {\it branching points} denoted as $x^{\text{\rm br}}$, where 
$$x^{\text{\rm br}}\in\Gamma_{\text{\rm tp}}\quad\text{and}\quad\text{the Lebesgue measure} \,\,\left|B_r(x^{\text{br}})\cap\{u=0\}\right|>0,\,\,\text{for any}\,\,r>0.$$
 The analysis of the free boundary near the branching points is still a long-standing and open problem. However, it is known that the functional $E_{\text{ACF}}$ does not possess any branching points when $\lambda_0\geqslant\min\left\{\lambda_+,\lambda_-\right\}$, as detailed in $\text{Section 6}$ of \cite{ACF1}. This naturally leads us to investigate whether branching points may arise when $\lambda_0<\min\left\{\lambda_+,\lambda_-\right\}$. Notably, an ingenious construction in \cite{DEST} demonstrated that there exists a family of minimizers to the functional $E_{\text{ACF}}$ with $\lambda_+^2=\lambda_-^2=1$ and $\lambda_0=0$ for which branching points necessarily exist. More recently, a major breakthrough was achieved by De Philippis, Spolaor and Velichkov (see \cite{PSV}) in the regularity theory for two-phase flow problem. To study the free boundaries near the branching points, they considered the two-phase functional
$$
E_{\text{tp}}(u,\Omega)=\int_{\Omega}\left(\left|\nabla u\right|^2+\lambda_+^2\chi_{\{u>0\}}+\lambda_-^2\chi_{\{u<0\}}\right)dx,
$$
and proved that in a neighborhood of every two-phase point, the free boundary is $C^{1,\alpha}$-regular curve provided $\lambda_+>0$ and $\lambda_->0$. Meanwhile, the regularity of the one-phase free boundaries $\Gamma_{\text{\rm op}}^+$ and $\Gamma_{\text{\rm op}}^-$ follows from the results in \cite{W}. 

In the context of this work, however, the parameters $\lambda_+$, $\lambda_-$ and $\lambda_0$ depend on the $x_2$-coordinate and degenerate at $x_2=x_2^0$. Consequently, we select them appropriately to construct a suitable EHD functional. Comparing the free boundary conditions in \eqref{ACF FB} with the transition condition 
\begin{equation}
\nonumber
	\left|\nabla u^-\left(x\right)\right|^2-\left|\nabla u^+\left(x\right)\right|^2=x_2^0-x_2\quad\text{on}\,\,\Gamma_{\text{\rm tp}}.
\end{equation}
	It is easy to see that the parameters must satisfy $\lambda_-^2-\lambda_+^2=x_2^0-x_2$. This motivates the following choice (without loss of generality),
\begin{equation}
	\label{parameter}
\lambda_+^2=\left(x_2-x_2^0\right)^+,\quad\lambda_-^2=\left(x_2-x_2^0\right)^++x_2^0-x_2,\quad\lambda_0^2=0,
\end{equation}
with $x_2^0=\frac{\mathcal{B}}{g}$.

Throughout this paper, we introduce the following EHD functional
$$
E_{\text{EHD}}(u,\Omega)=\int_{\Omega}\left(\left|\nabla u\right|^2+\left(x_2-x_2^0\right)^+\chi_{\{u>0\}}+\left(\left(x_2-x_2^0\right)^++x_2^0-x_2\right)\chi_{\{u<0\}}\right)dx.
$$
The corresponding free boundary problem satisfies (see Appendix A for a detailed derivation)
\begin{equation}
\label{1.20}
\begin{cases}
\begin{matrix}
\Delta u=0,&\hspace{40mm}  \text{in 
}\,\,\{u\ne0\},\\
\end{matrix}\\
\begin{matrix}
\left|\nabla u^-\right|^2-\left|\nabla  u^+\right|^2=x_2^0-x_2,&\hspace{6.9mm}  \text{on 
}\,\,\Gamma_{\text{\rm tp}}=\partial\{u>0\}\cap\partial\{u<0\},\\
\end{matrix}\\
\begin{matrix}
\left|\nabla u^-\right|^2=\left(x_2-x_2^0\right)^++x_2^0-x_2,&		 \text{on 
}\,\,\Gamma_{\text{\rm op}}^-=\partial\{u<0\}\backslash\partial\{u>0\},\\
\end{matrix}\\
\begin{matrix}
\left|\nabla u^+\right|^2=\left(x_2-x_2^0\right)^+,&	\hspace{16.7mm} 	 \text{on 
}\,\,\Gamma_{\text{\rm op}}^+=\partial\{u>0\}\backslash\partial\{u<0\}.\\
\end{matrix}\\
\end{cases}
\end{equation}

A function $u$ is called a {\textit{local minimizer} of $E_{\text{EHD}}(u,\Omega)$ in $\Omega$ if
$$
E_{\text{EHD}}(u,D)\leqslant E_{\text{EHD}}(v,D),
$$
for every open set $D$ and every function $v\in W^{1,2}(D)$ satisfying $\bar{D}\subset\Omega$ and $v=u$ on $\Omega\backslash D$.

Since the stagnation points are defined by the fluid region $\{u<0\}$, our analysis will focus on the free boundary $\partial\{u<0\}$. The following  proposition provides a key characterization. 
\begin{proposition} 
	Let $u$ be the local minimizer of $E_{\rm EHD}(u,\Omega)$. Then
	$$\Gamma_{\text{\rm op}}^+\cap\{x_2<x^0_2\}=\varnothing\quad\text{and}\quad\Gamma_{\text{\rm op}}^-\cap\{x_2>x^0_2\}=\varnothing.$$
\end{proposition}
\begin{pf}
We argue by contradiction. Suppose that there exists a point $\tilde{x}\in\Gamma_{\text{\rm op}}^+\cap\{x_2<x_2^0\}$. Then, for some small enough $r_0$,   
$$
u(x)\geqslant0\quad\text{in}\,\,B_r(\tilde{x}),
$$
as $0<r<r_0$. Hence, in the region $\Omega\cap\{x_2<x_2^0\}$, the functional simplifies to
$$
E_{\text{EHD}}(u,\Omega\cap\{x_2<x_2^0\})=\int_{\Omega\cap\{x_2<x_2^0\}}\left(\left|\nabla u\right|^2+\left(x_2^0-x_2\right)\chi_{\{u<0\}}\right)dx.
$$
We now claim that $\Delta u=0$ in $B_r(\tilde{x})$.
If not, there exists a $w(x)$ solving the Dirichlet problem
$$
\begin{cases}
\begin{matrix}
\Delta w=0,&\text{in}\,\,B_r(\tilde{x}),		\\
\end{matrix}\\
\begin{matrix}
w=u,&\hspace{2.3mm}\text{on}\,\,\partial B_r(\tilde{x}),		\\
\end{matrix}\\
\end{cases}
$$
such that 
\begin{equation}
\label{1.19}
\int_{B_r(\tilde{x})}\left|\nabla w\right|^2dx<\int_{B_r(\tilde{x})}\left|\nabla u\right|^2dx.
\end{equation}
It follows from $u\geqslant0$ that 
$$
\begin{cases}
\begin{matrix}
\Delta \left(w-u\right)\leqslant0,&\text{in}\,\,B_r(\tilde{x}),		\\
\end{matrix}\\
\begin{matrix}
w=u,&\hspace{12.6mm}\text{on}\,\,\partial B_r(\tilde{x}).	\\
\end{matrix}\\
\end{cases}
$$
Thus, the minimum principle implies $w\geqslant u\geqslant0$ in $B_r(\tilde{x})$. Together with \eqref{1.19} and the arbitrariness of $\tilde{x}$, this yields
$$
\int_{\Omega\cap\{x_2<x_2^0\}}\left(\left|\nabla w\right|^2+\left(x_2^0-x_2\right)\chi_{\{w<0\}}\right)dx<\int_{\Omega\cap\{x_2<x_2^0\}}\left(\left|\nabla u\right|^2+\left(x_2^0-x_2\right)\chi_{\{u<0\}}\right)dx,
$$
which contradicts to the definition of local minimizer.

Hence, $\Delta u=0$ in $B_r(\tilde{x})$. Combining this with the maximum principle gives $u\equiv0$ in $B_r(\tilde{x})$, which contradicts the assumption that $\tilde{x}\in\Gamma_{\text{\rm op}}^+$.

As a consequence, $\Gamma_{\text{\rm op}}^+\cap\{x_2<x_2^0\}=\varnothing$. The second statement $\Gamma_{\text{\rm op}}^-\cap\{x_2>x_2^0\}=\varnothing$ follows by a similar argument.
\end{pf}
\begin{remark}
	This property remains valid for weak solutions  (as defined by \text{Definition A} in Appendix A). Suppose, for contradiction, that there exist a point $\tilde{x}\in\Gamma_{\text{\rm op}}^+\cap\{x_2<x_2^0\}$ and a sufficiently small $r>0$ such that $u\geqslant0$ in $B_r(\tilde{x})$. From the definition of weak solution and integration by parts, we obtain
	$$
	0=\int_{ B_{r}(\tilde{x})}\left(\left|\nabla u^+\right|^2\mathrm{div}\tilde{\phi}-2\nabla u^+ D\tilde{\phi}\nabla u^+\right)dx,
	$$
where $\tilde{\phi}=\phi\cap B_{r}(\tilde{x})$ and $\phi\in C^1_0(\Omega,\mathbb{R}^2)$. By $\text{Theorem 3.1}$ in \cite{CF}, for any $\varepsilon>0$, one has $$\mathcal{H}^{\varepsilon}\left(\left\{
u^+=0,\,\left|\nabla u^+\right|=0\right\}\cap B_1\right)=0,$$ 
where $\mathcal{H}^{\varepsilon}$ is the $\varepsilon$-dimensional Hausdorff measure. It follows that $\Delta u^+=0$ in $B_{r}(\tilde{x})$. Hence, the maximum principle implies that $u^+\equiv0$. This contradicting the assumption that $\tilde{x}\in\Gamma_{\text{\rm op}}^+$.
\end{remark}

$\text{Proposition 1.1}$ implies that the one-phase free boundary $\Gamma_{\text{\rm op}}^+$ must lie above the uniform height $x_2^0=\frac{\mathcal{B}}{g}$, while $\Gamma_{\text{\rm op}}^-$ lies entirely below $x_2^0$. Hence,
\begin{equation}
	\label{1.27}
	\partial\{u<0\}\cap\{x_2<x_2^0\}=\left(\Gamma_{\text{\rm tp}}\cup\Gamma_{\text{\rm op}}^-\right)\cap\{x_2<x_2^0\}
\end{equation}
and
\begin{equation}
	\label{1.28}
	\partial\{u<0\}\cap\{x_2>x_2^0\}=\Gamma_{\text{\rm tp}}\cap\{x_2>x_2^0\}.
\end{equation}

The key advantage of the parameter choice in \eqref{parameter} lies in yielding the identity
$$
\lambda_0^2(x_2)=\min\left\{\left(x_2-x_2^0\right)^+,\left(x_2-x_2^0\right)^++x_2^0-x_2\right\},
$$
and in ensuring that the condition on $\Gamma_{\text{\rm tp}}$ matches the physical transition condition
\begin{equation}
	\nonumber
	\left|\nabla u^-\left(x\right)\right|^2-\left|\nabla u^+\left(x\right)\right|^2=x_2^0-x_2\quad\text{on}\,\,\Gamma_{\text{\rm tp}}.
\end{equation}
However, there are additional conditions in \eqref{1.20} 
\begin{equation}
\label{FB}
\begin{cases}
\begin{matrix}
	\left|\nabla u^-\right|^2=\left(x_2-x_2^0\right)^++x_2^0-x_2,&		 \text{on 
	}\,\,\Gamma_{\text{\rm op}}^-,\\
\end{matrix}\\
\begin{matrix}
	\left|\nabla u^+\right|^2=\left(x_2-x_2^0\right)^+,&	\hspace{16.7mm} 	 \text{on 
	}\,\,\Gamma_{\text{\rm op}}^+.\\
\end{matrix}\\
\end{cases}
\end{equation}
Critically, our analysis remains unaffected by the additional conditions $\eqref{FB}$. Indeed, the free boundary $\partial\{u<0\}$ depends solely on $\Gamma_{\text{\rm tp}}$ and $\Gamma_{\text{\rm op}}^-$, independent of $\Gamma_{\text{\rm op}}^+$. Moreover, for $x^{\text{\rm st}}\in\Gamma_{\text{\rm op}}^-$, the possible singular profile of $\Gamma_{\text{\rm op}}^-$ resembles that in one-phase case (as shown in \cite{VW1}). This work therefore focuses exclusively on the analysis of singularities for the two-phase free boundary $\Gamma_{\text{\rm tp}}$.
 
\subsection{The analysis of the transition condition}
\quad

The main purpose of this paper is to investigate the singular profiles of the interface near the possible stagnation points.

A pioneering analysis of singularities in the free boundary for axisymmetric EHD equations was conducted by Garcia, V\v{a}rv\v{a}ruc\v{a} and Weiss in \cite{GVW}. There, the authors analyzed a two-phase axisymmetric free boundary problem with gravity  neglected
\begin{equation}
	\label{GVW}
	\begin{cases}
		\begin{matrix}
			\mathrm{div}\left(x_1\nabla u(x)\right)=0,&	\hspace{14.2mm} 	 \text{in gas region},\\
		\end{matrix}\\
		\begin{matrix}
			\mathrm{div}\left(\frac{1}{x_1}\nabla u(x)\right)=0,&	\hspace{12.7mm} 	 \text{in fluid region},\\
		\end{matrix}\\
		\begin{matrix}
			x_1\left|\nabla u^+\right|^2-\frac{1}{x_1}\left|\nabla  u^-\right|^2=0,& \text{on the interface between fluid and gas,}\\
		\end{matrix}\\
	\end{cases}
\end{equation}
where $x=(x_1,x_2)$, $x_2$ is the symmetric axis and $\Omega$ is a connected open subset of the right half-plane $\{x\left|\,x_1\geqslant0\right.\}$. Note that the transition condition in \eqref{GVW} implies 
$$\left|\nabla u^-(x)\right|^2=x_1^2\left|\nabla u^+(x)\right|^2\quad\text{on the interface between fluid and gas.}$$
Consequently, points on the $x_2$-axis must be stagnation points, where the gradient of the fluid-phase vanishes, i.e., $\left|\nabla u^-\right|=0$. Specifically, under certain non-degeneracy condition and assuming an injective free surface, it was shown that the only possible singular profile of the free boundary near the stagnation points is cusp.  

However, when gravity effect is considered, a celebrated conjecture (the so-called Stokes conjecture) implies the corner type will occur in the one-phase case. This naturally raises the question of whether a similar structure applies to singular profiles in the two-phase flow problem. In particular, we highlight that the transition condition \eqref{1.9} provides some useful and important information for analyzing singularities on the free boundary, as it indicates that the decay rate of $\left|\nabla u^-(x)\right|$ depends on the positive-phase term $\left|\nabla u^+(x)\right|$ and the linear term $x_2^0-x_2$. A central finding of our work, guided by the transition condition in \eqref{1.9}, is the identification of a critical decay rate that dictates the classification of all possible singularities. 

More precisely, we define the set of stagnation points of the incompressible fluid by 
\begin{equation}
\nonumber
	S^u=\left\{x\left|\,x\in\partial\{u<0\},\,\,\left|\nabla u^-(x)\right|=0\right.\right\}
\end{equation}
and consider an arbitrary $x^{\text{\rm st}}=\left(x_1^{\text{st}},x_2^{\text{\rm st}}\right)\in S^u$. At any such point $x^{\text{\rm st}}$, the transition condition in \eqref{1.9}
\begin{equation}
	\label{FB1}
	\left|\nabla u^-\left(x\right)\right|^2=\left|\nabla u^+\left(x\right)\right|^2+x_2^0-x_2\quad\text{on the interface,}
\end{equation}
together with the condition $\left|\nabla u^-\left(x^{\text{\rm st}}\right)\right|=0$ implies that $x^{\text{\rm st}}_2\geqslant x_2^0$. Consequently, the height of any stagnation point is at least the uniform height $x_2^0=\frac{\mathcal{B}}{g}$.

This observation leads to a natural classification of stagnation points based on their height relative to $x_2^0$.
$$
\text{\bf Case 1.}\quad x_2^{\text{\rm st}}=x_2^0=\frac{\mathcal{B}}{g},\quad\quad\text{\bf Case 2.}\quad x_2^{\text{\rm st}}>x_2^0=\frac{\mathcal{B}}{g}.
$$
A schematic summarizing all possible scenarios is shown in Fig. 6.

\begin{figure}[!h]
	\includegraphics[width=110mm]{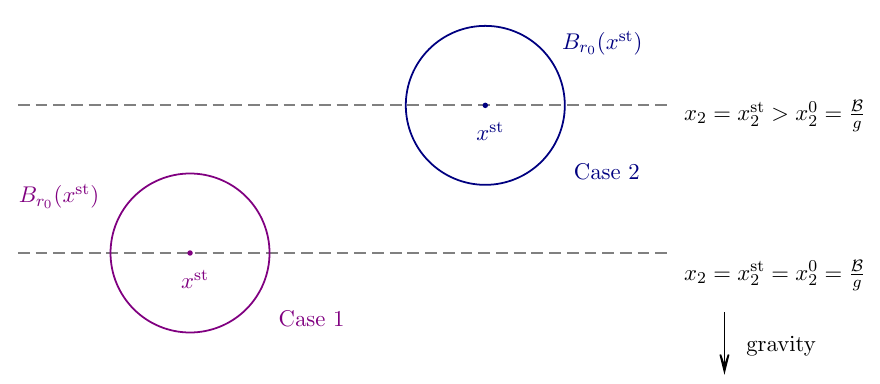}
	\caption{Possible stagnation points}
\end{figure} 

{\bf Case 1.} $x_2^{\text{\rm st}}=x_2^0=\frac{\mathcal{B}}{g}$. In this case, the definition of $u$ implies that the electric field $\left|\mathbf{E}\right|=\left|\nabla u^+\right|$ vanishes at $x^0$. A rigorous analysis of the transition condition \eqref{FB1} further reveals a characteristic decay rate proportional to the square root of the distance, i.e., $\left|x-x^{\text{\rm st}}\right|^{1/2}$. This decay rate plays a crucial role in our subsequent classification of the influence of the electric field on the singular analysis of incompressible fluid near stagnation points. Therefore, we refer to $\left|x-x^{\text{\rm st}}\right|^{1/2}$ as {\textit{critical decay rate}} and divide the Case 1 into the following three subcases.

{\bf Case 1.1.}  $x^{\text{\rm st}}_2=x_2^0$ and $\left|\mathbf{E}\right|$ decays faster than the critical decay rate $\left|x-x^{\text{\rm st}}\right|^{1/2}$ near $x^{\text{\rm st}}$. More precisely, for any $\varepsilon>0$ small enough, there exists $r_0>0$ such that
\begin{equation}
	\nonumber
	\left|\nabla u^+\right|\leqslant \varepsilon\left|x-x^{\text{\rm st}}\right|^{1/2},\quad \text{in}\,\,B_{r_0}(x^{\text{\rm st}})\cap\{u>0\}.
\end{equation}

{\bf Case 1.2.}  $x^{\text{\rm st}}_2=x_2^0$ and $\left|\mathbf{E}\right|$ decays as $\left|x-x^{\text{\rm st}}\right|^{1/2}$ near $x^{\text{\rm st}}$, namely,  there exist $0<C_1\leqslant C_2$ and $r_0>0$ such that
\begin{equation}
	\nonumber
C_1\left|x-x^{\text{\rm st}}\right|^{1/2}\leqslant	\left|\nabla u^+\right|\leqslant C_2\left|x-x^{\text{\rm st}}\right|^{1/2},\quad \text{in}\,\,B_{r_0}(x^{\text{\rm st}})\cap\{u>0\}.
\end{equation}

{\bf Case 1.3.}  $x^{\text{\rm st}}_2=x_2^0$ and $\left|\mathbf{E}\right|$ decays slower than $\left|x-x^{\text{\rm st}}\right|^{1/2}$ near $x^{\text{\rm st}}$, namely, for any $M>0$ large enough, there exists $r_0>0$ such that
\begin{equation}
	\nonumber
	\left|\nabla u^+\right|\geqslant M\left|x-x^{\text{\rm st}}\right|^{1/2},\quad \text{in}\,\,B_{r_0}(x^{\text{\rm st}})\cap\{u>0\}.
\end{equation}

This classification exhausts the possible asymptotic behaviors of $\left|\nabla u^+\right|$ near a stagnation point in Case 1. 

{\bf Case 2.} $x_2^{\text{\rm st}}>x_2^0=\frac{\mathcal{B}}{g}$. Suppose a stagnation point lies above the uniform height. At such a point, the transition condition in \eqref{1.9} and $\left|\nabla u^-\left(x^{\text{\rm st}}\right)\right|=0$ implies that $\left|\nabla u^+\left(x^{\text{\rm st}}\right)\right|^2=x_2^{\text{\rm st}}-x_2^0>0$. This  non-degeneracy of the positive phase, via the method of improvement of flatness (see \cite{PSV}), would force the two-phase free boundary $\Gamma_{\text{\rm tp}}=\partial\{u>0\}\cap\partial\{u<0\}$ to be a $C^{1,\beta}$ curve in a small neighborhood of $x^{\text{\rm st}}$ for some $\beta\in(0,1)$. A detailed proof of this result is provided in Appendix B for completeness.  

However, this leads to a contradiction. Recall from $\text{Proposition 1.1}$ and \eqref{1.28} that in the region $\{x_2>x_2^0\}$, the free boundary $\partial\{u<0\}$ consists solely of $\Gamma_{\text{\rm tp}}$ (since $\Gamma_{\text{\rm op}}^-$ is empty there). Consequently, the free boundary is entirely $C^{1,\beta}$ near $x^{\text{\rm st}}$, precluding the existence of a stagnation point. Therefore, Case 2 cannot occur.

In summary, the transition condition \eqref{1.9} forces any stagnation point to lie exactly at the uniform height $x_2=x_2^0$ (Case 1), where the electric field degenerates. Taking $x^{\text{\rm st}}=x^0=(x_1^0,x_2^0)$, the subsequent analysis therefore pivots on the detailed classification of subcases within Case 1, which is governed by the decay rate of $\left|\nabla u^+\right|$. 
\subsection{Main results}
\quad

This paper aims to characterize the singularity of the free boundary near the stagnation points $x^{\text{\rm st}}=x^0=(x_1^0,x_2^0)$ on $\partial\{u<0\}$. 

\begin{theorem 1}
	  Let $u$ be a weak solution of problem \eqref{1.20} (as defined by \text{Definition A} in Appendix A). For $r_0>0$ small enough, we assume the following statements hold
	 \begin{itemize}
	 	\item   $u\geqslant0$ in $\{x_2\geqslant x_2^0\}$.
	 	\item  For Case 1.1 and Case 1.2, we
	 	suppose that $\{u=0\}$ has locally only finitely many connected components and
	 	$$\left|\nabla u^-\right|\leqslant\left|x_2-x_2^0\right|^{1/2}\quad \text{in} \,\,B_{r_0}(x^0),$$
	  while for Case 1.3, we assume that there exists $\alpha\in\left(1,\frac{3}{2}\right)$ such that for some $C>0$
	   $$\left|\nabla u^-\right|\leqslant C\left|x-x^0\right|^{\alpha-1}\quad\text{and}\quad\left|\nabla u^+\right|\leqslant C\left|x-x^0\right|^{\alpha-1}\quad \text{in} \,\,B_{r_0}(x^0).$$
	 \end{itemize}	 	
	 Then the following asymptotic behaviors hold
	 \begin{itemize}
	 	\item For Case 1.1, the unique possible singular asymptotics is the Stokes corner (see Fig. 7 ($A_1$)), namely, 
	 	$$
	 \frac{u(x^0+r_mx)}{r_m^{3/2}}\rightarrow u_\infty(x)=	U_{\infty}(\rho,\theta)=\begin{cases}
	 		\begin{matrix}
	 			\frac{\sqrt{2}}{3}\rho^{3/2}\cos\left(\frac{3}{2}\theta-\frac{\pi}{4}\right),&	\text{in}\,\,\left(-\frac{5}{6}\pi,-\frac{\pi}{6}\right),	\\
	 		\end{matrix}\\
	 		\begin{matrix}
	 			0,&	\hspace{32mm}\text{otherwise,}	\\
	 		\end{matrix}\\
	 	\end{cases}
	 	$$
	  as $r_m\rightarrow 0^+$ strongly in $W^{1,2}_{\text{\rm loc}}(\mathbb{R}^2)$ and locally uniformly on $\mathbb{R}^2$ with corresponding density of the negative phase (as defined by \text{Definition C} in Appendix A)
	  	$$
	  M_{x^0,u^-}(0^+)=\frac{\sqrt{3}}{3},
	  $$ 
	   where $u_\infty(x)=U_\infty(\rho,\theta)$ in polar coordinates and $U_{\infty}(\rho,\theta)\leqslant0$.	 
	 
	 	\item For Case 1.2, the singular asymptotics must be one of the following Stokes-type corners (see Fig. 7 ($A_2$)$-$($A_4$)), and we have the convergence
	 $$
	 \frac{u(x^0+r_mx)}{r_m^{3/2}}\rightarrow U_\infty(
	 \rho,\theta)\quad\text{as}\,\,r_m\rightarrow0^+
	 $$
	  strongly in $W^{1,2}_{\text{\rm loc}}(\mathbb{R}^2)$ and locally uniformly on $\mathbb{R}^2$. Specifically,
	 
	 \item[($A_2$)]  Stokes-type corner with single-point contact, namely, $\left(\Gamma_{\text{\rm tp}}\cap B_r(x^0)\right)\backslash\{x^0\}=\varnothing$ for some small enough $r$.
	 $$
	 U_{\infty}(\rho,\theta)=\begin{cases}
	 	\begin{matrix}
	 		\frac{\sqrt{2}}{3}\rho^{3/2}\cos\left(\frac{3}{2}\theta-\frac{3\pi}{4}\right),&	\text{in}\,\,\left(\frac{\pi}{6},\frac{5\pi}{6}\right),	\\
	 	\end{matrix}\\
	 	\begin{matrix}
	 		\frac{\sqrt{2}}{3}\rho^{3/2}\cos\left(\frac{3}{2}\theta-\frac{\pi}{4}\right),&	\hspace{1.4mm}	\text{in}\,\,\left(-\frac{5\pi}{6},-\frac{\pi}{6}\right),	\\
	 	\end{matrix}\\
	 	\begin{matrix}
	 		0,&	\hspace{33.4mm}\text{otherwise,}	\\
	 	\end{matrix}\\
	 \end{cases}
	 $$
	 with the density of the negative phase 	
	 $$
	 M_{x^0,u^-}(0^+)=\frac{\sqrt{3}}{3}.
	 $$ 
	 
	 \item[($A_3$)]  Stokes-type corner with bilateral contact, namely, $\Gamma_{\text{\rm op}}^-\cap B_r(x^0)=\varnothing$.
	 $$
	 U_{\infty}(\rho,\theta)=\begin{cases}
	 	\begin{matrix}
	 		\frac{2}{3}\rho^{3/2}\cos\left(\frac{3}{2}\theta+\frac{3\pi}{4}\right),&	\hspace{0.5mm}	\text{in}\,\,\left[\frac{\pi}{2},\frac{7\pi}{6}\right),	\\
	 	\end{matrix}\\
	 	\begin{matrix}
	 		\frac{2}{3}\rho^{3/2}\cos\left(\frac{3}{2}\theta-\frac{\pi}{4}\right),&	\hspace{2.4mm}	\text{in}\,\,\left[-\frac{\pi}{6},\frac{\pi}{2}\right),	\\
	 	\end{matrix}\\
	 	\begin{matrix}
	 		\frac{\sqrt{6}}{3}\rho^{3/2}\cos\left(\frac{3}{2}\theta-\frac{\pi}{4}\right),&\text{in}\,\,\left[-\frac{5\pi}{6},-\frac{\pi}{6}\right),	\\
	 	\end{matrix}\\
	 \end{cases}
	 $$
	 with the density of the negative phase 	$$
	 M_{x^0,u^-}(0^+)=\frac{\sqrt{3}}{3}.
	 $$ 
	 	 
	 \item[($A_4$)]  Stokes-type corner with unilateral contact, namely, $\Gamma_{\text{\rm op}}^-\cap B_r(x^0)\ne\varnothing$. More precisely,

	 \item[($A_{4.1}$)] (The left-sided unilateral contact)

	 $$
	 U_{\infty}(\rho,\theta)=\begin{cases}
	 	\begin{matrix}
	 		\frac{\sqrt{2\sqrt{3}}}{3}\rho^{3/2}\cos\left(\frac{3}{2}\theta-\pi\right),&	\text{in}\,\,\left(\frac{\pi}{3},\pi\right),	\\
	 	\end{matrix}\\
	 	\begin{matrix}
	 		\frac{\sqrt{2\sqrt{3}}}{3}\rho^{3/2}\cos\left(\frac{3}{2}\theta\right),&	\hspace{6.9mm}	\text{in}\,\,\left(-\pi,-\frac{\pi}{3}\right),	\\
	 	\end{matrix}\\
	 	\begin{matrix}
	 		0,&	\hspace{36mm}\text{otherwise,}	\\
	 	\end{matrix}\\
	 \end{cases}
	 $$
	 	with the density of the negative phase 	
	 	$$
	 	M_{x^0,u^-}(0^+)=\frac{1}{2}.
	 	$$
 \item[($A_{4.2}$)] (The right-sided unilateral contact)
	 $$
	 U_{\infty}(\rho,\theta)=\begin{cases}
	 	\begin{matrix}
	 		\frac{\sqrt{2\sqrt{3}}}{3}\rho^{3/2}\cos\left(\frac{3}{2}\theta-\frac{\pi}{2}\right),&		\text{in}\,\,\left(0,\frac{2\pi}{3}\right),	\\
	 	\end{matrix}\\
	 	\begin{matrix}
	 		\frac{\sqrt{2\sqrt{3}}}{3}\rho^{3/2}\cos\left(\frac{3}{2}\theta-\frac{\pi}{2}\right),&		\text{in}\,\,\left(-\frac{2\pi}{3},0\right),	\\
	 	\end{matrix}\\
	 	\begin{matrix}
	 		0,&	\hspace{36.2mm}\text{otherwise,}	\\
	 	\end{matrix}\\
	 \end{cases}
	 $$
	with the density of the negative phase 	
	$$
	M_{x^0,u^-}(0^+)=\frac{1}{2}.
	$$
	
	 	\item For Case 1.3, the stagnation points must be degenerate, namely,
	 	$$
	 \frac{u(x^0+r_mx)}{r_m^{\alpha}}\rightarrow 0,
	 	$$
 as $r_m\rightarrow 0^+$ strongly in $W^{1,2}_{\text{\rm loc}}(\mathbb{R}^2)$ and locally uniformly on $\mathbb{R}^2$ with the density of the negative phase 	
 $$
 M_{x^0,u^-,\alpha}(0^+)=0,
 $$
  where the strong $L^1_{\text{\rm loc}}$ limit of $\chi_{\{u(x^0+r_mx)<0\}}$ is either $0$ or $1$.
	 \end{itemize} 

  	\begin{figure}[!h]
 	\includegraphics[width=175mm]{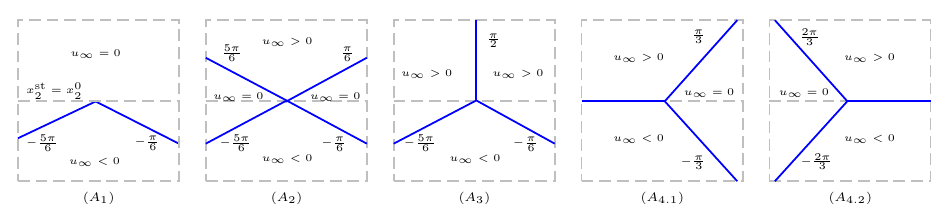}
 	\caption{The possible singular asymptotics}
 \end{figure}
\end{theorem 1} 

\begin{remark} 
	Let $x^0$ be a stagnation point. We classify such points as follows.
	
	$\bullet$ Non-degenerate stagnation points, which satisfy $u_\infty^-(x)\not\equiv0$. These correspond to cases ($A_1$)$-$($A_4$).
	
	$\bullet$ Degenerate stagnation points, which satisfy $u_\infty^-(x)\equiv0$ (see Case 1.3). In this case, the singular asymptotics of the free boundary is determined by the strong $L^1_{\text{\rm loc}}$ limit of $\chi_{\{u(x^0+r_mx)<0\}}$. Specifically, a cusp forms if the limit is $0$, whereas a horizontal point appears if the limit is $1$.
\end{remark}

\begin{remark}
A key contribution of our main results is to elucidate  the influence of the electric field on the singular analysis of the incompressible fluid near stagnation points. 
	
In general, the decay rate of the electric field $\left|\mathbf{E}\right|=\left|\nabla u^+\right|$ relative to the critical decay rate $\left|x-x^{\text{\rm st}}\right|^{1/2}$ governs the singular asymptotics.
\begin{itemize}
	\item [(1)] \textbf{Faster than critical decay rate.} Since the electric field effect is negligible, it follows that the results naturally coincide with the one-phase case in \cite{VW1}. Thus, the singular asymptotics near a stagnation point is Stokes corner.
	
	\item [(2)] \textbf{Critical decay rate.} The singular asymptotics near the stagnation points can take one of two distinct forms, dictated by the electric field distribution. If the electric field is either concentrated at a stagnation point or uniformly distributed along the entire free boundary $\partial\{u<0\}$, axial symmetry is maintained, resulting in the well-known Stokes corner (see  Fig. 7 $\left(A_{2}\right)$ and $\left(A_{3}\right)$). Conversely, if the field acts only on a portion of $\partial\{u<0\}$ near a stagnation point, the resulting force imbalance breaks the axial symmetry, leading to an asymmetric Stokes corner, as illustrated in Fig. 7 $\left(A_{4.1}\right)$ and $\left(A_{4.2}\right)$.
	
	\item [(3)] \textbf{Slower than critical decay rate.} The electric field effect dominates,  completely disrupting the local corner structure of fluid region and resulting in a degenerate stagnation point.
\end{itemize} 
		
		In the special case $\mathbf{E}\equiv0$, by defining $U=-u^->0$ in the fluid region and $U=u^+\equiv0$ in the gas region  within $\Omega$, the scenario ($A_1$) in $\text{Theorem 1}$ recovers the classical Stokes conjecture for steady water waves, which is consistent with $\text{Theorem A}$ in \cite{VW1}. 
\end{remark}

\begin{remark}
This remark completes the classification for Case 1.2. We distinguish two subcases based on the local intersection of the $\partial\{u>0\}$ and $\partial\{u<0\}$ near the stagnation point $x^0$.
\begin{itemize}
	\item [($a$)] {\bf Disjoint boundaries}. If the two boundaries intersect only at $x^0$ locally (i.e.,  $\partial\{u>0\}\cap\partial\{u<0\}\cap B_r\left(x^0\right)=\left\{x^0\right\}$ for some small $r>0$), then the electric field does not affect the symmetry of the corner domain. As shown in ($A_2$), the domain remains symmetric with respect to the direction of gravity. 
	
\item [($b$)] {\bf Touching boundaries}. Otherwise, the configuration depends on the degree of contact between $\partial\{u>0\}$ and $\partial\{u<0\}$, specifically on whether the part of $\Gamma_{\text{\rm op}}^-$ is empty in a neighborhood of $x^0$. 
	\begin{itemize}	
	\item [($b_1$)] {\bf Complete contact}. If $\Gamma_{\text{\rm op}}^-\cap B_r\left(x^0\right)=\varnothing$ for any $r>0$, then $\partial\{u<0\}=\Gamma_{\text{\rm tp}}\subset\partial\{u>0\}$. In this scenario, the electric field acts equally on both sides of the corner region $\{u<0\}$, thus preserving its symmetry.
	\item [($b_2$)] {\bf Partial contact}. If $\Gamma_{\text{\rm op}}^-\cap B_r\left(x^0\right)\ne\varnothing$, then the corner regions of $\{u>0\}$ and $\{u<0\}$ together form a cavity region subtending an angle of $\frac{2\pi}{3}$ (see Fig. 8). Under the assumption that $u\geqslant0$ in $\{x_2\geqslant x_2^0\}$, the free boundary $\partial\{u<0\}\cap B_r(x^0)$ (for sufficiently small $r>0$) must lie in either the third quadrant $\mathcal{Q}_3=\left\{(x_1,x_2)\left|\,x_1\leqslant x_1^0,x_2\leqslant x_2^0\right.\right\}$ or the fourth quadrant $\mathcal{Q}_4=\left\{(x_1,x_2)\left|\,x_1\geqslant x_1^0,x_2\leqslant x_2^0\right.\right\}$. Consequently, the dielectric gas and the incompressible fluid share a common free boundary $\Gamma_{\text{\rm tp}}$ within $B_r(x^0)\backslash\{x^0\}$. This shared boundary is located either in $\mathcal{Q}_3$ (resulting in left-sided unilateral contact ($A_{4.1}$)) or in $\mathcal{Q}_4$ (resulting in right-sided unilateral contact ($A_{4.2}$)). Physically, the shared boundary $\Gamma_{\text{\rm tp}}\cap B_r\left(x^0\right)$ is influenced by both the electric and gravity fields, whereas the one-phase boundary $\Gamma_{\text{\rm op}}^-\cap B_r\left(x^0\right)$ is affected by gravity alone. The resulting force imbalance breaks the axial symmetry. Subsequent calculations confirm that  $\Gamma_{\text{\rm tp}}$ becomes parallel to the $x_2$-axis (see Fig. 7 ($A_4$)).
\end{itemize}
\end{itemize}
\end{remark}
  	\begin{figure}[!h]
	\includegraphics[width=140mm]{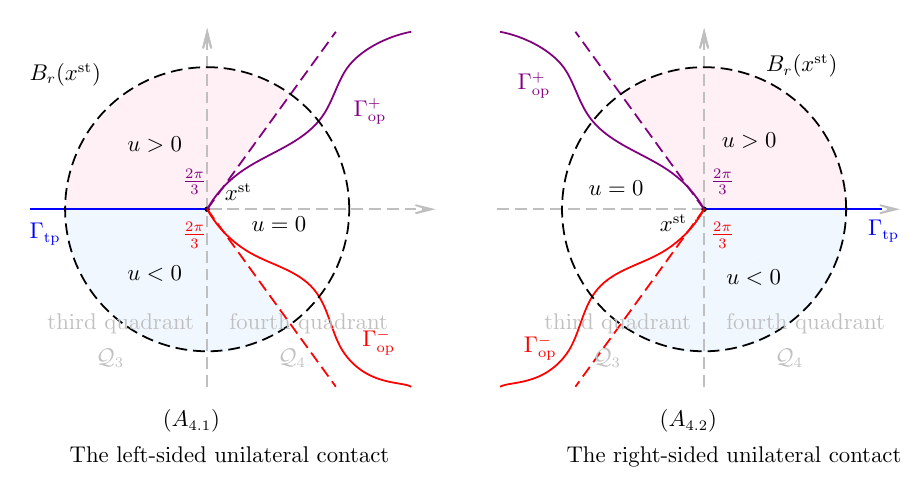}
	\caption{The division of left-sided and right-sided}
\end{figure}
		
\begin{remark}
Consider the following functions
$$
u_m^+=\frac{u^+(r_mx+x^0)}{r_m^{3/2}}\quad\text{and}\quad u_m^-=\frac{u^-(r_mx+x^0)}{r_m^{3/2}}.
$$
for Case 1.1, 1.2 and 
$$
u_m^+=\frac{u^+(r_mx+x^0)}{r_m^{\alpha}}\quad\text{and}\quad u_m^-=\frac{u^-(r_mx+x^0)}{r_m^{\alpha}},
$$
for Case 1.3. Without loss of generality, suppose $u_m^+\rightarrow q_1$ and $u_m^-\rightarrow q_2$ strongly in $W^{1,2}_{\text{\rm loc}}(\mathbb{R}^2)$ and locally uniformly on $\mathbb{R}^2$. Then $u_m=u_m^++u_m^-\rightarrow q_1+q_2$. We emphasize that the limit $u_\infty=u_\infty^++u_\infty^-$ is well-defined, where $q_1=u_\infty^+$ and $q_2=u_\infty^-$. Indeed, by the definition of $u$ (see \eqref{u}), the functions $u^+$ and $u^-$ have disjoint supports. Since this property is preserved under scaling and convergence process, we have $
	\text{supp}\,q_1\cap\text{supp}\, q_2=\varnothing$.
\end{remark}

\begin{remark}
	The asymptotic analysis of $u$ near a stagnation points $x^0$ relies on the transition condition at the free boundary. Moreover, these assumptions are derived from the complete classification of the electric field in $\text{Subsection 1.4}$. 
	\begin{itemize}
		\item In Case 1.1, the transition condition coincides with that of the one-phase flow problem. Consequently, we impose that $\{u=0\}$ has locally only finitely many connected components, that $u\geqslant0$ in $\{x_2\geqslant x_2^0\}$ and that $\left|\nabla u^-\right|^2\leqslant\left|x_2-x_2^0\right|$ in $B_{r_0}(x^0)$. 
		
		\item In Case 1.2, although the free boundary is influenced by both electric and gravity fields, their combined effect near $x^0$ still resembles the one-phase case. Hence, the assumptions of Case 1.1 remain applicable.
		
	\item In Case 1.3, the electric field dominates and decays more slowly than in the one-phase case. We therefore assume the existence of constants $\alpha\in\left(1,3/2\right)$, $r_0>0$ and $C>0$ such that $\left|\nabla u^+\right|\leqslant C\left|x-x^0\right|^{\alpha-1}$ in $B_{r_0}(x^0)$. The upper bound $\alpha<3/2$ arises from the classification, while the lower bound $\alpha>1$ follows from the Lipschitz continuity of $u$ ($\text{Theorem 5.3}$ in \cite{ACF1}), which guarantees a bounded gradient of $u$. Finally, the transition condition justifies the analogous assumption on the negative phase, i.e., $\left|\nabla u^-\right|\leqslant C\left|x-x^0\right|^{\alpha-1}$ in $B_{r_0}(x^0)$.
	\end{itemize}
\end{remark}

In order to analyze the possible shapes of the free boundary close to stagnation points, we establish the following result.

\begin{theorem 1}
	Let $u$ be a weak solution of \eqref{1.20} (as defined in \text{Definition A} in Appendix A). Assume that $u\geqslant0$ in $\{x_2\geqslant x_2^0\}$ and
\begin{itemize}
	\item For Case 1.1 and Case 1.2,
	$$
	\left|\nabla u^-\right|\leqslant \left|x_2-x_2^0\right|^{1/2}\quad \text{in}\,\,B_{r_0}(x^0),
	$$
	\item For Case 1.3, there exists $\alpha$ with $1<\alpha<3/2$ such that
			$$
	\left|\nabla u^+\right|\leqslant C\left|x-x^0\right|^{\alpha-1}\quad \text{and}\quad
\left|\nabla u^-\right|\leqslant C\left|x-x^0\right|^{\alpha-1}\quad \text{in}\,\,B_{r_0}(x^0),
	$$
	where $r_0>0$ small enough and $C>0$.
\end{itemize}	 
Moreover, suppose that $\partial\{u<0\}$ is locally a continuous injective curve $\sigma(t)=\left(\sigma_1(t),\sigma_2(t)\right)$, defined for $t\in(-t_0,t_0)$ with $t_0>0$, and that $\sigma(0)=\left(x_1^0,x_2^0\right)$. Then exactly one of the following holds (see Fig. 9).
	\begin{itemize}
		\item [(1)]{Stokes corner.}
		$$
		M_{x^0,u^-}(0^+)=\frac{\sqrt{3}}{3},
		$$
		in which case $\sigma_1(t)\ne x_1^0$ in $\left(-t_0,t_0\right)\backslash\{0\}$ and 
		$$
		\underset{t\rightarrow 0^+}{\lim}\frac{\sigma_2(t)-x_2^0}{\sigma_1(t)-x_1^0}=-\frac{\sqrt{3}}{3},\quad \underset{t\rightarrow 0^-}{\lim}\frac{\sigma_2(t)-x_2^0}{\sigma_1(t)-x_1^0}=\frac{\sqrt{3}}{3}.
		$$
		\item[(2)] {Asymmetric Stokes corner.}$$
		M_{x^0,u^-}(0^+)=\frac{1}{2}.$$
		
		$\bullet$ Left-sided Stokes-type corner.
		$$
		\underset{t\rightarrow 0^+}{\lim}\frac{\sigma_2(t)-x_2^0}{\sigma_1(t)-x_1^0}=-\sqrt{3},\quad \underset{t\rightarrow 0^-}{\lim}\frac{\sigma_2(t)-x_2^0}{\sigma_1(t)-x_1^0}=0,
		$$
		where $\sigma_1(t)\ne x_1^0$ for all $\left(-t_0,t_0\right)\backslash\{0\}$.
		
		$\bullet$ Right-sided Stokes-type corner.
		$$
		\underset{t\rightarrow 0^+}{\lim}\frac{\sigma_2(t)-x_2^0}{\sigma_1(t)-x_1^0}=0,\quad \underset{t\rightarrow 0^-}{\lim}\frac{\sigma_2(t)-x_2^0}{\sigma_1(t)-x_1^0}=\sqrt{3},
		$$
		where $\sigma_1(t)\ne x_1^0$ for all $\left(-t_0,t_0\right)\backslash\{0\}$.
		
			\item[(3)] {Cusp.}
		$$
		M_{x^0,u^-,\alpha}(0^+)=0,
		$$
		in which case $\sigma_1(t)\ne x_1^0$ in $\left(-t_0,t_0\right)\backslash\{0\}$ and there exists a $\theta_0\in\left[0,\pi\right]$ such that
			$$
		\underset{t\rightarrow 0^+}{\lim}\frac{\sigma_2(t)-x_2^0}{\sigma_1(t)-x_1^0}=	\underset{t\rightarrow 0^-}{\lim}\frac{\sigma_2(t)-x_2^0}{\sigma_1(t)-x_1^0}=\tan\theta_0.
$$
	\end{itemize}	
\end{theorem 1}
\begin{figure}[!h]
	\includegraphics[width=175mm]{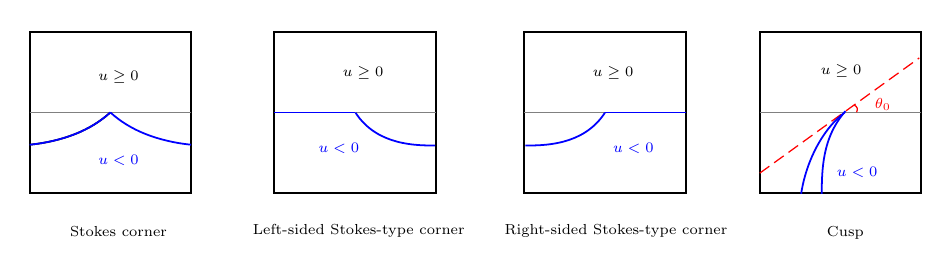}
	\caption{The possible singular profiles}
\end{figure}
\begin{remark}
This work presents an interesting finding in the singular analysis of the two-phase EHD flow under gravity. Specifically, a key insight is the identification of $\left|x-x^{\text{\rm st}}\right|^{1/2}$ as the critical decay rate of the electric field, which governs the formation of singular profiles near the stagnation points. Furthermore, our results bridge a theoretical gap between the purely electric stress-driven EHD problem studied in \cite{GVW} and the one-phase gravity wave problem without electric field examined in \cite{VW1}.
\end{remark}

In conclusion, a central contribution of this work is the rigorous characterization of how the electric field influences singular profiles near $x^{\text{\rm st}}$ with $x^{\text{\rm st}}_2=x_2^0$. The singular profile is entirely determined by the decay rate of the electric field relative to the critical decay rate $\left|x-x^{\text{\rm st}}\right|^{1/2}$.
\begin{itemize}
	\item [(i)] If the electric field decays faster than the critical decay rate $\left|x-x^{\text{\rm st}}\right|^{1/2}$, its influence becomes negligible, yielding Stokes corner.
	
	\item [(ii)] If the electric field decays as the critical decay rate, the singular profile near the stagnation point branches into two distinct types, determined by the extent of the electric field's influence. If the effect is either strictly localized at the stagnation point or uniformly distributed along the entire free boundary, the symmetry of the corner region remains intact, yielding the classical Stokes corner. When the effect acts only on a portion of the free boundary, the symmetry is broken, leading to an asymmetric Stokes corner.
	
	\item [(iii)] If the electric field decays slower than $\left|x-x^{\text{\rm st}}\right|^{1/2}$, the electric field fundamentally destroys the corner geometry of the fluid region and the singular profile must be a cusp.
\end{itemize}

For a clear overview of the possible singular asymptotics near $x^{\text{\rm st}}$ with $x_2^{\text{\rm st}}= x_2^0$, all results are summarized in Tab. 1. Schematics comparing the different cases and their profiles are provided in Tab. 2 and Fig. 10, illustrating the conclusions of $\text{Theorem 2}$.

\renewcommand\arraystretch{1.5}
\begin{table}[h]
	\caption{The possible singular asymptotics as $x_2^{\text{st}}= x_2^0$}\label{table1}
	\centering \scalebox{0.8}{
		\begin{tabular}{c|c|c|c|c|c|c|c}
			\hline
	\multicolumn{2}{c|}{\multirow{3}{*}{Asymptotics}}&\multicolumn{4}{c|}{Non-degenerate point} &\multicolumn{2}{c}{Degenerate point} \\\cline{3-8}
			\multicolumn{2}{c|}{}&\multicolumn{3}{c|}{Stokes corner}&{Asymmetric Stokes corner}&\multirow{2}{*}{Cusp}&\multirow{2}{*}{Horizontal point} \\\cline{3-6}
			\multicolumn{2}{c|}{}&$(A_1)$&$(A_2)$&$(A_3)$&$(A_4)$&&\\\cline{1-8}
		\multicolumn{2}{c|}{Density}&\multicolumn{3}{c|}{	$M_{x^0,u^-}(0^+)=\sqrt{3}/3$}&$M_{x^0,u^-}(0^+)=1/2$
		&\multicolumn{2}{c}{	$M_{x^0,u^-,\alpha}(0^+)=0$}\\\cline{1-8}	
		\multirow{3}{*}{$x_2^{\text{st}}=x_2^0$}&Faster case&$\surd$&$\times$&$\times$&$\times$&$\times$&$\times$	\\\cline{2-8}
		&Critical case&$\times$&$\surd$&$\surd$&$\surd$&$\times$&$\times$\\\cline{2-8}
		&Slower  case&$\times$&$\times$&$\times$&$\times$&$\surd$&$\surd$\\\cline{1-8}
	\end{tabular}}
\end{table}

	\renewcommand\arraystretch{1.5}
\begin{table}[h]
	\caption{The possible singular profiles as $x_2^{\text{st}}=x_2^0$}\label{table2}
	\centering	 	\scalebox{0.8}{
	\begin{tabular}{c|c|c|c|c|c}
		\hline
		\multicolumn{2}{c|}{\multirow{3}{*}{Profiles}}&\multicolumn{3}{c|}{Corner}&\multirow{3}{*}{Cusp}\\\cline{3-5}
	\multicolumn{2}{c|}{}&\multirow{2}{*}{Stokes}	 &Left-sided 	&Right-sided&\\
	\multicolumn{2}{c|}{}&	& Stokes-type&  Stokes-type&\\\cline{1-6}
		\multicolumn{2}{c|}{Density}&	$M_{x^0,u^-}(0^+)=\sqrt{3}/3$&\multicolumn{2}{c|}{$M_{x^0,u^-}(0^+)=1/2$}
	&	$M_{x^0,u^-,\alpha}(0^+)=0$\\\cline{1-6}	
	\multirow{3}{*}{$x_2^{\text{st}}=x_2^0$}&Faster case&$\surd$&$\times$&$\times$&$\times$\\\cline{2-6}
	&Critical case&$\surd$&$\surd$&$\surd$&$\times$\\\cline{2-6}
	&Slower case&$\times$&$\times$&$\times$&$\surd$\\\cline{1-6}	\end{tabular}}
\end{table}
				\begin{figure}[!h]
	\includegraphics[width=170mm]{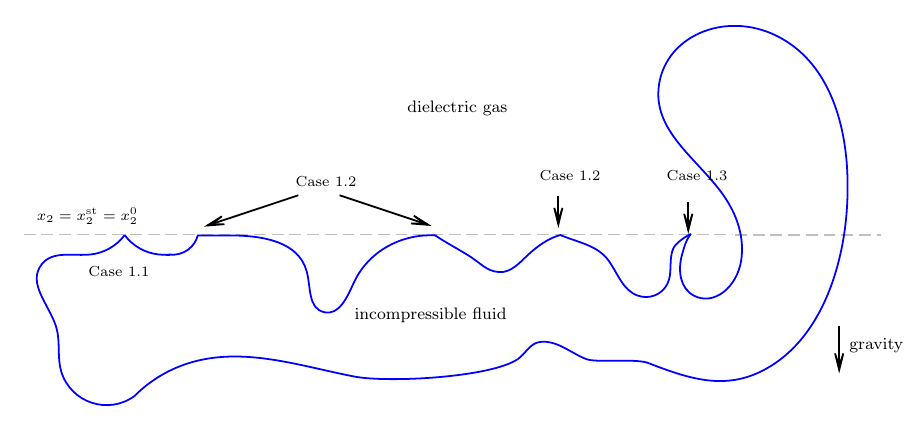}
	\caption{The possible singular profiles of $\Gamma_{\text{\rm tp}}$ for two-phase FBP}
\end{figure}

This work focuses on a two-dimensional gravity-driven fluid surrounded by a dielectric gas. The overall strategy of proof is structured as follows.

Section 2 provides the necessary preliminaries. To analyze the local behavior near stagnation points, we investigate the blowup limits using the Weiss type monotonicity formula. As this approach is standard, we refer to \cite{DHP,DJ2,GVW,VW1,VW2,VW3} for details.

 In Section 3, we perform the blowup analysis for Case 1 with $x_2^{\text{\rm st}}=x_2^0$, categorized by the decay rate of the electric field. Using technical estimates and a frequency formula (detailed in Subsection 3.1), we show that the possible asymptotics in Case 1.1 correspond to a Stokes corner, while in Case 1.2 the possible asymptotics may be either a Stokes corner or a Stokes-type corner. We further prove that the blowup limit must be degenerate near stagnation points in Case 1.3, with the strong $L^1_{\text{\rm loc}}$ limit of $\chi_{\{u(x^0+r_mx)<0\}}$ is either $0$ or $1$. Finally, under the assumption that the free boundary $\partial\{u<0\}$ is a continuous injective curve, we provide a complete classification of all possible singular profiles for stagnation points with $x_2^{\text{st}}=x_2^0$.

\section{Preliminaries}

The primary technique for investigating the local behavior of the free boundary is the blowup method. Specifically, we analyze the blowup limits near stagnation points (refer to \cite{DHP} and \cite{VW1}$-$\cite{W}).
 
A key observation is that the negative-phase solution $u^-$ is a weak solution to the one-phase flow problem
 \begin{equation}
 \label{2.0}
 	\begin{cases}
 		\begin{matrix}
 			\Delta u^-=0,&	\hspace{21.7mm}	\text{in}\,\,\{u<0\},\\
 		\end{matrix}\\
 		\begin{matrix}
 			\left| \nabla u^- \right|^2=\left(x_{2}^{0}-x_2\right)^+,&		\text{on}\,\,\partial \{u<0\}.\\
 		\end{matrix}\\
 	\end{cases}
 \end{equation}
 
 \begin{proposition}
 	Let $u$ be a weak solution of two-phase free boundary problem \eqref{1.20} in $\Omega$. Then $u^-$ is also a weak solution to the one-phase flow problem \eqref{2.0} in $\Omega\cap\{u\leqslant0\}$.
 \end{proposition}
 \begin{pf}
 	By $\text{Definition A and B}$ in $\text{Appendix A}$, we have
 	\begin{equation}
 		\label{2.1}
 		\begin{aligned}
 			0=&\int_{\Omega}\left(\left|\nabla u\right|^2\mathrm{div}\phi-2\nabla uD\phi\nabla u\right)dx+\int_{\Omega\cap\{x_2<x_2^0\}}\left(\left(x_2^0-x_2\right)\chi_{\{u<0\}}\mathrm{div}\phi-\chi_{\{u<0\}}\phi_2\right)dx\\
 			&+\int_{\Omega\cap\{x_2>x_2^0\}}\left(\left(x_2-x_2^0\right)\chi_{\{u>0\}}\mathrm{div}\phi+\chi_{\{u>0\}}\phi_2\right)dx,
 		\end{aligned}
 	\end{equation}
 	for all $\phi(x)=\left(\phi_1(x),\phi_2(x)\right)\in C^1_0\left(\Omega,\mathbb{R}^2\right)$.
 	It suffices to show that $u^-$ satisfies 
 	\begin{equation}
 		\nonumber
 		0=\int_{\Omega\cap\{u\leqslant0\}}\left(\left|\nabla u\right|^2\mathrm{div}\phi-2\nabla uD\phi\nabla u\right)dx+\int_{\Omega\cap\{u\leqslant0\}\cap\{x_2<x_2^0\}}\left(\left(x_2^0-x_2\right)\chi_{\{u<0\}}\mathrm{div}\phi-\chi_{\{u<0\}}\phi_2\right)dx,
 	\end{equation}
 	for all $\phi(x)\in C^1_0\left(\Omega\cap\{u\leqslant0\},\mathbb{R}^2\right)$.
 	
 To prove this, we take an arbitrary $\phi(x)\in C^1_0\left(\Omega\cap\{u\leqslant0\},\,\mathbb{R}^2\right)$ and extend it by zero to all of $\Omega$, denoting this extension by $\tilde{\phi}$. Since $\tilde{\phi}\in C^1_0\left(\Omega,\,\mathbb{R}^2\right)$, we may substitute it into \eqref{2.1} to obtain
 	\begin{equation}
 		\nonumber
 		\begin{aligned}
 			0
 			=&\int_{\Omega}\left(\left|\nabla u\right|^2\mathrm{div}\tilde\phi-2\nabla uD\tilde\phi\nabla u\right)dx+\int_{\Omega\cap\{x_2<x_2^0\}}\left(\left(x_2^0-x_2\right)\chi_{\{u<0\}}\mathrm{div}\tilde\phi-\chi_{\{u<0\}}\tilde\phi_2\right)dx\\
 			&+\int_{\Omega\cap\{x_2>x_2^0\}}\left(\left(x_2-x_2^0\right)\chi_{\{u>0\}}\mathrm{div}\tilde\phi+\chi_{\{u>0\}}\tilde\phi_2\right)dx\\
 			=&\int_{\Omega\cap\{u\leqslant0\}}\left(\left|\nabla u\right|^2\mathrm{div}\phi-2\nabla uD\phi\nabla u\right)dx\\
 			&+\int_{\Omega\cap\{u\leqslant0\}\cap\{x_2<x_2^0\}}\left(\left(x_2^0-x_2\right)\chi_{\{u<0\}}\mathrm{div}\phi-\chi_{\{u<0\}}\phi_2\right)dx.
 		\end{aligned}
 	\end{equation}
 This is precisely the desired identity for $u^-$, which completes the proof.
 \end{pf}

When the stagnation points occur at the uniform height, i.e., $x^{\text{\rm st}}_2=x_2^0$, the transition condition forces the decay rate of $\left|\nabla u^-\right|^2$ to align with either the linear term or $\left|\nabla u^+\right|^2$. The former scenario corresponds to Case 1.1 and 1.2, whereas the latter corresponds to Case 1.3. To formalize this distinction, we introduce two versions of the monotonicity formula for $u$, the proofs of which are provided in Appendix C.

\begin{proposition}	{\rm(For Case 1.1 and 1.2)}
	Let $r_0=\frac{1}{2}\text{dist}\left(x^0,\partial\Omega\right)$. For any $r\in(0,r_0)$ and every weak solution $u$ defined in $\text{Definition A}$ in $\text{Appendix A}$, we introduce the Weiss-type boundary adjusted energy
	$$
M_{x^0,u}(r)=r^{-3}I(r)-\frac{3}{2} r^{-4}J(r),
	$$
	where
	$$
	I(r)=\int_{B_r(x^0)}\left(\left|\nabla u\right|^2+\left(x_2-x_2^0\right)^+\chi_{\{u>0\}}+\left(\left(x_2-x_2^0\right)^++x_2^0-x_2\right)\chi_{\{u<0\}}\right)dx,
	$$
	and
	$$
	J(r)=\int_{\partial B_r(x^0)} u^2dS.
	$$
	Then for almost every $r\in(0,r_0)$,
	$$
	\frac{d}{dr}M_{x^0,u}(r)=2r^{-3}\int_{\partial B_r(x^0)}{\left(\left(\nabla u^+,\nu\right)-\frac{3}{2r}u^+\right)^2+\left(\left(\nabla u^-,\nu\right)-\frac{3}{2r}u^-\right)^2dS},
	$$
	where $\nu$ is the unit outer normal.
\end{proposition}

\begin{proposition}	{\rm(For Case 1.3)}
		For some $\alpha\in\left(1,3/2\right)$, consider the modified Weiss-type boundary adjusted energy for $u$ 
	$$
	M_{x^0,u,\alpha}(r)=r^{-2\alpha}I(r)-\alpha r^{-2\alpha-1}J(r),
	$$
	where $I(r)$, $J(r)$, $r_0$ and $u$ are defined as in \text{Proposition 2.2}. Then for almost every $r\in(0,r_0)$
	$$
	\frac{d}{dr}M_{x^0,u,\alpha}(r)=2r^{-2\alpha}\int_{\partial B_r(x^0)}{\left(\left(\nabla u^+,\nu\right)-\frac{\alpha}{r}u^+\right)^2+\left(\left(\nabla u^-,\nu\right)-\frac{\alpha}{r}u^-\right)^2dS}+K(r),
	$$
	where $\nu$ is the unit outer normal and
	\begin{equation}
		\nonumber
			K(r)=(3-2\alpha)r^{-2\alpha-1}\int_{ B_r(x^0)}\left(\left(x_2-x_2^0\right)^+\chi_{\{u>0\}}+\left(x_2^0-x_2\right)^+\chi_{\{u<0\}}\right)dx.
	\end{equation}
	\end{proposition}

To identify the blowups in a neighborhood of the stagnation point $x^0\in\partial\{u<0\}$, consider the following growth conditions

$\bullet$ For Case 1.1 and 1.2
\begin{equation}
	\label{3.1}
\left|\nabla u^+\right|\leqslant C|x-x^0|^{1/2}\quad \text{and}\quad	\left|\nabla u^-\right|\leqslant C|x-x^0|^{1/2}\quad \text{in}\,\,B_{r}(x^0),
\end{equation}

$\bullet$ For Case 1.3
\begin{equation}
	\label{3.2}
	\left|\nabla u^+\right|\leqslant C|x-x^0|^{\alpha-1}\quad \text{and}\quad\left|\nabla u^-\right|\leqslant C|x-x^0|^{\alpha-1}\quad \text{in}\,\,B_{r}(x^0),
\end{equation}
where $0<r<r_0$ and $1<\alpha<3/2$.
\begin{proposition}
	Let $u$ be a weak solution of \eqref{1.20}. Suppose that $x^0=\left(x^0_1,x^0_2\right)\in\partial\{u<0\}$ and $u$ satisfies \eqref{3.1}. Then the following hold.
	\item[(1)] The limit $M_{x^0,u}(0^+)=\underset{r\rightarrow0^+}{\lim}M_{x^0,u}(r)$ exists and is finite.
	
	\item[(2)] Let $r_m$ be a sequence such that  $0<r_m\rightarrow0^+$ as $m\rightarrow\infty$ and the blowup sequence
	$$
	u_m(x)=\frac{u(x^0+r_mx)}{r_m^{3/2}}
	$$
	converges weakly in $W_{\text{\rm loc}}^{1,2}(\mathbb{R}^2)$ to a blowup limit $u_\infty$ as $m\rightarrow\infty$. Then $u_\infty$ is a homogeneous function of degree $\frac{3}{2}$ in $\mathbb{R}^2$, i.e., $u_\infty(\lambda x)=\lambda^{3/2} u_\infty(x)$ for all $x\in\mathbb{R}^2$ and all $\lambda>0$.
	
	\item[(3)] $u_m$ converges strongly in $W_{\text{\rm loc}}^{1,2}(\mathbb{R}^2)$ and locally uniformly in $\mathbb{R}^2$.
	
	\item[(4)] The density satisfies
	$$M_{x^0,u}(0^+)=	\underset{m\rightarrow\infty}{\lim}\int_{B_1}{\left(x_2\right)^+\chi_{\{u_m>0\}}+\left(-x_2\right)^+\chi_{\{u_m<0\}}dx}.$$
\end{proposition}
\begin{pf}
	(1) From the Weiss-type boundary adjusted energy for $u$, it follows directly that $M_{x^0,u}(r)$ is a monotone increasing function. 
	
Moreover, the growth conditions \eqref{3.1} imply that
	$$
	\left|u^+(x)\right|\leqslant C\left|x-x^0\right|^{3/2}\quad \text{in}\,\, B_{r}(x^0)$$
	 and$$
	\left|u^-(x)\right|\leqslant C\left|x-x^0\right|^{3/2}\quad\text{in}  \,\,B_{r}(x^0).$$ 
	Due to $I(r)\geqslant0$, we have 
	$$
M_{x^0,u}(r)=r^{-3}I(r)- \frac{3}{2} r^{-4}J(r)\geqslant-\frac{3}{2} r^{-4}J(r),
	$$
which ensures that $M_{x^0,u}(r)$ is bounded below. Combining this with monotonicity, we conclude that the limit $M_{x^0,u}(0^+)$ exists and is finite.
	
	(2) The growth conditions imply that $u_m$ is uniformly bounded in $C^{0,1}(B_\delta(x^0))$ for any $0<\delta<\infty$. Consequently, $u_m$ converges locally uniformly to $u_\infty$ in $\mathbb{R}^2$. For any $0<\tau<\sigma<\infty$, we have
	\begin{equation}
		\nonumber
			\begin{aligned}
		&M_{x^0,u}(\sigma r_m)-M_{x^0,u}(\tau r_m)\\
		=&\int_{\tau r_m}^{\sigma r_m}2r^{-3}\int_{\partial B_r(x^0)}\left(\left(\nabla u^+,\nu\right)-\frac{3}{2r}u^+\right)^2+\left(\left(\nabla u^-,\nu\right)-\frac{3}{2r}u^-\right)^2dSdr,
			\end{aligned}
	\end{equation}
		where $\nu$ is the unit outer normal.
	
	On the other hand, a change of variable yields
	\begin{equation}
		\nonumber
		\begin{aligned}
			&\int_{\tau r_m}^{\sigma r_m}2r^{-3}\int_{\partial B_r(x^0)}\left(\left(\nabla u^+,\nu\right)-\frac{3}{2r}u^+\right)^2+\left(\left(\nabla u^-,\nu\right)-\frac{3}{2r}u^-\right)^2dSdr\\
		=&\int_{\tau}^{\sigma}2(r_ms)^{-3}r_m\int_{\partial B_{sr_m}(x^0)}\left(\left(\nabla u^+,\nu\right)-\frac{3}{2sr_m}u^+\right)^2+\left(\left(\nabla u^-,\nu\right)-\frac{3}{2sr_m}u^-\right)^2dSds\\
		=&\int_{\tau}^{\sigma}2s^{-3}\int_{\partial B_{s}}\left(\left(\nabla u_m^+,\nu\right)-\frac{3}{2s}u_m^+\right)^2+\left(\left(\nabla u_m^-,\nu\right)-\frac{3}{2s}u_m^-\right)^2dSds\\
		=&2\int_{B_\sigma\backslash B_\tau}\frac{1}{|x|^3}\left(\left(\left(\nabla u_m^+,x\right)-\frac{3}{2} u_m^+\right)^2+\left(\left(\nabla u_m^-,x\right)-\frac{3}{2} u_m^-\right)^2\right)dx.
		\end{aligned}
	\end{equation}
	Since the limit $M_{x^0,u}(0^+)$ exists, the left-hand side tends to zero as $m\rightarrow\infty$, which implies
	$$\int_{B_\sigma\backslash B_\tau}\frac{1}{|x|^3}\left(\left(\left(\nabla u_m^+,x\right)-\frac{3}{2} u_m^+\right)^2+\left(\left(\nabla u_m^-,x\right)-\frac{3}{2} u_m^-\right)^2\right)dx\rightarrow0\quad\text{as}\,\,m\rightarrow\infty.$$ 
	This is equivalent to
	$$
	\nabla u_\infty^+\cdot x-\frac{3}{2} \, u_\infty^+=0\quad\text{and}\quad	\nabla u_\infty^-\cdot x-\frac{3}{2} \, u_\infty^-=0,\quad\text{for}\,\,x\in\mathbb{R}^2,
	$$
	namely, $u_\infty(\lambda x)=\lambda^{3/2} u_\infty(x)$, where $u_\infty=u_\infty^++u_\infty^-$.
	
	(3) To establish the strong convergence of 
	$u_m$ in $W_{\text{\rm loc}}^{1,2}(\mathbb{R}^2)$ (see $\text{Proposition 3.32}$ in \cite{B}), it suffices to prove the inequalities
	$$
	\int_{\mathbb{R}^2}\left|\nabla u_m^+\right|^2\eta dx\leqslant   \int_{\mathbb{R}^2}\left|\nabla u_\infty^+\right|^2\eta dx,
	$$
	and
		$$
	\int_{\mathbb{R}^2}\left|\nabla u_m^-\right|^2\eta dx\leqslant   \int_{\mathbb{R}^2}\left|\nabla u_\infty^-\right|^2\eta dx,
	$$
		 where $\eta\in C_0^1(\mathbb{R}^2,\mathbb{R}^2)$.
	
To begin, let $\delta=\frac{\text{dist}\left(x^0,\partial\Omega\right)}{2}$. Then, by direct computation,
	\begin{equation}
		\nonumber
		\Delta u_m=\Delta\left(\frac{u(x^0+r_mx)}{r_m^{3/2}}\right)=0\quad\quad\text{in}\,\,B_{\frac{\delta}{r_m}}\cap\{u_m\ne0\},
	\end{equation}
which implies that
	\begin{equation}
		\label{3.5}
		\Delta u_\infty=0\quad \text{in}\,\,\{u_\infty\ne0\}.
	\end{equation}
We now analyze the negative part $u_m^-$. For sufficiently small $\delta>0$, we decompose the integral as follows,
		\begin{equation}
	\label{3.6}
	\begin{aligned}
		&\int_{\mathbb{R}^2}\left|\nabla u_m^-\right|^2\eta dx\\
		=&\int_{\mathbb{R}^2}\left(\nabla u_m^-,\nabla\left(\min\left\{u_m^-+\delta,0\right\}\right)^{1+\delta}\right)\eta dx+
		\int_{\mathbb{R}^2}\left(\nabla u_m^-,\left(\nabla u_m^--\nabla\left(\min\left\{u_m^-+\delta,0\right\}\right)^{1+\delta}\right)\right)\eta dx \\
		=&-\int_{\mathbb{R}^2} \left(\min\left\{u_m^-+\delta,0\right\}\right)^{1+\delta}\Delta u_m^-\eta dx-\int_{\mathbb{R}^2}\left(\min\left\{u_m^-+\delta,0\right\}\right)^{1+\delta}\left(\nabla u_m^-,\nabla\eta\right) dx\\
		&+\int_{\mathbb{R}^2}\left(\nabla u_m^-,\left(\nabla u_m^--\nabla\left(\min\left\{u_m^-+\delta,0\right\}\right)^{1+\delta}\right)\right)\eta dx\\
		=&I_1+I_2+I_3,
	\end{aligned}
\end{equation}
	where
	$$
	I_1=-\int_{\mathbb{R}^2} \left(\min\left\{u_m^-+\delta,0\right\}\right)^{1+\delta}\Delta u_m^-\eta dx,
	$$
	$$
	I_2=-\int_{\mathbb{R}^2}\left(\min\left\{u_m^-+\delta,0\right\}\right)^{1+\delta}\left(\nabla u_m^-,\nabla\eta\right) dx,
	$$
	and
	$$
	I_3=\int_{\mathbb{R}^2}\left(\nabla u_m^-,\left(\nabla u_m^--\nabla\left(\min\left\{u_m^-+\delta,0\right\}\right)^{1+\delta}\right)\right)\eta dx.
	$$
	
We focus on the second term,
	\begin{equation}
		\nonumber
		\begin{aligned}
			I_2=&-\int_{\mathbb{R}^2}\left(\min\left\{u_m^-+\delta,0\right\}\right)^{1+\delta}\left(\nabla u_m^-,\nabla\eta\right) dx\\
			&-\int_{\mathbb{R}^2}\left(\min\left\{u^-_\infty+\frac{\delta}{2},0\right\}\right)^{1+\delta}\left(\left(\nabla u_m^-+\nabla u_\infty^-\right),\nabla\eta\right) dx\\
			&+\int_{\mathbb{R}^2}\left(\nabla u_m^-,\nabla\eta\right) \left(\min\left\{u_\infty^-+\frac{\delta}{2},0\right\}\right)^{1+\delta}dx\\
			&+\int_{\mathbb{R}^2}\left(\nabla u^-_\infty,\nabla\eta\right)\left(\min\left\{u^-_\infty+\frac{\delta}{2},0\right\}\right)^{1+\delta} dx\\
			=&\int_{\mathbb{R}^2}\left(\nabla u_m^-,\nabla\eta\right) \left(\left(\min\left\{u^-_\infty+\frac{\delta}{2},0\right\}\right)^{1+\delta}-\left(\min\left\{u^-_m+\delta,0\right\}\right)^{1+\delta}\right)dx\\
			&-\int_{\mathbb{R}^2}\left(\min\left\{u^-_\infty+\frac{\delta}{2},0\right\}\right)^{1+\delta}\left(\left(\nabla u^-_m-\nabla u^-_\infty\right),\nabla\eta\right) dx\\
			&-\int_{\mathbb{R}^2}\left(\min\left\{u^-_\infty+\frac{\delta}{2},0\right\}\right)^{1+\delta}\left(\nabla u^-_\infty,\nabla\eta \right)dx\\
			=&J_1+J_2+J_3,
		\end{aligned}
	\end{equation}
	where
	$$
	J_1=\int_{\mathbb{R}^2}\left(\nabla u^-_m,\nabla\eta\right) \left(\left(\min\left\{ u_\infty^-+\frac{\delta}{2},0\right\}\right)^{1+\delta}-\left(\min\left\{u_m^-+\delta,0\right\}\right)^{1+\delta}\right)dx,
	$$
	$$
	J_2=-\int_{\mathbb{R}^2}\left(\min\left\{ u^-_\infty+\frac{\delta}{2},0\right\}\right)^{1+\delta}\left(\left(\nabla u_m^--\nabla  u^-_\infty\right),\nabla\eta\right) dx,
	$$
	and
	$$
	J_3=-\int_{\mathbb{R}^2}\left(\min\left\{ u^-_\infty+\frac{\delta}{2},0\right\}\right)^{1+\delta}\left(\nabla  u^-_\infty,\nabla\eta\right) dx.
	$$
	
	By uniform convergence and the continuity of $ u_\infty$, for any $\varepsilon>0$, we can deduce that 
	$$\left|\left(\min\left\{ u^-_\infty+\frac{\delta}{2},0\right\}\right)^{1+\delta}-\left(\min\left\{u_m^-+\delta,0\right\}\right)^{1+\delta}\right|<\varepsilon\quad\text{and}\quad\left|\nabla u_m^--\nabla  u^-_\infty\right|<\varepsilon.$$
	Hence, there exist positive constants $C_1$ and $C_2$ such that
	$$
	|J_1|=\left|\int_{\mathbb{R}^2}\left(\nabla u_m^-,\nabla\eta\right) \left(\left(\min\left\{ u^-_\infty+\frac{\delta}{2},0\right\}\right)^{1+\delta}-\left(\min\left\{u_m^-+\delta,0\right\}\right)^{1+\delta}\right)dx\right|\leqslant C_1\varepsilon,
	$$
	and
	$$
	|J_2|=\left|-\int_{\mathbb{R}^2}\left(\min\left\{ u^-_\infty+\frac{\delta}{2},0\right\}\right)^{1+\delta}\left(\left(\nabla u^-_m-\nabla  u^-_\infty\right),\nabla\eta\right) dx\right|\leqslant C_2\varepsilon.
	$$
	
	For $J_3$, we compute
	\begin{equation}
		\nonumber
		\begin{aligned}
			\left|J_3\right|=&\left|-\int_{\mathbb{R}^2}\left(\min\left\{ u^-_\infty+\frac{\delta}{2},0\right\}\right)^{1+\delta}\left(\nabla  u^-_\infty,\nabla\eta\right) dx\right|\\
			=&\left|-\int_{\mathbb{R}^2}\left(\nabla  u^-_\infty,\nabla\left(\left(\min\left\{ u^-_\infty+\frac{\delta}{2},0\right\}\right)^{1+\delta}\right)\right)\eta dx-\int_{\mathbb{R}^2}\Delta  u^-_\infty\left(\min\left\{ u^-_\infty+\frac{\delta}{2},0\right\}\right)^{1+\delta}\eta dx\right|\\
			\leqslant&\int_{\mathbb{R}_+^2}\left|\nabla  u^-_\infty\right|^2\eta dx\quad\quad\text{as} \,\,\delta\rightarrow0.
		\end{aligned}
	\end{equation}
	
Next, it follows from $\Delta u_m=0$ in $\{u_m\ne0\}$ that $|I_1|=0$.
	
	For any $\varepsilon>0$, since $u_m$ is bounded in $W^{1,\infty}_{\text{loc}}(\mathbb{R}^2)$, there exist $\delta_0(\varepsilon)>0$ and $C_3>0$ such that
	$$
	|I_3|=\left|\int_{\mathbb{R}^2}\left(\nabla u_m^-,\left(\nabla u_m^--\nabla\left(\min\left\{u_m^-+\delta,0\right\}\right)^{1+\delta}\right)\right)\eta dx\right|\leqslant C_3\varepsilon,\quad\text{as}\,\,\delta<\delta_0(\varepsilon).
	$$
	
	Combining the above estimates in \eqref{3.6}, we conclude that
	$$
	\int_{\mathbb{R}^2}\left|\nabla u^-_m\right|^2\eta dx\leqslant\int_{\mathbb{R}^2}\left|\nabla  u^-_\infty\right|^2\eta dx.
	$$
	A similar analysis for $u_m^+$, starting from the decomposition
	\begin{equation}
		\nonumber	
		\begin{aligned}
			&\int_{\mathbb{R}^2}\left|\nabla u_m^+\right|^2\eta dx\\
			=&\int_{\mathbb{R}^2}\left(\nabla u_m^+,\nabla\left(\max\left\{u_m^+-\delta,0\right\}\right)^{1+\delta}\right)\eta dx+
			\int_{\mathbb{R}^2}\left(\nabla u_m^+,\left(\nabla u_m^+-\nabla\left(\max\left\{u_m^+-\delta,0\right\}\right)^{1+\delta}\right)\right)\eta dx \\
			=&-\int_{\mathbb{R}^2} \left(\max\left\{u_m^+-\delta,0\right\}\right)^{1+\delta}\Delta u_m^+\eta dx-\int_{\mathbb{R}^2}\left(\max\left\{u_m^+-\delta,0\right\}\right)^{1+\delta}\left(\nabla u_m^+,\nabla\eta\right) dx\\
			&+\int_{\mathbb{R}^2}\left(\nabla u_m^+,\left(\nabla u_m^+-\nabla\left(\max\left\{u_m^+-\delta,0\right\}\right)^{1+\delta}\right)\right)\eta dx,
		\end{aligned}
	\end{equation}
yields
	$$
	\int_{\mathbb{R}^2}\left|\nabla u^+_m\right|^2\eta dx\leqslant\int_{\mathbb{R}^2}\left|\nabla  u^+_\infty\right|^2\eta dx.
	$$
 Therefore, the sequence $\{u_m\}$ converges strongly in $W^{1,2}_{\text{loc}}(\mathbb{R}^2)$.
	
	(4) To compute $M_{x^0,u}(0^+)$, let $\{u_m\}$ be the sequence defined in part (2). Then
	\begin{equation}
		\nonumber
		\begin{aligned}
		&M_{x^0,u}(0^+)\\=&\underset{m\rightarrow\infty}{\lim}M_{x^0,u}(r_m)\\
			=&\underset{m\rightarrow\infty}{\lim}r_m^{-3}\int_{B_{r_m}(x^0)}\left|\nabla u\right|^2+\left(x_2-x_2^0\right)^+\chi_{\{u>0\}}+\left(x_2^0-x_2\right)^+\chi_{\{u<0\}}dx-\frac{3}{2} r_m^{-4}\int_{\partial B_{r_m}(x^0)}u^2dS\\
			=&\underset{m\rightarrow\infty}{\lim}\int_{B_{1}}\left|\nabla u_m^+\right|^2+\left|\nabla u_m^-\right|^2dx-\frac{3}{2}\int_{\partial B_1}\left(u^+_m\right)^2+\left(u^-_m\right)^2dS\\
			&+\underset{m\rightarrow\infty}{\lim}\int_{  B_1}\left(x_2\right)^+\chi_{\{u_m>0\}}+\left(-x_2\right)^+\chi_{\{u_m<0\}}dx\\
			=&\int_{B_{1}}\left|\nabla u_\infty^+\right|^2+\left|\nabla u_\infty^-\right|^2dx-\frac{3}{2}\int_{\partial B_1}\left(u_\infty^+\right)^2+\left(u_\infty^-\right)^2dS\\
			&+\underset{m\rightarrow\infty}{\lim}\int_{  B_1}\left(x_2\right)^+\chi_{\{u_m>0\}}+\left(-x_2\right)^+\chi_{\{u_m<0\}}dx.
		\end{aligned}
	\end{equation}
	Using an argument analogous to \eqref{3.6} and recalling \eqref{3.5}, we arrive at
	\begin{equation}
		\nonumber
		\begin{aligned}
			&\int_{B_{1}}\left(\nabla u_\infty^+,\nabla\left(\max\left\{u_\infty^+-\delta,0\right\}\right)^{1+\delta}\right)dx\\
			&=\int_{\partial B_1}\left(\nabla u_\infty^+,\nu\right)\left(\max\left\{u_\infty^+-\delta,0\right\}\right)^{1+\delta}dS-\int_{B_{1}}\left(\max\left\{u_\infty^+-\delta,0\right\}\right)^{1+\delta}\Delta u_\infty^+ dx\\
			&=\int_{\partial B_1}\left(\nabla u_\infty^+,\nu\right)\left(\max\left\{u_\infty^+-\delta,0\right\}\right)^{1+\delta}dS\\
			&=\frac{3}{2}\int_{\partial B_1} u_\infty^+\left(\max\left\{u_\infty^+-\delta,0\right\}\right)^{1+\delta}dS,
		\end{aligned}
	\end{equation}
	and similarly,
		\begin{equation}
		\nonumber
		\begin{aligned}
			&\int_{B_{1}}\left(\nabla u_\infty^-,\nabla\left(\min\left\{u_\infty^-+\delta,0\right\}\right)^{1+\delta}\right)dx\\
			&=\int_{\partial B_1}\left(\nabla u_\infty^-,\nu\right)\left(\min\left\{u_\infty^-+\delta,0\right\}\right)^{1+\delta}dS-\int_{B_{1}}\left(\max\left\{u_\infty^-+\delta,0\right\}\right)^{1+\delta}\Delta u_\infty^- dx\\
			&=\int_{\partial B_1}\left(\nabla u_\infty^-,\nu\right)\left(\min\left\{u_\infty^-+\delta,0\right\}\right)^{1+\delta}dS\\
			&=\frac{3}{2}\int_{\partial B_1} u_\infty\left(\min\left\{u_\infty^-+\delta,0\right\}\right)^{1+\delta}dS.
		\end{aligned}
	\end{equation}
	Passing to the limit as $\delta\rightarrow0$, we obtain
	$$
	\int_{B_{1}}\left|\nabla u_\infty^+\right|^2dx=\frac{3}{2}\int_{\partial B_1}\left(u_\infty^+\right)^2dS,
	$$
	and
		$$
	\int_{B_{1}}\left|\nabla u_\infty^-\right|^2dx=\frac{3}{2}\int_{\partial B_1}\left(u_\infty^-\right)^2dS.
	$$
Therefore,
	\begin{equation}
	\nonumber
	\begin{aligned}
M_{x^0,u}(0^+)=&M_{x^0,u^+}(0^+)+M_{x^0,u^-}(0^+)\\
=&\underset{m\rightarrow\infty}{\lim}\int_{ B_1}\left(x_2\right)^+\chi_{\{u_m>0\}}+ \left(-x_2\right)^+\chi_{\{u_m<0\}}dx.
	\end{aligned}
\end{equation}
\end{pf}

 By arguments analogous to those in the proof of $\text{Proposition 2.4}$, we immediately obtain the following result (details are omitted for brevity).
 
\begin{proposition}
	Let $u$ be a weak solution of \eqref{1.20} with $1<\alpha<3/2$, and assume that $u$ satisfies \eqref{3.2}. Then the following statements hold.
	\begin{itemize}
		\item [(i)] The limit $M_{x^0,u,\alpha}(0^+)$ exists and is finite.
		
		\item [(ii)] Let $r_m$ be a sequence such that $0<r_m\rightarrow0^+$ and the blowup sequence
		$$
		u_m(x)=\frac{u(x^0+r_mx)}{r_m^{\alpha}}
		$$
		converges weakly in $W_{\text{\rm loc}}^{1,2}(\mathbb{R}^2)$ to a blowup limit $u_\infty$ as $m\rightarrow\infty$. Then $u_\infty$ is a homogeneous function of degree $\alpha$ in $\mathbb{R}^2$, namely, $u_\infty(\lambda x)=\lambda^\alpha u_\infty(x)$ for all $x\in\mathbb{R}^2$ and $\lambda>0$.
		
		\item [(iii)] $u_m$ converges strongly in $W^{1,2}_{\text{\rm loc}}(\mathbb{R}^2)$ and locally uniformly in $\mathbb{R}^2$.
		
		\item [(iv)] The density $$M_{x^0,u,\alpha}(0^+)=\lim\limits_{m\rightarrow\infty}r_m^{3-2\alpha}\int_{ B_1}\left(x_2\right)^+\chi_{\{u_m>0\}}+\left(-x_2\right)^+\chi_{\{u_m<0\}}dx=0.$$
	\end{itemize}
\end{proposition}

\section{Singular profiles near the stagnation points}
\quad

The main contribution of this section is the rigorous derivation of the possible singular asymptotics of the interface between the incompressible fluid and the dielectric gas near a stagnation point $x^{\text{\rm st}}=x^0=(x_1^0,x_2^0)$ in the three distinct cases. This analysis reveals how the decay rate of the electric field dictates the local singular structure of the free boundary near the stagnation points. 

Under the assumption that the free boundary is a continuous injective curve, we provide a detailed characterization of the singular geometry in each case. Specifically, a key finding is the identification of the critical decay rate $|x-x^{\text{\rm st}}|^{1/2}$. When the electric field decays faster than $|x-x^{\text{\rm st}}|^{1/2}$, its influence is negligible and the singular profile is the classical Stokes corner (a symmetric $\frac{2\pi}{3}$ corner about the vertical axis). When the electric field decays as $|x-x^{\text{\rm st}}|^{1/2}$, it may break the symmetry of the wave profile, resulting in either a Stokes corner or an asymmetric Stokes corner. When the decays slower than $|x-x^{\text{\rm st}}|^{1/2}$, the dominant electric field completely disrupts the corner structure, forming a cusp near the stagnation point.

Following the approaches in \cite{CS} and \cite{FH}, we determine the blowup limits using the eigenfunctions of the Laplace-Beltrami operator on the sphere, which are linked to the homogeneity of  harmonic functions.

In dimension $n=2$, the Laplacian in polar coordinates reads
$$
u_{\rho\rho}+\frac{1}{\rho}u_\rho+\frac{1}{\rho^2}\Delta_\theta u,
$$
where $\Delta_\theta$ denotes the Laplace-Beltrami operator.

Assuming that $u_\infty$ is a homogeneous harmonic function of degree $\frac{3}{2}$, there exists an eigenfunction $f(\theta)$ and we can write $u_\infty(x)=U_\infty(\rho,\theta)=\rho^{3/2}f(\theta)$, where $u_\infty$ represents the blowup limit of
$u_m$. Substituting this form into $\Delta u_\infty=0$ within the region $\{ u_\infty\ne 0\}$ yields
$$
\frac{3}{4}f(\theta)+\frac{3}{2} f(\theta)+f''(\theta)=0.
$$
Solving this equation leads to the general form
\begin{equation}
	\label{U^-}
	u_\infty(x)=U_\infty(\rho,\theta)=C\rho^{3/2}\cos\left(\frac{3}{2}\theta+\varphi\right)
\end{equation} 
 in each connected component of $\{u_\infty\ne0\}$. More precisely, every connected component of $\{u_\infty\ne0$\} constitutes a cone with vertex at the stagnation point $x^0$ and an opening angle of $\frac{2\pi}{3}$.

To simplify the subsequent calculations of the corresponding densities across different cases, we first observe that the characteristic function satisfies the following property.
	\begin{proposition}
Let $u$ be a weak solution of \eqref{1.20} satisfying the assumptions in $\text{Theorem 1}$, and let $u_m$ be the blowup sequence defined in $\text{Proposition 2.4}$. Then $\chi_-$ is a constant in the interior of $\{u_\infty=0\}\cap\{x_2\leqslant0\}$ in Case 1.1 and 1.2, where $\chi_-$ is the strong $L^1_{\text{\rm loc}}$ limit of $\chi_{\{u_m<0\}}$. 
\end{proposition}
\begin{pf} 
	 By the definition of weak solution, for any $\phi(x)=\left(\phi_1,\phi_2\right)\in C_0^1\left(\mathbb{R}^2,\mathbb{R}^2\right)$ and $\phi_m=\left(\phi_{m,1},\phi_{m,2}\right)=\phi\left(\frac{x-x^0}{r_m}\right)$, we have
	
	\begin{equation}
\label{def of weak}
\begin{aligned}
0=&\int_{ B_{r_m}(x^0)}\left(\left|\nabla u\right|^2\mathrm{div}\phi_m-2\nabla uD\phi_m\nabla u\right)dx\\
&+\int_{ B_{r_m}(x^0)\cap\{x_2<x_2^0\}}\left(\left(x_2^0-x_2\right)\chi_{\{u<0\}}\mathrm{div}\phi_m-\phi_{m,2}\chi_{\{u<0\}}\right)dx\\
&+\int_{ B_{r_m}(x^0)\cap\{x_2>x_2^0\}}\left(\left(x_2-x_2^0\right)\chi_{\{u>0\}}\mathrm{div}\phi_m+\phi_{m,2}\chi_{\{u>0\}}\right)dx\\
=&\int_{ B_1}r_m^2\left(\left|\nabla u\left(x^0+r_mx\right)\right|^2\frac{1}{r_m}\mathrm{div}\phi-2\nabla u\left(x^0+r_mx\right)\frac{1}{r_m}D\phi_m\nabla u\left(x^0+r_mx\right)\right)dx\\
&+\int_{ B_1\cap\{x_2<0\}}r_m^2\left(-r_mx_2\chi_{\{u\left(x^0+r_mx\right)<0\}}\frac{1}{r_m}\mathrm{div}\phi-\phi_2\chi_{\{u\left(x^0+r_mx\right)<0\}}\right)dx\\
&+\int_{ B_1\cap\{x_2>0\}}r_m^2\left(r_mx_2\chi_{\{u\left(x^0+r_mx\right)>0\}}\frac{1}{r_m}\mathrm{div}\phi+\phi_2\chi_{\{u\left(x^0+r_mx\right)>0\}}\right)dx\\
=&r_m^{2}\left(\int_{ B_1}\left(\left|\nabla u_m\right|^2\mathrm{div}\phi-2\nabla u_mD\phi\nabla u_m\right)dx+\int_{ B_1\cap\{x_2<0\}}\left(-x_2\chi_{\{u_m<0\}}\mathrm{div}\phi-\phi_2\chi_{\{u_m<0\}}\right)dx\right.\\
&+\left.\int_{ B_1\cap\{x_2>0\}}\left(x_2\chi_{\{u_m>0\}}\mathrm{div}\phi+\phi_2\chi_{\{u_m>0\}}\right)dx\right).
\end{aligned}
	\end{equation}
Dividing by $r_m^2$ and passing to the limit $r_m\rightarrow0$, we obtain
\begin{equation}
\label{3.7}
\begin{aligned}
0=&\int_{ \mathbb{R}^2}\left(\left|\nabla u_\infty\right|^2\mathrm{div}\phi-2\nabla u_\infty D\phi\nabla u_\infty\right)dx\\
&+\int_{ \mathbb{R}^2\cap\{x_2<0\}}\left(-x_2\chi_-\mathrm{div}\phi-\chi_-\phi_2\right)dx+\int_{ \mathbb{R}^2\cap\{x_2>0\}}\left(x_2\chi_+\mathrm{div}\phi+\chi_+\phi_2\right)dx,
\end{aligned}
\end{equation}
where $\chi_+$ and $\chi_-$ are the strong $L^1_{\text{loc}}$ limits of $\chi_{\{u_m>0\}}$ and $\chi_{\{u_m<0\}}$, respectively.

Recall that $M_{x^0,u^-}(0^+)$ denotes the density of $u^-$ in $\{x_2\leqslant x_2^0\}$. By the arbitrariness of $\phi$, we choose $\phi\in C_0^1\left(\mathbb{R}^2\cap\{u_\infty=0\}\cap\{x_2<0\},\mathbb{R}^2\right)$. It follows from \eqref{3.7} that 
$$
\int_{ \mathbb{R}^2\cap\{x_2<0\}}\left(-x_2\chi_-\mathrm{div}\phi-\chi_-\phi_2\right)dx=0,
$$
for all such $\phi$. This implies that $\chi_-$ is a constant in the interior of $\{u_\infty=0\}\cap\{x_2\leqslant0\}$. 
\end{pf}

 \subsection{The possible singular asymptotics for Case 1.1}
 \quad
 
 To complete the proof of Case 1.1 in $\text{Theorem 1}$, we observe that the influence of the electric field can be negligible and $u_\infty^+\equiv0$. Therefore, the analysis reduces to that of $u^-$, where two scenarios arise depending on whether $\{u_\infty<0\}$ is empty. In view of the assumption that $u\geqslant 0$ in $\{x_2\geqslant x_2^0\}$, it follows that $\{u_\infty<0\}$ contains at most one connected component. 
 
Consequently, when $\{u_{\infty}<0\}$ is not empty, the only possible singular asymptotics is Stokes corner.
 \begin{theorem}
 Let $u$ be the weak solution of \eqref{1.20} satisfying the assumptions in $\text{Theorem 1}$, and let $u_m$ be the blowup sequence defined in $\text{Proposition 2.4}$. Then 	
 	$$
 	U_{\infty}^-(\rho,\theta)=\begin{cases}
 		\begin{matrix}
 			\frac{\sqrt{2}}{3}\rho^{3/2}\cos\left(\frac{3}{2}\theta-\frac{\pi}{4}\right),&	\text{in}\,\,\left(-\frac{5}{6}\pi,-\frac{\pi}{6}\right),	\\
 		\end{matrix}\\
 		\begin{matrix}
 			0,&	\hspace{32mm}\text{otherwise,}	\\
 		\end{matrix}\\
 	\end{cases}
 	$$
 	with corresponding density 
 	$$
 	M_{x^0,u^-}(0^+)=\frac{\sqrt{3}}{3},
 	$$
 	where $u_\infty(x)=U_\infty(\rho,\theta)$ in polar coordinates.	 
 	\end{theorem}
 	\begin{pf}
 		By \eqref{U^-} and the assumption $u\geqslant 0$ in $\{x_2\geqslant x_2^0\}$, we have  $$U_\infty^-(\rho,\theta)=C\rho^{\frac{3}{2}}\cos\left(\frac{3}{2}\theta+\varphi\right)\quad\text{in} \,\,\{u_\infty<0\}.$$
 	We now consider the case where $\{u_\infty<0\}$ has one connected component. Since $u_m=\frac{u(x^0+r_mx)}{r_m^{3/2}}$, it follows from \eqref{def of weak} that 
 		\begin{equation}
 	\label{3.4}
 			\begin{aligned}
 				0=&\int_{ B_1}\left(\left|\nabla u_m\right|^2\mathrm{div}\phi-\nabla u_m D\phi\nabla u_m\right)dx+\int_{ B_1\cap\{x_2<0\}}\left(-x_2\chi_{\{u_m<0\}}-\chi_{\{u_m<0\}}\phi_2\right)dx\\
 				&+\int_{ B_1\cap\{x_2>0\}}\left(x_2\chi_{\{u_m>0\}}+\chi_{\{u_m>0\}}\phi_2\right)dx.
 			\end{aligned}
 		\end{equation}
 		Since $\left|\nabla u^+\right|\leqslant\varepsilon\left|x-x^0\right|^{1/2}$ holds in $B_r(x^0)\cap\{u>0\}$ for any sufficiently small $\varepsilon$ in Case 1.1, the blowup sequence $u_m$ satisfies
 		$$\left|\nabla u_m^+\right|=\frac{\left|\nabla u^+(x^0+r_mx)\right|}{r_m^{1/2}}\leqslant\varepsilon\frac{\left|x^0+r_mx-x^0\right|^{1/2}}{r_m^{1/2}}\leqslant\varepsilon|x_2|^{1/2},\quad\text{in}\,\,B_1\cap\{u_m>0\}.$$
	 	Consequently, $u_m^+\rightarrow 0$ strongly in $W^{1,2}_\text{\rm loc}\left(\mathbb{R}^2\right)$ as $m\rightarrow \infty$. Passing to the limit as $m\rightarrow\infty$ in \eqref{3.4} yields
 			\begin{equation}
 				\label{1}
 					\begin{aligned}
 		0=&\int_{ \mathbb{R}^2}\left(\left|\nabla u_\infty^-\right|^2\mathrm{div}\phi-2\nabla u_\infty^- D\phi\nabla u_\infty^-\right)dx\\
 		&+\int_{ \mathbb{R}^2\cap\{x_2<0\}}\left(-x_2\chi_-\mathrm{div}\phi-\chi_-\phi_2\right)dx.
 			\end{aligned}
 		\end{equation}
 		Let $y$ be an arbitrary point in $\partial\{u_\infty<0\}\backslash \{O=(0,0)\}$ and choose $\delta$ such that $B_\delta(y)\subset\{x_2<0\}$. Note that $\partial\{u_\infty<0\}$ has the constant value $\nu(y)$ on $\partial\{u_\infty<0\}\cap B_\delta(y)$. Substituting $\phi(x)=\eta(x)\nu(y)$ into \eqref{1}, where $\eta\in C_0^1(B_\delta(y),\mathbb{R}^2)$ is arbitrary, and integrating by parts, we obtain
 		$$
 		0=\int_{B_\delta(y)\cap\partial\{u_\infty<0\}}\left(-\left|\nabla u_\infty^-\right|^2-x_2\left(1-\chi_-\right)\eta\right)dS.
 		$$
 		Now, $\chi_{\{u_m<0\}}\rightarrow\chi_-$ strongly in $L_{\text{loc}}^1\left(\mathbb{R}^2\cap\{x_2<0\}\right)$ and $\chi_-$ takes values 0 or 1 in the interior of $\{u_\infty=0\}\cap\{x_2\leqslant0\}$. By Hopf's lemma, we must have $\chi_-=0$, and thus $\left|\nabla u_\infty^-\right|^2=-x_2$ on $\partial\{u_\infty<0\}$.
 		
The specific form of the blowup limit satisfying these conditions is well-known
 $$
 u_\infty(x)=U_{\infty}(\rho,\theta)=\begin{cases}
 	\begin{matrix}
 		\frac{\sqrt{2}}{3}\rho^{3/2}\cos\left(\frac{3}{2}\theta-\frac{\pi}{4}\right),&	\text{in}\,\,\left(-\frac{5}{6}\pi,-\frac{\pi}{6}\right),	\\
 	\end{matrix}\\
 	\begin{matrix}
 		0,&	\hspace{32mm}\text{otherwise,}	\\
 	\end{matrix}\\
 \end{cases}
 $$
and we refer to \cite{VW1} for details. Moreover, the density is
	\begin{equation}
	\nonumber
	\begin{aligned}
	M_{x^0,u^-}(0^+)&=\int_{ B_1}-x_2\chi_{\left\{\left(\rho\cos\theta,\rho\sin\theta\right)\left|\, -\frac{5\pi}{6}<\theta<-\frac{\pi}{6}\right.\right\}}dx\\
	&=\int_{0}^{1}\int_{-\frac{5\pi}{6}}^{-\frac{\pi}{6}}-\rho^2\sin\theta d\theta d\rho =\frac{\sqrt{3}}{3}.
	\end{aligned}
\end{equation}
 	\end{pf}
 
 On the other hand, if $u_\infty^-\equiv0$, $\text{Proposition 3.1}$ implies that the possible values of density are limited to the set
$$
M_{x^0,u^-}(0^+)\in\left\{\frac{2}{3},0\right\},
$$
where 
$$\frac{2}{3}=\int_{  B_1}\left(-x_2\right)^+dx=r^{-3}\int_{ B_{r}(x^0)}\left(x_2^0-x_2\right)^+dx.$$

Following the approaches in \cite{VW1}$-$\cite{VW3}, we observe that under a strong version of the Bernstein condition, the density $M_{x^0,u^-}(0^+)=0$ can be ruled out. It is crucial to note that this conclusion depends on the assumptions that $u\geqslant0$ in $\{x_2\geqslant x_2^0\}$ and that $\left|\nabla u^-\right|^2\leqslant\left|x_2-x_2^0\right|$ holds in $B_{r_0}(x^0)$ for some sufficiently small $r_0>0$. Consequently, the result remains valid in Case 1.2.
\begin{theorem}
	Let $u$ be a weak solution of \eqref{1.20} with $x^0\in\partial\{u<0\}$. Suppose $u$ satisfies $u\geqslant0$ in $\{x_2\geqslant x_2^0\}$ and
	$$
	\left|\nabla u^-\right|^2\leqslant |x_2-x_2^0|\quad\text{in}\,\,B_{r_0}(x^0),
	$$
	for some sufficiently small $r_0>0$. Then
	$M_{x^0,u^-}(0^+)\ne0$.
\end{theorem}
\begin{pf}
	Since $x^0\in\partial\{u>0\}$ and $u_m$ converges weakly to $u_\infty$ in $W^{1,2}_{\text{loc}}(\mathbb{R}^2)$, suppose for contradiction that $M_{x^0,u^-}(0^+)=0$. Note that for any $\varphi\in C^\infty_0(\Omega)$
	\begin{equation}
	\nonumber
 \int_{\Omega}u^-\Delta\varphi dx=\int_{\Omega\cap\{u<0\}}\Delta u^-\varphi dx+\int_{\Omega\cap\partial\{u<0\}}\left(u^-\frac{\partial\varphi}{\partial\nu}-\varphi\frac{\partial u^-}{\partial\nu}\right)dS,
	\end{equation}
which implies
	$$
 \int_{\Omega}u^-\Delta\varphi dx=\int_{\Omega\cap\partial\{u<0\}}\varphi\left|\nabla u^-\right|dS.
	$$
	Applying the free boundary conditions $$\left|\nabla u^-\right|^2=x_2^0-x_2\quad\text{on}\,\,\Gamma^-_{\text{op}}\quad\text{and} \quad\left|\nabla u^-\right|^2=\left|\nabla u^+\right|^2+x_2^0-x_2\geqslant x_2^0-x_2\quad\text{on}\,\,\Gamma_{\text{\rm tp}}$$
	to the blowup sequence $u_m$ yields
	$$
	\left|\nabla u^-_m\right|^2=\frac{1}{r_m}\left|\nabla u^-(x^0+r_mx)\right|^2\geqslant\frac{1}{r_m}(-r_mx_2)=-x_2\quad \text{on}\,\,\partial_{\text{red}}\{u_m<0\},
	$$
	where $\partial_{\text{red}}\{u_m<0\}$ denotes the reduced boundary of $\partial\{u_m<0\}$.
	
	Hence, passing to the limit in the sense of measures and using the lower bound for $\left|\nabla u_m^-\right|$, a contradiction arises if the limiting measure of $\Delta u_m^-$ vanishes, i.e.,
	\begin{equation}
	\nonumber
	0\leftarrow\Delta u_m^-(B_1)\geqslant\int_{B_1\cap\partial_{\text{red}}\{u_m<0\}}\sqrt{-x_2}>0\quad \text{as}\,\,m\rightarrow\infty,
	\end{equation}
	where $B_1$ is the unit ball in $\mathbb{R}^2$ and $\Delta u_m^-(B_1)$ is a non-negative Radon measure (see \cite{EG}) in $B_1$.
	
	Indeed, by the maximum principle, there exists at least one connected component $V_m$ of $\{u_m<0\}$ touching the origin and containing a point $x^m\in\partial A$, where $A=(-1,1)\times(-1,0)$. Consider the two possibilities
	$$
	\max\left\{-x_2\left|\,x\in V_m\cap\partial A\right.\right\}\nrightarrow 0\quad\text{as}\,\,m\rightarrow\infty,
	$$ 
	and
	\begin{equation}
\nonumber
		\max\left\{-x_2\left|\,x\in V_m\cap\partial A\right.\right\}\rightarrow 0\quad\text{as}\,\,m\rightarrow\infty.
	\end{equation}
Observe that
		\begin{equation}
	\nonumber
		\begin{aligned}
			0&=\Delta u_m^-\left(V_m\cap A\right)\\
			&=\int_{\partial V_m\cap A}\left(\nabla u_m^-,\nu\right) dS+\int_{V_m\cap \partial A}\left(\nabla u_m^-,\nu\right) dS
			\\
			&\leqslant\int_{\partial V_m\cap A}-\left|\nabla u_m^-\right|dS+\int_{V_m\cap \partial A}\left|\nabla u_m^-\right|dS\\
			&\leqslant\int_{\partial V_m\cap A}-\sqrt{-x_2}dS+\int_{ V_m\cap\partial A}\sqrt{-x_2}dS.
		\end{aligned}
	\end{equation}
We complete the proof by analyzing the two cases, each of which leads to a contradiction, following arguments that parallel the proof of $\text{Lemma 4.4}$ in \cite{VW1}.
\end{pf}

 To complete the proof of Case 1.1 in $\text{Theorem 1}$, we introduce a powerful tool, the so-called frequency formula, which will enable us to prove that the case
$$
		M_{x^0,u^-}(0^+)	=\int_{ B_{r}(x^0)}\left(x_2^0-x_2\right)^+dx= \int_{ B_1} \left(-x_2\right)^+dx
		=\int_{0}^{1}\int_{-\pi}^{0}-\rho^2\sin\theta d\theta d\rho=\frac{2}{3}.		
$$
is impossible.

\begin{definition}
	Let $u$ be a weak solution of \eqref{1.20}. Define the set
	$$
	\Sigma^u=\left\{x\in \partial\{u<0\}\left|\,x_2=x_2^0,\,\,\left|\nabla u^-(x)\right|=0\,\, \text{and}\,\, M_{x^0,u^-}(0^+)=\frac{2}{3}\right.\right\}.
	$$
\end{definition}

\begin{remark}
	By $\text{Remark 6.2}$ in \cite{VW1}, the mapping  $x\longmapsto M_{x,u^-}(0^+)$ is upper semicontinuous, which implies that the set $\Sigma^u$ is closed.
\end{remark}

\begin{theorem}
	Let $x^0$ be a point of $\Sigma^u$ and set $\delta=\frac{1}{2}\text{dist}\left(x^0,\partial\Omega\right)$. Define the functions
	$$
D_{x^0,u^-}(r)=\frac{r\int_{ B_r(x^0)}{\left|\nabla u^-\right|^2dx}}{\int_{\partial B_r(x^0)}\left(u^-\right)^2dS},
	$$
	and
	$$
V_{x^0,u^-}(r)=\frac{r\int_{ B_r(x^0)}\left(x_2^0-x_2\right)^+\left(1-\chi_{\{u<0\}}\right)dx}{\int_{\partial B_r(x^0)}\left(u^-\right)^2dS}.
	$$
	Then the ``frequency"
	$$
H_{x^0,u^-}(r)=D_{x^0,u^-}(r)-V_{x^0,u^-}(r)=\frac{r\int_{ B_r(x^0)}\left|\nabla u^-\right|^2dx}{\int_{\partial B_r(x^0)}\left(u^-\right)^2dS}-\frac{r\int_{ B_r(x^0)}\left(x_2^0-x_2\right)^+\left(1-\chi_{\{u<0\}}\right)dx}{\int_{\partial B_r(x^0)}\left(u^-\right)^2dS},
	$$
	satisfies
	\begin{equation}
	\label{4.1}
	\begin{aligned}
	\frac{d}{dr}H_{x^0,u^-}(r)=&\frac{2}{r}\int_{\partial B_r(x^0)}{\left(\frac{r\left(\nabla u^-,\nu\right)}{\left(\int_{\partial B_r(x^0)}\left(u^-\right)^2dS\right)^{1/2}}-D_{x^0,u^-}(r)\frac{u^-}{\left(\int_{\partial B_r(x^0)}\left(u^-\right)^2dS\right)^{1/2}}\right)^2dS}\\
	&+\frac{2}{r}\left(D_{x^0,u^-}(r)V_{x^0,u^-}(r)-\frac{3}{2}V_{x^0,u^-}(r)\right),
	\end{aligned}
	\end{equation}
	for a.e. $r\in\left(0,\delta\right)$.
\end{theorem}
\begin{pf}
	Define
	$$\tilde{I}(r)=r^{-3}\int_{ B_r(x^0)}\left(\left|\nabla u^-\right|^2+\left(x_2^0-x_2\right)^+\chi_{\{u<0\}}\right)dx,$$
	 and 
	 $$\tilde{J}(r)=r^{-4}\int_{\partial B_r(x^0)}\left(u^-\right)^2dS.$$
	  By repeating the calculation in $\text{Proposition 2.2}$, we obtain
	  	\begin{equation}
	  	\nonumber
	  	\begin{aligned}
	  \frac{d}{dr}\tilde{I}(r)=&-3r^{-4}\int_{ B_r(x^0)}\left(\left|\nabla u^-\right|^2+\left(x_2^0-x_2\right)^+\chi_{\{u<0\}}\right)dx\\
	  &+r^{-3}\left(2\int_{\partial B_r(x^0)}\left(\nabla u^-,\nu\right)^2dS+\frac{3}{r}\int_{ B_r(x^0)}\left(x_2^0-x_2\right)^+\chi_{\{u<0\}}dx\right)\\
	  =&2r^{-3}\int_{\partial B_r(x^0)}\left(\nabla u^-,\nu\right)^2dS-3r^{-4}\int_{\partial B_r(x^0)}u^-\left(\nabla u^-,\nu\right)dS,
	 	\end{aligned}
 \end{equation}
 and
 	\begin{equation}
 	\nonumber
 	\begin{aligned}
 		\frac{d}{dr}\tilde{J}(r)=&-4r^{-5}\int_{\partial B_r(x^0)}\left(u^-\right)^2dS+r^{-4}\left(\frac{1}{r}\int_{\partial B_r(x^0)}\left(u^-\right)^2dS+2\int_{\partial B_r(x^0)}u^-\left(\nabla u^-,\nu\right)dS\right)\\
 		=&-3r^{-5}\int_{\partial B_r(x^0)}\left(u^-\right)^2dS+2r^{-4}\int_{\partial B_r(x^0)}u^-\left(\nabla u^-,\nu\right)dS.
 			 	\end{aligned}
 	\end{equation}
	   Hence, for a.e. $r\in\left(0,\delta\right)$ we have
	\begin{equation}
	\label{4.2}
	\begin{aligned}
	&\frac{d}{dr}H_{x^0,u^-}(r)\\=&\frac{	\frac{d}{dr}\tilde{I}(r)}{\tilde{J}(r)}-\frac{\left(\tilde{I}(r)-\int_{B_{1}}\left(-x_2\right)^+dx\right)	\frac{d}{dr}\tilde{J}(r)}{\tilde{J}^2(r)}\\
	=&\frac{2\int_{\partial B_r(x^0)}{\left(\nabla u^-,\nu\right)^2dS}-3r\int_{\partial B_r(x^0)}u^-\left(\nabla u^-,\nu\right)dS}{r\int_{\partial B_r(x^0)}\left(u^-\right)^2dS}\\
	&-\left(D_{x^0,u^-}(r)-V_{x^0,u^-}(r)\right)\frac{-3r\int_{\partial B_r(x^0)}\left(u^-\right)^2dS+2\int_{\partial B_r(x^0)}u^-\left(\nabla u^-,\nu\right)dS}{r\int_{\partial B_r(x^0)}\left(u^-\right)^2dS}\\
	=&\frac{2}{r}\left(\frac{r^2\int_{\partial B_r(x^0)}\left(\nabla u^-,\nu\right)^2dS}{\int_{\partial B_r(x^0)}\left(u^-\right)^2dS}-\frac{3}{2}D_{x^0,u^-}(r)\right)-\frac{2}{r}\left(D_{x^0,u^-}(r)-V_{x^0,u^-}(r)\right)\left(D_{x^0,u^-}(r)-\frac{3}{2}\right),
	\end{aligned}
	\end{equation}
	where the last equation uses the identity
	\begin{equation}
	\label{4.3}
D_{x^0,u^-}(r)=\frac{r\int_{\partial B_r(x^0)}u^-\left(\nabla u^-,\nu\right)dS}{\int_{\partial B_r(x^0)}\left(u^-\right)^2dS}.
	\end{equation}
	Substituting \eqref{4.3} and the relation $D_{x^0,u^-}(r)=H_{x^0,u^-}(r)+V_{x^0,u^-}(r)$ into \eqref{4.2}, we obtain \eqref{4.1}.
\end{pf}
\begin{remark}
	By $D_{x^0,u^-}(r)=H_{x^0,u^-}(r)+V_{x^0,u^-}(r)$, it is easy to verify that
	\begin{equation}
	\label{4.4}
	\begin{aligned}
		\frac{d}{dr}H_{x^0,u^-}(r)=&\frac{2}{r}\int_{\partial B_r(x^0)}{\left(\frac{r\left(\nabla u^-,\nu\right)}{\left(\int_{\partial B_r(x^0)}\left(u^-\right)^2dS\right)^{1/2}}-H_{x^0,u^-}(r)\frac{u^-}{\left(\int_{\partial B_r(x^0)}\left(u^-\right)^2dS\right)^{1/2}}\right)^2dS}\\
	&+\frac{2}{r}V_{x^0,u^-}(r)\left(H_{x^0,u^-}(r)-\frac{3}{2}\right),
	\end{aligned}
	\end{equation}
	for a.e. $r\in\left(0,\delta\right)$.
\end{remark}

Next, we establish several key properties of $H_{x^0,u^-}(r)$.

\begin{proposition}
	Let $u$ be a weak solution of \eqref{1.20} satisfying the assumptions in $\text{Theorem 1}$, with $x^0\in\Sigma^u$ and $\delta=\frac{1}{2}\text{dist}\left(x^0,\partial\Omega\right)$. Then for $r\in\left(0,\delta\right)$ sufficiently small, the following hold.
	
	\item[(1)] $	\frac{d}{dr}H_{x^0,u^-}(r)\geqslant\frac{2}{r}\left(V_{x^0,u^-}(r)\right)^2$.
	
	\item[(2)] The limit  $\underset{r\rightarrow0^+}{\lim}H_{x^0,u^-}(r)=H_{x^0,u^-}(0^+)$ exists.
	
	\item[(3)] The function $r\longmapsto \frac{2}{r}\left(V_{x^0,u^-}(r)\right)^2$ is integrable on $\left(0,\delta\right)$.
	
	\item[(4)] The function $r\longmapsto r^{-4}\int_{\partial B_r(x^0)}\left(u^-\right)^2dS$ is non-decreasing.
	
	\item[(5)] The function $r\longmapsto r^{-4}\int_{B_r(x^0)}\left(x^0_2-x_2\right)^+\left(1-\chi_{\{u<0\}}\right)dx$ is integrable on $\left(0,\delta\right)$.
\end{proposition}
\begin{pf}
	(1) From $\text{Proposition 2.2}$, we have 
	$$
	M_{x^0,u^-}(r)=r^{-3}\int_{ B_r(x^0)}\left(\left|\nabla u^-\right|^2+\left(x^0_2-x_2\right)^+\chi_{\{u<0\}}\right)dx-\frac{3}{2}\int_{\partial B_r(x^0)}\left(u^-\right)^2dS,$$
	and
	$$
	\frac{d}{dr}M_{x^0,u^-}(r)=2r^{-3}\int_{\partial B_r(x^0)}\left(\left(\nabla u^-,\nu\right)-\frac{3}{2r}\right)^2dS.
	$$
 Since 
	$$M_{x^0,u^-}(r)-M_{x^0,u^-}(0^+)=\int_{0}^{r}	\frac{d}{dr}M_{x^0,u^-}(s)ds\geqslant0,$$
  multiplying both sides by $r^4/\int_{\partial B_r(x^0)}\left(u^-\right)^2dS$ yields
	\begin{equation}
	\label{4.5}
	\frac{r\int_{ B_r(x^0)}\left|\nabla u^-\right|^2dx}{\int_{\partial B_r(x^0)}\left(u^-\right)^2dS}-\frac{r\int_{ B_r(x^0)}\left(x^0_2-x_2\right)^+\left(1-\chi_{\{u<0\}}\right)dx}{\int_{\partial B_r(x^0)}\left(u^-\right)^2dS}-\frac{3}{2}\geqslant0,
	\end{equation}
	which is equivalent to $H_{x^0,u^-}(r)=D_{x^0,u^-}(r)-V_{x^0,u^-}(r)\geqslant\frac{3}{2}$.
	
Now, from \eqref{4.1} and \eqref{4.5}, we obtain 
	\begin{equation}
	\label{4.6}
	\frac{d}{dr}H_{x^0,u^-}(r)\geqslant\frac{2}{r}V_{x^0,u^-}(r)\left(D_{x^0,u^-}(r)-\frac{3}{2}\right)\geqslant\frac{2}{r}\left(V_{x^0,u^-}(r)\right)^2.
	\end{equation}
	
	(2) Combining the monotonicity $	\frac{d}{dr}H_{x^0,u^-}(r)\geqslant0$ from (1) and  $H_{x^0,u^-}(r)\geqslant\frac{3}{2}$, it easy to show that the limit $H_{x^0,u^-}(0^+)$ exists.
	
	(3) By \eqref{4.6}, the integrability of $r\longmapsto \frac{2}{r}\left(V_{x^0,u^-}(r)\right)^2$  on $\left(0,\delta\right)$ is equivalent to that of $\frac{d}{dr}H_{x^0,u^-}(r)$, since
	$$
	\int_{0}^{\delta}\frac{2}{r}\left(V_{x^0,u^-}(r)\right)^2dr\leqslant	\int_{0}^{\delta}	\frac{d}{dr}H_{x^0,u^-}(r)dr<\infty.
	$$
This inequality also establishes the existence of the limit $H_{x^0,u^-}(0^+)$, thereby completing the proof of (2).
	
	(4) The proof is similar to $\text{Corollary 7.4}$ in \cite{VW1} and is omitted here.
	
	(5) It suffices to show that
	$$
	\frac{d}{dr}\left(r^{-4}\int_{\partial B_r(x^0)}\left(u^-\right)^2dS\right)\geqslant0.
	$$
	Indeed, using $M_{x^0,u^-}(r)\geqslant M_{x^0,u^-}(0^+)$, we have
	\begin{equation}
	\label{4.7}
	\begin{aligned}
	\frac{d}{dr}\left(r^{-4}\int_{\partial B_r(x^0)}\left(u^-\right)^2dS\right)&=2r^{-5}\left(r\int_{ B_r(x^0)}\left|\nabla u^-\right|^2dx-\frac{3}{2}\int_{\partial B_r(x^0)}\left(u^-\right)^2dS\right)\\
	&\geqslant2r^{-5}\left(r\int_{ B_r(x^0)}\left(x^0_2-x_2\right)^+\left(1-\chi_{\{u<0\}}\right)dx\right)\geqslant0.
	\end{aligned}
	\end{equation}
	
	(6) Consider $$u_m^-=\frac{u^-(x^0+r_mx)}{r_m^{3/2}},$$
	and observe that
	\begin{equation}
	\nonumber
	\begin{aligned}	
	\underset{r_m\rightarrow0^+}{\lim}r_m^{-4}\int_{\partial B_{r_m}(x^0)}\left(u^-\right)^2dS&=	\underset{r_m\rightarrow0^+}{\lim}r_m^{-4}\int_{\partial B_1}\left(u^-\left(x^0+r_my\right) \right)^2r_mdS\\
	&=\underset{r_m\rightarrow0^+}{\lim}\int_{\partial B_1}\left(u_m^-\right)^2dS\\&=\int_{\partial B_1}\left(u_\infty^-\right)^2dS=0,
	\end{aligned}
	\end{equation}
	as $x^0\in\Sigma^u$. Combining this with\eqref{4.7}, we deduce that
	$$
	\int_{0}^{\delta}t^{-4}\int_{ B_t(x^0)}\left(x^0_2-x_2\right)^+\left(1-\chi_{\{u<0\}}\right)dxdt\leqslant\frac{1}{2}\int_{0}^{\delta}\frac{d}{dt} \left(t^{-4}\int_{\partial B_t(x^0)}\left(u^-\right)^2dS\right)dt<\infty,
	$$
	since the fact that $r^{-4}\int_{\partial B_r(x^0)}\left(u^-\right)^2dS$ is bounded as $r\rightarrow 0^+$.
\end{pf}

Based on the established framework, we now prove that the frequency formula passes to the limit.

\begin{proposition}
	Let $u$ be a weak solution of \eqref{1.20} with $x^0\in\Sigma^u$. Then
	
	\item[(1)] The limits 
	$$
	\underset{r\rightarrow0^+}{\lim}V_{x^0,u^-}(r)=0\quad\text{and}\quad\underset{r\rightarrow0^+}{\lim}D_{x^0,u^-}(r)=H_{x^0,u^-}(0^+)\quad\text{exist.}
	$$
	
	\item[(2)] For any sequence $r_m\rightarrow0^+$ as $m\rightarrow\infty$, the blowup sequence
	\begin{equation}
	\label{4.8}
	v_m(x)=\frac{u^-\left(x^0+r_mx\right)}{\sqrt{r_m^{-1}\int_{\partial B_{r_m}(x^0)}\left(u^-\right)^2dS}}
	\end{equation}
	is bounded in $W^{1,2}(B_1)$.
	
	\item[(3)] For any sequence $r_m\rightarrow0^+$ as $m\rightarrow\infty$ such that $v_m$ in \eqref{4.8} converges weakly in $W^{1,2}(B_1)$ to a blowup limit $w$ and the function $w$ is homogeneous with degree $H_{x^0,u^-}(0^+)$, which satisfies
	$$
	w\leqslant0\quad\text{in}\,\,B_1,\quad w\equiv0 \quad\text{in}\,\,B_1\cap\{x_2\geqslant0\}\quad\text{and}\quad \int_{\partial B_{1}}w^2dS=1.
	$$
\end{proposition}

\begin{pf}
	The proof relies on an key estimate
	\begin{equation}
	\label{4.9}
	\int_{B_\sigma\backslash B_\rho}2|x|^{-5}\left(\left(\nabla v_m, x\right)-H_{x^0,u^-}(0^+)v_m\right)^2dx\rightarrow0,
	\end{equation}
	for any sequence $r_m\rightarrow0$ as $m\rightarrow\infty$ and $0<\rho<\sigma<1$.
	
	Indeed, integrating \eqref{4.4} with respect to $r$ from $r_m\rho$ to $r_m\sigma$ and employing \eqref{4.5} yields 
	\begin{equation}
	\nonumber
	\begin{aligned}	
	&\int_{r_m\rho}^{r_m\sigma}{\frac{2}{r}\int_{\partial B_r(x^0)}\left(\frac{r\left(\nabla u^-,\nu\right)}{\left(\int_{\partial B_r(x^0)}\left(u^-\right)^2dS\right)^{1/2}}-H_{x^0,u^-}(r)\frac{u^-}{\left(\int_{\partial B_r(x^0)}\left(u^-\right)^2dS\right)^{1/2}}\right)^2dSdr}\\
	=&\int_{\rho}^{\sigma}\frac{2}{s}\int_{ \partial B_s}\left(\frac{s\left(\nabla v_m,\nu\right)}{\left(\int_{\partial B_s}\left(v_m\right)^2dS\right)^{1/2}}-H_{x^0,u^-}(r_ms)\frac{v_m}{\left(\int_{\partial B_s}\left(v_m\right)^2dS\right)^{1/2}}\right)^2dSds\\
	\leqslant&\int_{r_m\rho}^{r_m\sigma}	\frac{d}{dr}H_{x^0,u^-}(r)dr\rightarrow0,\quad\quad\text{as}\,\,m\rightarrow\infty.
	\end{aligned}	
	\end{equation}
From (4) in $\text{Proposition 3.5}$, we have
	$$
	\int_{\partial B_s}\left(v_m\right)^2dS=\frac{\int_{\partial B_{r_ms}(x^0)}\left(u^-\right)^2dS}{\int_{\partial B_{r_m}(x^0)}\left(u^-\right)^2dS}=\frac{\left(r_ms\right)^{-4}\int_{\partial B_{r_ms}(x^0)}\left(u^-\right)^2dS}{r_m^{-4}\int_{\partial B_{r_m}(x^0)}\left(u^-\right)^2dS}\, s^4\leqslant s^4,
	$$	
	which proves \eqref{4.9}.
	
From (2) in $\text{Proposition 3.5}$, we first note that the existence of $V_{x^0,u^-}(0^+)$ is equivalent to that of $D_{x^0,u^-}(0^+)$. 
	
	To prove $\underset{r\rightarrow0^+}{\lim}V_{x^0,u^-}(r)=0$, assume, to the contrary, that there exists a sequence $s_m\rightarrow0$ such that $V_{x^0,u^-}(s_m)$ is bounded away from $0$. Applying (3) in $\text{Proposition 3.5}$, we obtain
	$$
	\underset{r\in\left[s_m,2s_m\right]}{\min}V_{x^0,u^-}(r)\rightarrow0\quad \text{as}\,\,m\rightarrow\infty.
	$$
	Choose $r_m\in\left[s_m,2s_m\right]$ such that $V_{x^0,u^-}(r_m)\rightarrow0$ as $m\rightarrow0$. Note that
	$$
D_{x^0,u^-}(r_m)=\frac{r_m\int_{ B_{r_m}(x^0)}\left|\nabla u^-\right|^2dx}{\int_{\partial B_{r_m}(x^0)}\left(u^-\right)^2dS}=\int_{B_{1}}\left|\nabla v_m\right|^2dx.
	$$
	Passing to the limit as $m\rightarrow\infty$, we conclude that $D_{x^0,u^-}(r_m)$ is bounded. The same argument as in \cite{VW1}$-$\cite{VW3} shows that $v_m\in W^{1,2}(B_1)$ and this completes the proof of (2).
	
	Furthermore, let $w$ be the weak limit of $v_m$. By the compact embedding $W^{1,2}\hookrightarrow L^2(\partial B_1)$, we have $\lVert v_m \rVert_{L^2(\partial B_1)}\rightarrow\lVert w \rVert_{L^2(\partial B_1)}=1$ as $m\rightarrow\infty$. The homogeneity degree of $w$ follows from \eqref{4.9}. Moreover, since $u\leqslant0$ in $B_{r_m}(x^0)$ and $u\geqslant0$ in $B_{r_m}(x^0)\cap\{x_2\geqslant x_2^0\}$, we conclude $w\leqslant0$ in $B_1$ and $w\equiv0$ in $B_1\cap\{x_2\geqslant0\}$. This completes the proof of the statement (3).
	
	Finally, by part (5) in $\text{Proposition 3.5}$, for $s_m\leqslant r_m\leqslant2s_m$, we have
	\begin{equation}
	\label{4.10}
	\begin{aligned}
V_{x^0,u^-}(s_m)&\leqslant\frac{s_m\int_{ B_{r_m}(x^0)}\left(x^0_2-x_2\right)^+\left(1-\chi_{\{u<0\}}\right)dx}{\int_{\partial B_{r_m/2}(x^0)}\left(u^-\right)^2dS}\\
	&\leqslant\frac{r_m\int_{ B_{r_m}(x^0)}\left(x^0_2-x_2\right)^+\left(1-\chi_{\{u<0\}}\right)dx}{\int_{\partial B_{1/2}}\left(v_m\right)^2dS\int_{\partial B_{r_m/2}(x^0)}\left(u^-\right)^2dS}\\
	&\leqslant\frac{V_{x^0,u^-}(r_m)}{\int_{\partial B_{1/2}}\left(v_m\right)^2dS}.
	\end{aligned}
	\end{equation}
	By the homogeneity of $w$ and $\int_{\partial B_1}w^2=1$, there exists a subsequence satisfying
	$$
	\int_{\partial B_{1/2}}\left(v_m\right)^2dS\rightarrow\int_{\partial B_{1/2}}w^2dS>0.
	$$
	Hence, \eqref{4.10} contradicts our assumption, which implies $\underset{r\rightarrow0^+}{\lim}V_{x^0,u^-}(r)=0$.
\end{pf}

For the subsequent blowup analysis, we introduce the following concentration compactness in two dimensions, which ensures the strong convergence of our blowup sequence. We omit the proof here and refer the reader to \cite{EM} and \cite{VW1} for details.

\begin{proposition}
	Let $u$ be a weak solution of \eqref{1.20}  satisfying the assumptions in $\text{Theorem 1}$ with $x^0\in\Sigma^u$. Let $r_m\rightarrow0^+$ be such that the sequence $v_m$ defined in \eqref{4.8} converges weakly to $w$ in $W^{1,2}(B_1)$. Then 
	\item[(1)]  $v_m$ converges to $w$ strongly in $W^{1,2}_{\text{\rm loc}}(B_1\backslash\{O\})$, where $O=(0,0)$.
	\item[(2)] $w$ is continuous in $B_1$.
	\item[(3)]  $\Delta w\geqslant0$ in $B_1$.
	\item[(4)]  $w\Delta w=0$ in the sense of Radon measure on $B_1$.
\end{proposition}

We now derive the explicit form of the blowup limit $w$.

\begin{proposition}
	Let $u$ be a weak solution of \eqref{1.20}  satisfying the assumptions in $\text{Theorem 1}$. Then at each point $x^0\in\Sigma^u$ there exists an integer $N_0\geqslant2$ such that $
H_{x^0,u^-}(0^+)=N_0$, and
	$$
	\frac{u^-(x^0+r_mx)}{\sqrt{r_m^{-1}\int_{\partial B_{r_m}(x^0)}\left(u^-\right)^2dS}}\rightarrow	w\left(\rho,\theta\right)
	$$
	strongly in $W^{1,2}_{\text{\rm loc}}\left(B_1\backslash\{O\}\right)$ and weakly in $W^{1,2}(B_1)$ as $r\rightarrow0^+$, where
	$$
	w\left(\rho,\theta\right)=\begin{cases}
		\begin{matrix}
			-\frac{\rho^{N_0}\left|\sin\left(N_0\,\theta\right)\right|}{\sqrt{\left|\int_{-\pi}^{0}\left(\sin\left(N_0\,\theta\right)\right)^2d\theta\right|}},&		\theta \in\left(-\pi,0\right),\\
		\end{matrix}\\
		\begin{matrix}
			0,&	\hspace{30.6mm}\text{otherwise},\\
		\end{matrix}\\
	\end{cases}
	$$
	and $x=\left(\rho\cos\theta,\rho\sin\theta\right)$.
\end{proposition}
\begin{pf}
	Let $r_m\rightarrow0^+$ be an arbitrary sequence such that $v_m$ defined in \eqref{4.8} converges weakly in $W^{1,2}(B_1)$ to $w$. Set $N_0=H_{x^0,u^-}(0^+)$ and it will be shown below
	 that $N_0$ is an integer. By $\text{Proposition 3.6}$ (3) and $\text{Proposition 3.7}$, we have $N_0\geqslant\frac{3}{2}$ and $\Delta w=0$ in $B_1\cap\{w<0\}$. Moreover,
\begin{equation}
	\label{3.18}
	w\left(\rho,\theta\right)=\begin{cases}
	\begin{matrix}
	-C\rho^{N_0}\left|\sin\left(N_0\theta+\varphi_-\right)\right|,&		\theta \in\left(\theta_1,\theta_2\right),\\
	\end{matrix}\\
	\begin{matrix}
	0,&	\hspace{38.3mm}\theta\notin\left(\theta_1,\theta_2\right),\\
	\end{matrix}\\
	\end{cases}
\end{equation}
	where $\theta_1<\theta_2$ and the condition $\int_{\partial B_1}w^2dS=1$ implies $$
	C=\frac{1}{\sqrt{\left|\int_{-\pi}^{0}\left(\sin\left(N_0\theta+\varphi_-\right)\right)^2d\theta\right|}}.$$
	
	We next claim that 
	$$\theta_1=-\pi\quad\text{and} \quad\theta_2=0.$$ Suppose towards a contradiction that $-\pi<\theta_1<0$. Without loss of generality, we take $y_0=\left(r_0\cos{\theta_1},r_0\sin{\theta_1}\right)$ and $r_m=\frac{1}{2}r_0\sin{\theta_1}$ such that $B_{r_m}(y^0)\cap\{x_2=0\}=\varnothing$. Since $u_\infty\equiv0$, we can deduce that
	\begin{equation}
	\label{4.11}
	\Delta u_m^-(B_{r_m}(y^0))\rightarrow0\quad \text{as}\,\, r_m\rightarrow0^+.
	\end{equation}
	However, the assumption on connected components implies that $\partial_{red}\{u\left(x^0+r_mx\right)<0\}$ contains a continuous curve whose image converges, as $r_m\rightarrow 0^+$, locally in $\{x_2<0\}$ to the half-line $\{\beta x\left|\,\beta>0,\,\,x_2<0\right.\}$. Consequently,
	$$
	\Delta u_m^-(B_{r_m}(y^0))=\int_{ B_{r_m}(y^0)\cap\partial\{u_m<0\}}\left|\nabla u_m^-\right|dS\geqslant\int_{ B_{r_m}(y^0)\cap\partial_{\text{red}}\{u_m<0\}}\sqrt{-x_2}dS>0,
	$$
which contradicts \eqref{4.11}. Hence, we can obtain $\theta_1=-\pi$ and $\theta_2=0$.
	
Since $w=0$ in $B_1\cap\{x_2\geqslant0\}$ and $w$ is continuous in $B_1$, it follows from \eqref{3.18} that $H_{x^0,u^-}(0^+)$ is an integer and $N_0\geqslant2$. Consequently, 
	$$
	w\left(\rho,\theta\right)=\begin{cases}
		\begin{matrix}
			-\frac{\rho^{N_0}\left|\sin\left(N_0\,\theta\right)\right|}{\sqrt{\left|\int_{-\pi}^{0}\left(\sin\left(N_0\,\theta\right)\right)^2d\theta\right|}},&		\theta \in\left(-\pi,0\right),\\
		\end{matrix}\\
		\begin{matrix}
			0,&	\hspace{30.6mm}\text{otherwise}.\\
		\end{matrix}\\
	\end{cases}
	$$
\end{pf}

\begin{proposition}
	Let $u$ be a weak solution of \eqref{1.20}  satisfying the assumptions in $\text{Theorem 1}$. Then $\Sigma^u$ is locally in $\Omega$ a finite set.
\end{proposition}
\begin{pf}
	Suppose this fails. There exists $x^m\in\Sigma^u$ converging to $x^0\in\Omega$, where $x^m\ne x^0$ for all $m$. Since $\Sigma^u$ is closed by $\text{Remark 3.1}$, we have $x^0\in\Sigma^u$. Let $r_m=2\left|x^m-x^0\right|$ and $y=\frac{x^m-x^0}{r_m}$. Using the fact \eqref{4.8} and observing that
	$$
	\frac{u_m^-\left(y+rx\right)}{r^\frac{3}{2}}=\frac{u^-(x^0+r_m(y+rx))}{r^\frac{3}{2}r_m^\frac{3}{2}}=\frac{u^-(x^m+r_mrx)}{(rr_m)^\frac{3}{2}}\rightarrow0\quad \text{as} \,\,m\rightarrow\infty,
	$$
	we conclude that $y\in\Sigma^{u_m}$ for any $m\in\mathbb{N}$.
	
	Now, since $H_{x^0,u^-}(r)\geqslant\frac{3}{2}$ and $V_{x^0,u^-}(r)\geqslant0$, it follows from $y\in\Sigma^{u_m}$ that
	\begin{equation}
	\nonumber
	r\int_{ B_r(y)}\left|\nabla v_m\right|^2dx\geqslant\frac{3}{2}\int_{\partial B_r(y)}\left(v_m\right)^2dS\quad \text{for\,\,all}\,\,r\in\left(0,\frac{1}{2}\right).
	\end{equation}
	By $\text{Proposition 3.7}$ and passing to the limit, we arrive at
	\begin{equation}
	\label{4.12}
	r\int_{ B_r(y)}\left|\nabla w\right|^2dx\geqslant\frac{3}{2}\int_{\partial B_r(y)}w^2dS\quad \text{for\,\,all}\,\,r\in\left(0,\frac{1}{4}\right).
	\end{equation}
	
	On the other hand, by Taylor formula, we can find that
	$$
	\underset{r\rightarrow0^+}{\lim}\frac{r\int_{ B_r(y)}\left|\nabla w\right|^2dx}{\int_{\partial B_r(y)}w^2dS}=	\underset{r\rightarrow0^+}{\lim}\frac{r\int_{\partial B_r(y)}w\left(\nabla w,r\nu\right) dS}{\int_{\partial B_r(y)}w^2dS}=1,
	$$
	contradicting \eqref{4.12}.
\end{pf}

Finally, we complete the proof that $M_{x^0,u^-}(0^+)\ne\frac{2}{3}$ in Case 1.1.

\begin{theorem}
	Let $u$ be a weak solution of \eqref{1.20}  satisfying the assumptions in $\text{Theorem 1}$. Then $\Sigma^u$ is empty.
\end{theorem}
\begin{pf}
	Suppose towards a contradiction that $\Sigma^u\ne\varnothing$ and there exists a point $x^0\in\Sigma^u$. By $\text{Proposition 3.8}$, as $r\rightarrow0^+$
	$$
	\frac{u^-(x^0+r_mx)}{\sqrt{r_m^{-1}\int_{\partial B_{r_m}(x^0)}\left(u^-\right)^2dS}}\rightarrow	w\left(\rho,\theta\right)=\begin{cases}
		\begin{matrix}
			-\frac{\rho^{N_0}\left|\sin\left(N_0\,\theta\right)\right|}{\sqrt{\left|\int_{-\pi}^{0}\left(\sin\left(N_0\,\theta\right)\right)^2d\theta\right|}},&		\theta \in\left(-\pi,0\right),\\
		\end{matrix}\\
		\begin{matrix}
			0,&	\hspace{30.6mm}\text{otherwise},\\
		\end{matrix}\\
	\end{cases}
	$$
	strongly in $W^{1,2}_{\text{loc}}\left(B_1\backslash\{O\}\right)$ and weakly in $W^{1,2}(B_1)$, where $x=\left(\rho\cos\theta,\rho\sin\theta\right)$.
	
	In fact, $N_0\geqslant2$ implies that there exists $\tilde{\theta}$ such that $w(x)=0$ on $\{\tilde{\theta} r_m\left|\,r_m>0\right.\}$. However, since $x^0\in\Sigma^u$, for any $\tilde{B}\subset \subset B_1\cap\{x_2<0\}$, we have $v_m<0$ in $\tilde{B}$ as $r_m$ sufficiently small. Thus, we can choose $\tilde{B}$ such that $\tilde{B}\cap\{\tilde{\theta} r_m\left|\,r_m>0\right.\}\ne\varnothing$. 
	
	Similar to the one in the proof of $\text{Theorem 3.3}$ and we obtain
	$$
	0\leftarrow\Delta u_m^-(\tilde{B})\geqslant\int_{ \tilde{B}\cap\partial_{red}\left\{u(x^0+r_mx)<0\right\}}\sqrt{-x_2}dS>0,
	$$
	which is a contradiction.
\end{pf}

\subsection{The possible singular asymptotics for Case 1.2}
\quad

In this subsection, we analyze the blowup limits for Case 1.2. Consider the rescaled function
$$
u_m=\frac{u\left(x^0+r_mx\right)}{r_m^{3/2}}.
$$
satisfying
\begin{equation}
	\label{4.19}
	C_1\left|x-x^0\right|\leqslant \left|\nabla u^+\right|^2\leqslant C_2\left|x-x^0\right|\quad \text{in}\,\,B_{r_m}\left(x^0\right)\cap\{u>0\},
\end{equation}	
for some $0<C_1<C_2$. Note that the blowup limit $u^-_\infty$ is non-trivial, i.e., $u_\infty^-\ne 0$. This follows from the estimate
\begin{equation}
	\nonumber
	\begin{aligned}
\left|\nabla u^-_m\right|^2=&\frac{1}{r_m}\left|\nabla u^-(r_mx+x^0)\right|^2\\
\geqslant&-\frac{1}{r_m}r_mx_2=-x_2>0\quad\quad \text{on}\,\,\partial\{u_m<0\}\cap\{x_2<0\},
	\end{aligned}
\end{equation}
which ensures that the frequency formula and $\text{Theorem 3.3}$ remain applicable in Case 1.2.

It is crucial that $u_\infty^+\not\equiv0$ in this case. Otherwise, the classification \eqref{4.19} yields 
$$
\left|\nabla u_m^+\right|^2=\left|\frac{\nabla u^+(x^0+r_mx)}{r^{1/2}}\right|^2\geqslant\frac{C_1r|x|}{r}=C_1|x|.
$$
Combined with the strong convergence of $u_m$ to $u_\infty$ in $W^{1,2}_{\text{loc}}\left(
\mathbb{R}^2\right)$, this implies
$$
\int_{A}\left|\nabla u_\infty^+\right|^2dx\leftarrow\int_{A}\left|\nabla u_m^+\right|^2dx\geqslant\int_{A\cap\{u_m>0\}}C_1|x|dx\geqslant C\rho>0,
$$
where $A=\left\{x\in\mathbb{R}^2\left|\,|x|\geqslant\rho>0\right.\right\}$, which contradicts $u_\infty^+=0$. 

 The remainder of the proof closely parallels that of Case 1.1, leading to the asymptotic $$u_\infty^-(x)=U_\infty^-(\rho,\theta)=C_-\rho^{\frac{3}{2}}\cos\left(\frac{3}{2}\theta+\varphi_-\right)\quad\text{in} \,\,\{u_\infty<0\}.$$
Similarly, by $\text{Proposition 2.4}$, we also have
$$u_\infty^+(x)=U_\infty^+(\rho,\theta)=C_+\rho^{\frac{3}{2}}\cos\left(\frac{3}{2}\theta+\varphi_+\right)\quad\text{in} \,\,\{u_\infty>0\}.$$
To complete the proof of Case 1.2 in $\text{Theorem 1}$, we proceed to determine $u_\infty$ explicitly.

First, we derive the free boundary condition for $u_\infty$. Following the argument in Case 1.1, let $z=\left(\rho\cos\theta_1,\rho\sin\theta_1\right)\in\partial\{u_\infty<0\}\backslash \{O\}$ and choose $\delta>0$ small enough such that $B_\delta(z)\subset\{x_2<0\}$. Since the normal to $\partial\{u_\infty<0\}$ has the constant value $\nu(z)$ on $\partial\{u_\infty<0\}\cap B_\delta(z)$, we substitute the test function $\phi(x)=\eta(x)\nu(z)$ with an arbitrary $\eta\in C_0^1(B_\delta(z),\mathbb{R}^2)$ into \eqref{3.7} and integrating by parts. This leads to the following two possibilities.
	
(1) $z\in\partial\{u_\infty<0\}\cap\partial\{u_\infty>0\}$.
\begin{equation}
	\nonumber
	\begin{aligned}
		0=&\int_{ B_\delta(z)}\left(\left|\nabla u_\infty\right|^2\mathrm{div}\phi-2\nabla u_\infty D\phi\nabla u_\infty-x_2\chi_-\mathrm{div}\phi-\chi_-\phi_2\right)dx\\
		=&\int_{ B_\delta(z)\cap\{u_\infty\geqslant 0\}}\left(\mathrm{div}\left(\left|\nabla u_\infty^+\right|^2\phi-2\left(\phi,\nabla u_\infty^+\right)\nabla u_\infty^+\right)+\mathrm{div}\left(-x_2\phi\right)\chi_-\right)dx\\
		&+\int_{ B_\delta(z)\cap\{u_\infty< 0\}}\left(\mathrm{div}\left(\left|\nabla u_\infty^-\right|^2\phi-2\left(\phi,\nabla u_\infty^-\right)\nabla u_\infty^-\right)+\mathrm{div}\left(-x_2\phi\right)\right)dx\\
		=&\int_{ B_\delta(z)\cap\partial\{u_\infty<0\}\cap\partial\{u_\infty>0\}}\left(-\left(\left|\nabla u_\infty^-\right|^2-\left|\nabla u_\infty^+\right|^2\right)-x_2\left(1-\chi_-\right)\right)\left(\eta,\nu\right) dS.
	\end{aligned}
\end{equation}
 Note that $\chi_-$ is the strong $L_{\text{loc}}^1$ limit of $\chi_{\{u_m<0\}}$ and $\chi_-\in\{0,1\}$. If $\chi_-=1$, we have 
$$
\left|\nabla u_\infty^-\right|^2-\left|\nabla u_\infty^+\right|^2=0\quad\text{on}\,\,\partial\{u_\infty<0\}\cap\partial\{u_\infty>0\}\cap B_\delta(y).
$$
In polar coordinates, this becomes
$$
\frac{9}{4}C_+^2\rho=\left|\nabla u_\infty^+\right|^2=\left|\nabla u_\infty^-\right|^2=\frac{9}{4}C_-^2\rho,
$$
 which implies that $C_+=C_-$. On the other hand, by
 \begin{equation}
\nonumber
 	\lim\limits_{\tau\rightarrow\theta_1+}\partial_\theta\left(1,\tau\right)=-\lim\limits_{\tau\rightarrow\theta_1-}\partial_\theta\left(1,\tau\right).
 \end{equation}
 we obtain that
\begin{equation}
	\label{3.19}
	\sin\left(\frac{3}{2}\theta_1+\varphi_+\right)+\sin\left(\frac{3}{2}\theta_1+\varphi_-\right)=2\sin\left(\frac{3}{2}\theta_1+\frac{\varphi_++\varphi_-}{2}\right)\cos\left(\frac{\varphi_+-\varphi_-}{2}\right)=0.
\end{equation}
However, since $U_\infty^-(\rho,\theta_1)=U_\infty^+(\rho,\theta_1)=0$, we have (see Fig. 11 (1))
$$
\frac{3}{2}\theta_1+\varphi_+=\frac{\pi}{2}+2k_1\pi,\quad\frac{3}{2}\theta_1+\varphi_-=-\frac{3\pi}{2}+2k_2\pi,
$$
or (see Fig. 11 (2))
$$
\frac{3}{2}\theta_1+\varphi_+=-\frac{\pi}{2}+2k_3\pi,\quad\frac{3}{2}\theta_1+\varphi_-=-\frac{\pi}{2}+2k_4\pi,
$$
where $k_i\in\mathbb{Z}$ for $i\in\{1,2,3,4\}$. 
	\begin{figure}[!h]
	\includegraphics[width=120mm]{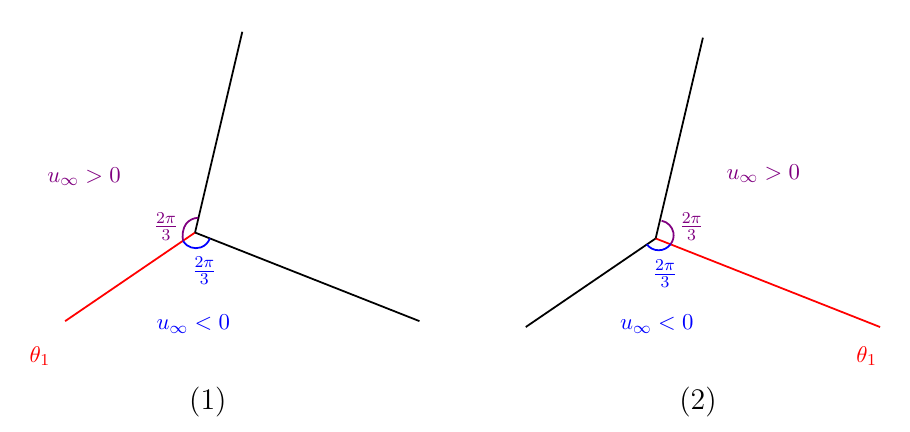}
	\caption{The possible scenarios}
\end{figure}

Consequently,  $$\sin\left(\frac{3}{2}\theta_1+\frac{\varphi_++\varphi_-}{2}\right)\ne0\quad\text{and} \quad \cos\left(\frac{\varphi_+-\varphi_-}{2}\right)\ne0,$$ which contradicts \eqref{3.19}. Therefore, $\chi_-=0$ and
\begin{equation}
\label{fb1}
\left|\nabla u_\infty^-\right|^2-\left|\nabla u_\infty^+\right|^2=-x_2\quad\text{on}\,\,\partial\{u_\infty<0\}\cap\partial\{u_\infty>0\}.
\end{equation}

(2) If $z\in\partial\{u_\infty<0\}\backslash\partial\{u_\infty>0\}$, we take $\delta$ small enough such that $u_\infty\leqslant0$ in $B_\delta(z)$. Similarly, if $z\in\partial\{u_\infty>0\}\backslash\partial\{u_\infty<0\}$, we choose  $\delta$ small enough such that $u_\infty\geqslant0$ in $B_\delta(z)$. Hence, argue similarly to the proof of Case 1.1, we obtain
\begin{equation}
\label{fb2}
\left|\nabla u^-_\infty\right|^2=-x_2\quad\text{on}\,\,\partial\{u_\infty<0\}\backslash\partial\{u_\infty>0\},
\end{equation}
and
\begin{equation}
\label{fb3}
\left|\nabla u^+_\infty\right|^2=x_2\quad\text{on}\,\,\partial\{u_\infty>0\}\backslash\partial\{u_\infty<0\}.
\end{equation}

Note that each connected component of $\{u_\infty<0\}$ and $\{u_\infty>0\}$ is a cone with vertex at the origin and opening angle $\frac{2\pi}{3}$, and that $\{u_\infty>0\}$ contains at most two connected components. We now examine the possible singular asymptotics configurations.

(1) Suppose that $\{u_\infty>0\}$ has only one connected component.

(1a) If $\partial\{u_\infty>0\}$ and  $\partial\{u_\infty<0\}$ only touch at one point (that is $\partial\{u_\infty>0\}\cap\partial\{u_\infty<0\}=\{O\}$), then by \eqref{fb2} and \eqref{fb3}, the solution $u_\infty$ takes the form
	$$
U_{\infty}(\rho,\theta)=\begin{cases}
	\begin{matrix}
		\frac{\sqrt{2}}{3}\rho^{3/2}\cos\left(\frac{3}{2}\theta-\frac{3\pi}{4}\right),&	\text{in}\,\,\left(\frac{\pi}{6},\frac{5\pi}{6}\right),	\\
	\end{matrix}\\
	\begin{matrix}
		\frac{\sqrt{2}}{3}\rho^{3/2}\cos\left(\frac{3}{2}\theta-\frac{\pi}{4}\right),&	\hspace{1.4mm}	\text{in}\,\,\left(-\frac{5\pi}{6},-\frac{\pi}{6}\right),	\\
	\end{matrix}\\
	\begin{matrix}
		0,&	\hspace{33.4mm}\text{otherwise,}	\\
	\end{matrix}\\
\end{cases}
$$
with the density $M_{x^0,u^-}(0^+)=\frac{\sqrt{3}}{3}$, where $u_\infty(x)=U_\infty(\rho,\theta)$ in polar coordinates.

(1b) If $\partial\{u_\infty>0\}$ and  $\partial\{u_\infty<0\}$ share a common free boundary segment (see Fig. 12), namely $\partial\{u_\infty<0\}\backslash\partial\{u_\infty>0\}\ne\varnothing$, then the free boundary conditions are
\begin{equation}
\label{3.24}
	\begin{cases}
		\begin{matrix}
			\left|\nabla u^-_\infty\right|^2-\left|\nabla  u^+_\infty\right|^2=-x_2,& \text{on 
			}\,\,\partial\{u_\infty>0\}\cap\partial\{u_\infty<0\},\\
		\end{matrix}\\
		\begin{matrix}
			\left|\nabla u^-_\infty\right|^2=-x_2,&	\hspace{17.2mm} 	 \text{on 
			}\,\,\partial\{u_\infty<0\}\backslash\partial\{u_\infty>0\},\\
		\end{matrix}\\
		\begin{matrix}
			\left|\nabla u^+_\infty\right|^2=x_2,&	\hspace{20.3mm} 	 \text{on 
			}\,\,\partial\{u_\infty>0\}\backslash\partial\{u_\infty<0\}.\\
		\end{matrix}\\
	\end{cases}
\end{equation}
	\begin{figure}[!h]
	\includegraphics[width=100mm]{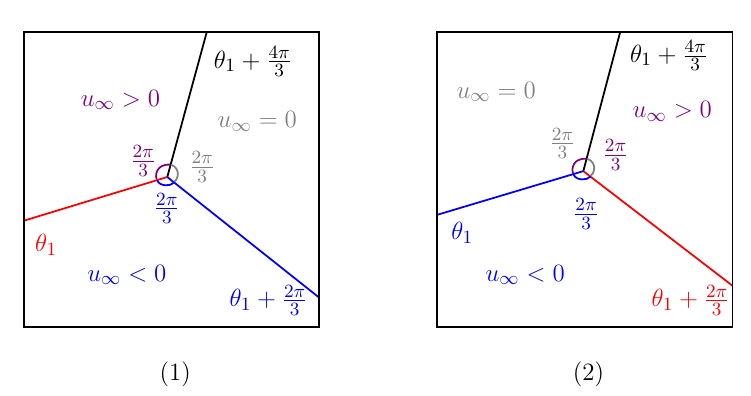}
	\caption{$\partial\{u_\infty>0\}\backslash\partial\{u_\infty<0\}\ne\varnothing$}
\end{figure}

Suppose that $\{u_\infty<0\}=\{(\rho,\theta)\left|\, \theta_1<\theta<\theta_1+\frac{2\pi}{3}\right.\}$, where $-\pi\leqslant\theta_1\leqslant-\frac{2\pi}{3}$. Substituting the explicit forms
 $$U_\infty^-(\rho,\theta)=C_-\rho^{\frac{3}{2}}\cos\left(\frac{3}{2}\theta+\varphi_-\right)\quad\text{in}\,\,\{u_\infty<0\},$$ 
$$U_\infty^-(\rho,\theta)=C_+\rho^{\frac{3}{2}}\cos\left(\frac{3}{2}\theta+\varphi_+\right)\quad\text{in}\,\, \{u_\infty>0\},$$ 
into \eqref{3.24} yields, for the configuration in Fig. 12 (1),
\begin{equation}
	\nonumber
	\begin{cases}
		\begin{matrix}
		\frac{9}{4}C_-^2-\frac{9}{4}C_+^2=-\sin\theta_1,\\
		\end{matrix}\\
		\begin{matrix}
			\frac{9}{4}C_-^2=-\sin\left(\theta_1+\frac{2\pi}{3}\right),\\
		\end{matrix}\\
		\begin{matrix}
			\frac{9}{4}C_+^2=\sin\left(\theta_1+\frac{4\pi}{3}\right),\\
		\end{matrix}\\
	\end{cases}
\end{equation}
or for Fig. 12 (2),
\begin{equation}
	\nonumber
	\begin{cases}
		\begin{matrix}
			\frac{9}{4}C_-^2-\frac{9}{4}C_+^2=-\sin\left(\theta_1+\frac{2\pi}{3}\right),\\
		\end{matrix}\\
		\begin{matrix}
			\frac{9}{4}C_-^2=-\sin\theta_1,\\
		\end{matrix}\\
		\begin{matrix}
			\frac{9}{4}C_+^2=\sin\left(\theta_1+\frac{4\pi}{3}\right),\\
		\end{matrix}\\
	\end{cases}
\end{equation}
which implies $\theta_1=-\pi$ or $\theta_1=-\frac{2\pi}{3}$.

Thus, $u_\infty$ is given by
	$$
U_{\infty}(\rho,\theta)=\begin{cases}
	\begin{matrix}
		\frac{\sqrt{2\sqrt{3}}}{3}\rho^{3/2}\cos\left(\frac{3}{2}\theta-\pi\right),&	\text{in}\,\,\left(\frac{\pi}{3},\pi\right),	\\
	\end{matrix}\\
	\begin{matrix}
		\frac{\sqrt{2\sqrt{3}}}{3}\rho^{3/2}\cos\left(\frac{3}{2}\theta\right),&	\hspace{6.9mm}	\text{in}\,\,\left(-\pi,-\frac{\pi}{3}\right),	\\
	\end{matrix}\\
	\begin{matrix}
		0,&	\hspace{36mm}\text{otherwise,}	\\
	\end{matrix}\\
\end{cases}
$$
or
$$
U_{\infty}(\rho,\theta)=\begin{cases}
	\begin{matrix}
		\frac{\sqrt{2\sqrt{3}}}{3}\rho^{3/2}\cos\left(\frac{3}{2}\theta-\frac{\pi}{2}\right),&		\text{in}\,\,\left(0,\frac{2\pi}{3}\right),	\\
	\end{matrix}\\
	\begin{matrix}
		\frac{\sqrt{2\sqrt{3}}}{3}\rho^{3/2}\cos\left(\frac{3}{2}\theta-\frac{\pi}{2}\right),&		\text{in}\,\,\left(-\frac{2\pi}{3},0\right),	\\
	\end{matrix}\\
	\begin{matrix}
		0,&	\hspace{36.2mm}\text{otherwise,}	\\
	\end{matrix}\\
\end{cases}
$$
with corresponding density 
	\begin{equation}
	\nonumber
	\begin{aligned}
		M_{x^0,u^-}(0^+)&=\int_{ B_1}-x_2\chi_{\left\{\left(\rho\cos\theta,\rho\sin\theta\right)\left|\, -\pi<\theta<-\frac{\pi}{3}\right.\right\}}dx\\
		&=\int_{0}^{1}\int_{-\pi}^{-\frac{\pi}{3}}-\rho^2\sin\theta  d\theta d\rho=\frac{1}{2},
	\end{aligned}
\end{equation}
or
	\begin{equation}
	\nonumber
	\begin{aligned}
		M_{x^0,u^-}(0^+)&=\int_{ B_1}-x_2\chi_{\left\{\left(\rho\cos\theta,\rho\sin\theta\right)\left|\, -\frac{2\pi}{3}<\theta<0\right.\right\}}dx\\
		&=\int_{0}^{1}\int_{-\frac{2\pi}{3}}^{0}-\rho^2\sin\theta d\theta d\rho=\frac{1}{2}.
	\end{aligned}
\end{equation}

(2) Suppose that $\{u_\infty>0\}$ consists two connected components and let the negative phase be given by $\{u_\infty<0\}=\left\{\left(\rho,\theta\right)\left|\,\theta_1<\theta<\theta_1+\frac{2\pi}{3}\right.\right\}$, where $-\pi\leqslant\theta_1\leqslant-\frac{2\pi}{3}$. Then the solution has the form 
\begin{equation}
\nonumber
\begin{cases}
	\begin{matrix}
		U_\infty^-=C_1\rho^{3/2}\cos\left(\frac{3}{2}\theta+\varphi_1\right),&		\text{in}\,\,\left[\theta_1,\theta_1+\frac{2\pi}{3}\right),\\
	\end{matrix}\\
	\begin{matrix}
		U_\infty^+=C_2\rho^{3/2}\cos\left(\frac{3}{2}\theta+\varphi_2\right),&		\text{in}\,\,\left[\theta_1+\frac{2\pi}{3},\theta_1+\frac{4\pi}{3}\right),\\
	\end{matrix}\\
	\begin{matrix}
		U_\infty^+=C_3\rho^{3/2}\cos\left(\frac{3}{2}\theta+\varphi_3\right),&			\text{in}\,\,\left[\theta_1+\frac{4\pi}{3},\theta_1+2\pi\right).	\\
	\end{matrix}\\
\end{cases}\quad\text{(see Fig. 13)}
\end{equation}
	\begin{figure}[!h]
	\includegraphics[width=50mm]{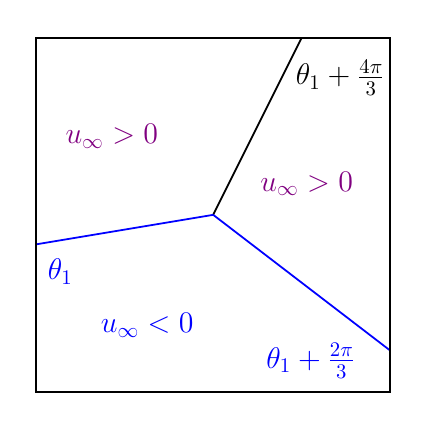}
	\caption{$\partial\{u_\infty>0\}\backslash\partial\{u_\infty<0\}=\varnothing$}
\end{figure}

 We argue similarly to (1) that
\begin{equation}
	\nonumber
	\begin{cases}
		\begin{matrix}
			\frac{9}{4}C_1^2-\frac{9}{4}C_3^2=-\sin\theta_1,\\
		\end{matrix}\\
		\begin{matrix}
			\frac{9}{4}C_1^2-\frac{9}{4}C_2^2=-\sin\left(\theta_1+\frac{2\pi}{3}\right),\\
		\end{matrix}\\
		\begin{matrix}
			\frac{9}{4}C_2^2=\frac{9}{4}C_3^2=\sin\left(\theta_1+\frac{4\pi}{3}\right),\\
		\end{matrix}\\
	\end{cases}
\end{equation}
which implies $\theta_1=-\frac{5\pi}{6}$. Hence
$$
U_{\infty}(\rho,\theta)=\begin{cases}
	\begin{matrix}
		\frac{2}{3}\rho^{3/2}\cos\left(\frac{3}{2}\theta+\frac{3}{4}\pi\right),&	\hspace{0.5mm}	\text{in}\,\,\left[\frac{\pi}{2},\frac{7\pi}{6}\right),	\\
	\end{matrix}\\
	\begin{matrix}
		\frac{2}{3}\rho^{3/2}\cos\left(\frac{3}{2}\theta-\frac{\pi}{4}\right),&	\hspace{2.4mm}	\text{in}\,\,\left[-\frac{\pi}{6},\frac{\pi}{2}\right),	\\
	\end{matrix}\\
	\begin{matrix}
		\frac{\sqrt{6}}{3}\rho^{3/2}\cos\left(\frac{3}{2}\theta-\frac{\pi}{4}\right),&\text{in}\,\,\left[-\frac{5\pi}{6},-\frac{\pi}{6}\right),	\\
	\end{matrix}\\
\end{cases}
$$
with the density 
$$
M_{x^0,u^-}(0^+)=\frac{\sqrt{3}}{3}.
$$

This completes the proof of Case 1.2 in $\text{Theorem 1}$.

\subsection{The possible singular asymptotics for Case 1.3}
\quad

 According to $\text{Proposition 2.1}$, we analyze $u^-$ by following the derivation in Subsection 3.2, which yields the following result.

\begin{proposition}
Let $u$ be a weak solution of \eqref{1.20} satisfying the assumptions in $\text{Theorem 1}$ and consider the blowup sequence 
	$$
u_m=\frac{u\left(x^0+r_mx\right)}{r_m^{\alpha}},\quad1<\alpha<3/2,
$$
as defied by $\text{Definition C}$ in $\text{Appendix A}$. Then $u_m\rightarrow u_\infty\equiv0$ strongly in $W^{1,2}_{\text{\rm loc}}(\mathbb{R}^2)$ and locally uniformly in $\mathbb{R}^2$ as $m\rightarrow\infty$. 
\end{proposition}

\begin{pf}
We argue by contradiction. Suppose that $\{u_\infty<0\}\ne\varnothing$ (the case $\{u_\infty>0\}\ne\varnothing$ can be handled by a similar argument, the detail of which we omit here). By the homogeneity of $u_\infty$ stated in $\text{Proposition 2.5}$, every connected component of $\{u_\infty>0\}$ and $\{u_\infty<0\}$ is a corner region with opening angle $\pi/\alpha$, where $1<\alpha<\frac{3}{2}$. From the previous two subsections, if $\{u_\infty<0\}$ is not empty, it consists of exactly one connected component. Hence, in polar coordinates, 
	\begin{equation}
		\nonumber
		u_{\infty}^-(x)=U_{\infty,\alpha}^-(\rho,\theta)=	\begin{cases}
			\begin{matrix}
			C\rho^\alpha\cos\left(\alpha\theta+\varphi\right),&	\text{in}\,\,\left(\theta_1,\theta_1+\frac{\pi}{\alpha}\right),	\\
			\end{matrix}\\
			\begin{matrix}
		0,&\hspace{27.2mm}	\text{otherwise,}	\\
			\end{matrix}\\
		\end{cases}
	\end{equation}
for some $C>0$, $\varphi\in\left[-\pi,\pi\right)$ and $\theta_1\in\left[-\pi,-\frac{\pi}{\alpha}\right]$.
	
Now, following the proof of $\text{Proposition 3.1}$, take $\phi\in C_0^1\left(\mathbb{R}^2,\mathbb{R}^2\right)$. From \eqref{def of weak} we obtain that
	\begin{equation}
	\nonumber
		\begin{aligned}
			0&=r_m^{2\alpha-1}\int_{ B_1}\left(\left|\nabla u_m\right|^2\mathrm{div}\phi-2\nabla u_mD\phi_m\nabla u_m\right)dx\\
			&+r_m^2\int_{ B_1\cap\{x_2<0\}}\left(-x_2\chi_{\{u_m<0\}}\mathrm{div}\phi-\phi_2\chi_{\{u_m<0\}}\right)dx\\
			&+r_m^2\int_{ B_1\cap\{x_2>0\}}\left(x_2\chi_{\{u_m>0\}}\mathrm{div}\phi+\phi_2\chi_{\{u_m>0\}}\right)dx.
		\end{aligned}
	\end{equation}
Dividing by $r_m^{2\alpha-1}$ gives 
	\begin{equation}
	\nonumber
	\begin{aligned}
		0&=\int_{ B_1}\left(\left|\nabla u_m\right|^2\mathrm{div}\phi-2\nabla u_mD\phi_m\nabla u_m\right)dx\\
		&+r_m^{3-2\alpha}\int_{ B_1\cap\{x_2<0\}}\left(-x_2\chi_{\{u_m<0\}}\mathrm{div}\phi-\phi_2\chi_{\{u_m<0\}}\right)dx\\
		&+r_m^{3-2\alpha}\int_{ B_1\cap\{x_2>0\}}\left(x_2\chi_{\{u_m>0\}}\mathrm{div}\phi+\phi_2\chi_{\{u_m>0\}}\right)dx.
	\end{aligned}
\end{equation}
Since $1<\alpha<\frac{3}{2}$, the exponent $\left(3-2\alpha\right)$ is positive. Moreover 
 $$\left|\int_{ B_1\cap\{x_2<0\}}\left(-x_2\chi_{\{u_m<0\}}\mathrm{div}\phi-\phi_2\chi_{\{u_m<0\}}\right)dx+\int_{ B_1\cap\{x_2>0\}}\left(x_2\chi_{\{u_m>0\}}\mathrm{div}\phi+\phi_2\chi_{\{u_m>0\}}\right)dx\right|<\infty.$$Consequently, passing to the limit $m\rightarrow\infty$ yields
$$
\int_{ \mathbb{R}^2}\left(\left|\nabla u_\infty\right|^2\mathrm{div}\phi-2\nabla u_\infty D\phi\nabla u_\infty\right)dx=0.
$$

Now, we claim that there exists a one-phase free boundary point $z\in\Gamma_{\text{\rm op},\infty}^-\backslash\{O\}$, where $\Gamma_{\text{\rm op},\infty}^-=\partial\{u_\infty<0\}\backslash\partial\{u_\infty>0\}$. Indeed, since $1<\alpha<3/2$, every connected component of $\{u_\infty>0\}$ and $\{u_\infty<0\}$ is a corner region with opening angle $\pi/\alpha$, which lies strictly between $2\pi/3$ and $\pi$. The hypothesis that $\{u_\infty<0\}$ has exactly one component implies that $\{u_\infty>0\}$ has at most one component. 

(1) If $\{u_\infty>0\}=\varnothing$ (see Fig.14 (1)), then the entire free boundary of $\partial\{u_\infty<0\}$ coincides with $\Gamma_{\text{\rm op},\infty}^-$. Hence, the claim holds trivially.

(2) If $\{u_\infty>0\}\ne\varnothing$ (see Fig.14 (2) and (3)), the contact set $\partial\{u_\infty<0\}\cap\partial\{u_\infty>0\}$ can occupy at most a portion of $\partial\{u_\infty<0\}$. Consequently, the remaining part of $\partial\{u_\infty<0\}$ consists of one-phase free boundary points, which again proves the claim.

	\begin{figure}[!h]
	\includegraphics[width=140mm]{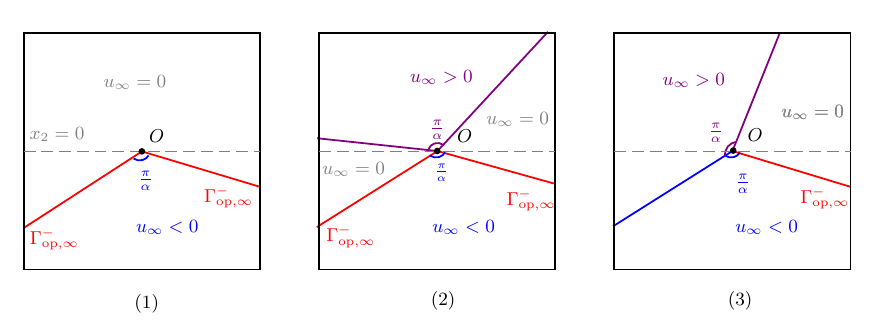}
	\caption{$\Gamma_{\text{\rm op},\infty}^-$}
\end{figure}

Take such a point $z$ and choose $\delta>0$ small enough such that $B_\delta(z)\subset\{x_2<0\}$. Argue similarly to the proof of Case 1.1, we obtain
$$
\left|\nabla u_\infty^-\right|=0,\quad\text{on}\,\,\Gamma_{\text{\rm op},\infty}^-\cap B_\delta(z).
$$
However, it follows from Hopf's lemma that $(\nabla u_\infty^-,\nu)\ne0$ on $\Gamma_{\text{\rm op},\infty}^-\cap B_\delta(z)$. This is a contradiction.
\end{pf}

We can still derive a result analogous to $\text{Proposition 3.1}$ in Case 1.3.

	\begin{proposition}
Let $u_m$ be the blowup sequence defined in $\text{Proposition 2.5}$. Then $\chi_-$ is a constant in the interior of $\{u_\infty=0\}\cap\{x_2\leqslant0\}$ in Case 1.3, where $\chi_-$ is the strong $L^1_{\text{\rm loc}}$ limit of $\chi_{\{u_m<0\}}$. 
\end{proposition}

\begin{pf}
For Case 1.3, let $\phi\in C_0^1(\mathbb{R}^2\cap\{x_2\leqslant0\},\mathbb{R}^2)$. From the definition of weak solution \eqref{def of weak}, we rewrite
	\begin{equation}
		\nonumber
		\begin{aligned}
			0=&r_m^{2\alpha-3}\int_{ B_1}\left(\left|\nabla u_m\right|^2\mathrm{div}\phi-2\nabla u_mD\phi\nabla u_m\right)dx\\
			&+\int_{ B_1\cap\{x_2<0\}}\left(-x_2\chi_{\{u_m<0\}}\mathrm{div}\phi-\phi_2\chi_{\{u_m<0\}}\right)dx\\
			=&A_m+B_m,
		\end{aligned}
	\end{equation}
	where
	$$
	A_m=r_m^{2\alpha-3}\int_{ B_1}\left(\left|\nabla u_m\right|^2\mathrm{div}\phi-2\nabla u_mD\phi\nabla u_m\right)dx,
	$$
	and
	$$
	B_m=\int_{ B_1\cap\{x_2<0\}}\left(-x_2\chi_{\{u_m<0\}}\mathrm{div}\phi-\phi_2\chi_{\{u_m<0\}}\right)dx.
	$$
Clearly, $|B_m|$ is bounded. Moreover, by the compact embedding of BV into $L^1$, $\chi_{\{u_m<0\}}$ converges strongly in $L^1_{\text{\rm loc}}(\mathbb{R}^2)$, which we denote by $\chi_-$. Consequently, the limit of $B_m$ exists and is finite, i.e., 
	\begin{equation}
	\nonumber
	\lim\limits_{m\rightarrow\infty}\left|B_m\right|=\left|\int_{ B_1\cap\{x_2<0\}}\left(-x_2\chi_-\mathrm{div}\phi-\phi_2\chi_-\right)dx
\right|=C,\quad\text{for some}\,\,C\geqslant0.
\end{equation}

We claim that $C=0$. Then, we obtain
$$
0= \int_{ \mathbb{R}^2\cap\{x_2<0\}}\left(-x_2\chi_-\mathrm{div}\phi-\phi_2\chi_-\right)dx.
$$
By repeating the argument from Case 1.1 and 1.2, we conclude that $\chi_-$ is a constant in the interior of $\{u_\infty=0\}\cap\{x_2\leqslant0\}$.

Assume, for contradiction, that $C\ne0$. Consider a connected open set $U$ contained in the interior of $\{u_\infty=0\}\cap\{x_2\leqslant0\}$. Since $u_m$ converges uniformly to $0$ in $\mathbb{R}^2$, for some $d_0>0$ and sufficiently large $m$, there exists $V=\{x\in U\left|\, u_m<-\delta\right.\}$ such that $\text{dist}\left(V,\{u_m=0\}\right)\geqslant d_0$ for $\delta>0$ small enough. Now consider a horizontal line segment $L$ connecting $V$ to $\{u_m=0\}$. We then derive that the tangential derivative of $u_m^-$ along $L$ is strictly positive. By the Cauchy-Schwarz inequality, we obtain that
	$$
	\int_L\left|\nabla u_m^-\right|^2d\ell\geqslant \tilde{C}>0, \quad \text{for some}\,\,\tilde{C}>0,
	$$
	where $d\ell$ is the arc-length element along $L$.
	
	Since the choice of $L$ is arbitrary, we conclude that
	$$
	0=\int_U\left|\nabla u_\infty^-\right|^2dx\leftarrow\int_U \left|\nabla u_m^-\right|^2dx\geqslant\tilde{C}>0,
	$$
which is a contradiction.
\end{pf}
	
\subsection{The possible singular profiles for Case 1}
\quad

To characterize the singular profiles formed by electric field under different configurations as $x^{\text{\rm st}}=x^0=(x_1^0,x_2^0)$, we establish that
\begin{itemize}
	\item [(1)] Singular profiles with density $M_{x^0,u^-}(0^+)=\frac{\sqrt{3}}{3}$ is Stokes corner.
	
	\item [(2)] Singular profiles with density $M_{x^0,u^-}(0^+)=\frac{1}{2}$ is asymmetric Stokes corner.
	
		\item [(3)] Singular profiles with density $M_{x^0,u^-,\alpha}(0^+)=0$ is cusp.
\end{itemize} 

To be more specific,

$\bullet$ Stokes corner

Firstly, we provide the proof for  $M_{x^0,u^-}(0^+)=\frac{\sqrt{3}}{3}$ as a representative example. The proof of $\text{Theorem 2}$ will be split into three steps. The case $M_{x^0,u^-}(0^+)=\frac{1}{2}$ follows by analogous arguments with minor modifications.

{\bf Step 1.} Assume that $\partial\{u<0\}$ is given by a continuous injective curve 
$$\sigma(t)=\left(\sigma_1(t),\sigma_2(t)\right),\quad t\in(-t_0,t_0)\quad\text{with}\,\, t_0>0 \quad\text{and} \quad\sigma(0)=\left(x_1^0,x_2^0\right).$$
Define the sets of possible tangential directions 
$$
\mathfrak{L}^\pm=\left\{\theta^*\in\left[-\pi,0\right]\left|\,\text{there is}\,\, t_m\rightarrow0^\pm\,\,\text{such that}\,\,\arg\left(\sigma(t_m)-\sigma(0)\right)\rightarrow\theta^*\,\,\text{as}\,\,m\rightarrow\infty\right. \right\}.
$$
We claim that $$\mathfrak{L}^+\,\,\text{and}\,\,\mathfrak{L}^-\,\,\text{are subsets of}\,\,\left\{-\pi,-\frac{5\pi}{6},-\frac{\pi}{6},0\right\}.$$

Suppose the claim is false. Then there exist a sequence $t_m\rightarrow0$ and a direction $\theta^*\in\left(\mathfrak{L}^+\cup\,\mathfrak{L}^-\backslash\left\{-\pi,-\frac{5\pi}{6},-\frac{\pi}{6},0\right\}\right)$ such that $\arg\left(\sigma(t_m)-\sigma(0)\right)\rightarrow\theta^*$ as $m\rightarrow\infty$. Let $$r_m=\left|\sigma(t_m)-\sigma(0)\right|\quad\text{and}\quad u_m=\frac{u(x^0+r_mx)}{r_m^{\frac{3}{2}}}.$$
Choose $\rho>0$ small enough such that $B^*=B_\rho\left(\cos\theta^*,\sin\theta^*\right)$ (see Fig. 15) satisfies
$$
B^*\cap\left(\left\{(x,0)\left|\,x\in\mathbb{R}\right.\right\}\cap\left\{\left(x,-\frac{\sqrt{3}}{3}\,|x|\right)\left|\,x\in\mathbb{R}\right.\right\}\right)=\varnothing.
$$
	\begin{figure}[!h]
	\includegraphics[width=100mm]{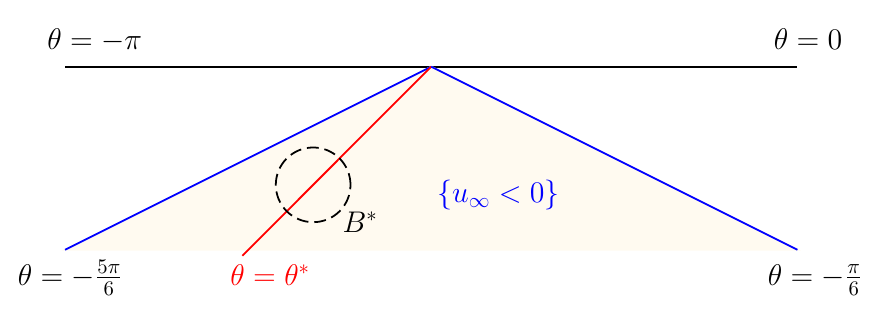}
	\caption{$B^*$}
\end{figure}

By $u_m\rightarrow u_\infty$ strongly in $W^{1,2}_{\text{\rm loc}}(\mathbb{R}^2)$, we have
$$
\Delta u_m^-(B^*)\rightarrow\Delta u_\infty^-(B^*)=0.
$$
However, by the assumption $u\geqslant0$ in $\{x_2\geqslant x_2^0\}$
\begin{equation}
	\label{3.25}
\Delta u_m^-(B^*)=\int_{ B^*\cap\partial\{u_m<0\}}\left|\nabla u_m^-\right|dS\geqslant\int_{ B^*\cap\partial\{u_m<0\}\cap\partial\{u_m>0\}}-x_2dS>0,
\end{equation}
which is a contradiction.

{\bf Step 2}. $\mathfrak{L^+}$ and $\mathfrak{L^-}$ are connected sets.

This is equivalent to saying that if there exist $\theta_1$, $\theta_2\in\mathfrak{L}^+$ with $\theta_1<\theta_2$ and ${t_m^{(1)}}$, ${t_m^{(2)}}$ such that $\arg\sigma\left(t_m^{(1)}\right)\rightarrow\theta_1$ and $\arg\sigma\left(t_m^{(2)}\right)\rightarrow\theta_2$ as $t_m^{(1)}\rightarrow0$, $t_m^{(2)}\rightarrow0$, we claim $\theta_3\in\mathfrak{L^+}$ for any $\theta_3\in\left(\theta_1,\theta_2\right)$.

Indeed, by the continuity of the free boundary and the intermediate value theorem, we conclude that there exists $t_m^{(3)}$ with $t_m^{(1)}<t_m^{(3)}<t_m^{(2)}$ such that $\arg\sigma\left(t_m^{(3)}\right)\rightarrow\theta_3$. Similarly, we obtain $\mathfrak{L^+}$ and $\mathfrak{L^-}$ are connected sets. 

 Moreover, define
$$
l^+=\lim\limits_{t\rightarrow0^+}\arg\left(\sigma(t)-\sigma(0)\right)\quad\text{and}\quad
l^-=\lim\limits_{t\rightarrow0-}\arg\left(\sigma(t)-\sigma(0)\right),
$$
where the limits $l^+$ and $l^-$ exist and take values in the set $\left\{-\pi,-\frac{5\pi}{6},-\frac{\pi}{6},0\right\}$.

{\bf Step 3}. We determine the values of $l^+$ and $l^-$ for each configuration.

 $\bullet$ Stokes corner.
 
Note that $M_{x^0,u^-}(0^+)=\int_{ B_1}(-x_2)\chi_{\left\{\left(\rho\cos\theta,\rho\sin\theta\right)\left|\, -\frac{5\pi}{6}<\theta<-\frac{\pi}{6}\right.\right\}}dx$ in Case 1.1, which yields that $\Delta u_\infty^-\left(B_{1/10}\left(\sqrt{3}/2,-1/2\right)\right)>0$ and $\Delta u_\infty^-\left(B_{1/10}\left(-\sqrt{3}/2,-1/2\right)\right)>0$. It follows that $\left\{l^+,l^-\right\}$ contains both $-\frac{\pi}{6}$ and $-\frac{5\pi}{6}$, and the fact implies the result in $\text{Theorem 2}$ for Case 1.1 (see Fig. 16).
	\begin{figure}[!h]
	\includegraphics[width=60mm]{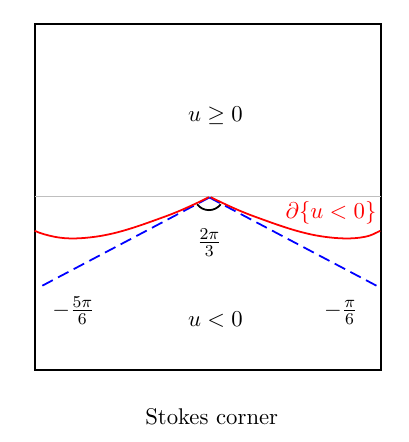}
	\caption{The singular profile of Case 1.1}
\end{figure}

$\bullet$ Asymmetric Stokes corner.

We now consider the case $M_{x^0,u^-}(0^+)=\frac{1}{2}$. Following Step 1, we claim that
$$\mathfrak{L}^+\,\,\text{and}\,\,\mathfrak{L}^-\,\,\text{are subsets of}\,\,\left\{-\pi,-\frac{\pi}{3},0\right\}$$
or
$$ \mathfrak{L}^+\,\,\text{and}\,\,\mathfrak{L}^-\,\,\text{are subsets of}\,\,\left\{-\pi,-\frac{2\pi}{3},0\right\}.$$
Suppose that there exists $0\ne t_m\rightarrow0$ as $m\rightarrow\infty$ such that $\arg\sigma(t_m)\rightarrow\theta^*$, where $$\theta^*\in\left(\mathfrak{L}^+\cup\,\mathfrak{L}^-\backslash\left\{-\pi,-\frac{\pi}{3},0\right\}\right)\quad\text{or}\quad\theta^*\in\left(\mathfrak{L}^+\cup\,\mathfrak{L}^-\backslash\left\{-\pi,-\frac{2\pi}{3},0\right\}\right).$$
 Hence, for each $\rho>0$ such that $B^*=B_\rho\left(\cos\theta^*,\sin\theta^*\right)$ satisfies
$$
B^*\cap\left(\left\{(x,0)\left|\,x\in\mathbb{R}\right.\right\}\cap\left\{\left.\left(x,-\sqrt{3}\,|x|\right)\right|\,x\in\mathbb{R}\right\}\right)=\varnothing,
$$
we observe that \eqref{3.25} remains valid. By Step 2, the limits $l^+$ and $l^-$ exist and take values in either the set $\left\{-\pi,-\frac{2\pi}{3},0\right\}$ or the set $\left\{-\pi,-\frac{\pi}{3},0\right\}$. Consequently, if 
$$M_{x^0,u^-}(0^+)=\int_{ B_1}(-x_2)\chi_{\left\{\left(\rho\cos\theta,\rho\sin\theta\right)\left|\,-\pi<\theta<-\frac{\pi}{3}\right.\right\}}dx,$$
 we obtain that $\Delta u^-_\infty\left(B_{1/10}\left(1/2,-\sqrt{3}/2\right)\right)>0$ and $\left\{l^+,l^-\right\}$ contains $-\frac{\pi}{3}$. By the homogeneity of $u_\infty^-$, it is straightforward to show that $l^-=-\pi$ and $l^+=-\frac{\pi}{3}$. Similarly, if $$M_{x^0,u^-}(0^+)=\int_{ B_1}(-x_2)\chi_{\left\{\left(\rho\cos\theta,\rho\sin\theta\right)\left|\,-\frac{2\pi}{3}<\theta<0\right.\right\}}dx,$$ then $l^-=-\frac{2\pi}{3}$ and $l^+=0$ (see Fig. 17).
 \begin{figure}[!h]
 	\includegraphics[width=115mm]{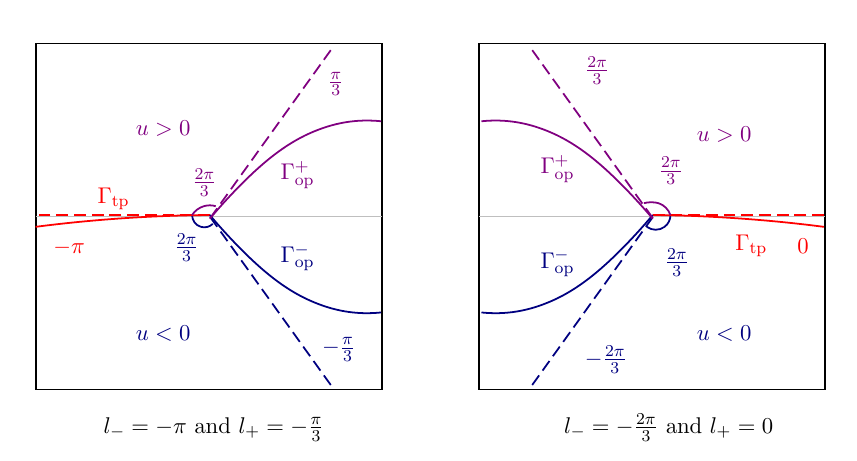}
 	\caption{The singular profile of Case 1.2}
 \end{figure}
 
 $\bullet$ Cusp.

 In Case 1.3, the electric field effect dominates, which completely disrupts the corner structure observed in the other scenarios. Consequently, the singular profile of the free boundary $\partial\{u<0\}$ near $x^0$ must be a cusp.
 
From $\text{Proposition 3.12}$ $$\lim\limits_{m\rightarrow\infty}\int_{ B_1\cap\{x_2<0\}}\chi_{\{u_m<0\}}dx\in\left\{0,\frac{\pi}{2}\right\}.$$ 
We now demonstrate that the singular profiles of the free boundary $\partial\{u<0\}$ is cusp by proving  \begin{equation}
\label{3.27}
\lim\limits_{m\rightarrow\infty}\int_{ B_1\cap\{x_2<0\}}\chi_{\{u_m<0\}}dx=0.\end{equation}
  
Suppose, for the sake of contradiction, that
$$\lim\limits_{m\rightarrow\infty}\int_{ B_1\cap\{x_2<0\}}\chi_{\{u_m<0\}}dx=\frac{\pi}{2}.$$ For brevity of the proof, assume the stagnation point is at the origin $O=(0,0)$. Then there exist $r_0>0$, $\theta_1$ and $\theta_2$ such that $$\Delta u=0\quad\text{in}\,\,\{u<0\}\cap B_{r_0}\quad\text{and}\quad\overline{G}\backslash\{O\}\subset\{u<0\}\cap B_{r_0},$$ 
where $-\pi\leqslant\theta_1<\theta_2\leqslant0$ and $\theta_2-\theta_1>\frac{\pi}{\alpha}$ and $$G=\left\{\left(\rho\cos\theta,\rho\sin\theta\right)\left|\,0<\rho<r_0,\theta_1<\theta<\theta_2\right.\right\}.$$
 Following the strategy in $\text{Proposition 6.5}$ of \cite{VW2}, there exists $C_0>0$ such that
$$
u^-(0,x_2)\geqslant C_0|x_2|^\gamma\quad\text{for all}\,\,x_2\in(0,r_0),
$$
where $\gamma=\frac{\pi}{\theta_2-\theta_1}$. It follows from $\theta_2-\theta_1>\frac{\pi}{\alpha}$ that $1<\gamma<\alpha$, as $1<\alpha<3/2$. This contradicts the following estimate
$$
u^-(0,x_2)\leqslant C|x_2|^\alpha,
$$
which is a consequence of the assumed growth condition $\left|\nabla u\right|\leqslant C \left|x\right|^{\alpha-1}$ in $B_{r_0}$}. 

Moreover, recall that $$
\mathfrak{L}^\pm=\left\{\theta^*\in\left[-\pi,0\right]\left|\,\text{there is}\,\, t_m\rightarrow0^\pm\,\,\text{such that}\,\,\arg\left(\sigma(t_m)-\sigma(0)\right)\rightarrow\theta^*\,\,\text{as}\,\,m\rightarrow\infty\right. \right\}.
$$
By the continuity of $\partial\{u<0\}$, the sets $\mathfrak{L^+}$ and $\mathfrak{L^-}$ are connected. Now define 
$$
\theta_+=\sup\mathfrak{L^+},\quad\theta_-=\inf\mathfrak{L^-},\quad\text{with}\,\,\theta_-\leqslant\theta_+.
$$

We claim that $\theta_+=\theta_-$. If not, then by definition there exist $t_m^+\rightarrow0^+$ and $t_m^-\rightarrow0^-$ such that $\arg\left(\sigma(t_m^\pm)-\sigma(0)\right)\rightarrow \theta_\pm$ and  $\theta_-<\theta_+$. Set $\delta=\frac{\theta_+-\theta_-}{3}>0$ and
$$
r_m=\min\left\{\left|\sigma(t_m^+)-\sigma(0)\right|,\left|\sigma(t_m^-)-\sigma(0)\right|\right\}\quad\text{and}\quad u_m=\frac{u(x^0+r_mx)}{r_m^\alpha}.
$$
Choose any $\bar{\theta}\in I=\left(\theta_-+\delta,\theta_+-\delta\right)\subset B_1\cap\left\{x_2<0\right\}$. Then take $\varepsilon<\frac{\delta}{4}$ and $\rho>0$ small enough such that $u_m<0$ in $D$ (see Fig. 18), where
$$
D=\left\{x\in A\left|\,\arg(x)\in\left(\bar{\theta}-\varepsilon,\bar{\theta}+\varepsilon\right)\right.\right\}\subset I,\quad A=\left\{x\left|\,\frac{\rho}{2}<|x|<\rho\right.\right\}.
$$
 \begin{figure}[!h]
	\includegraphics[width=120mm]{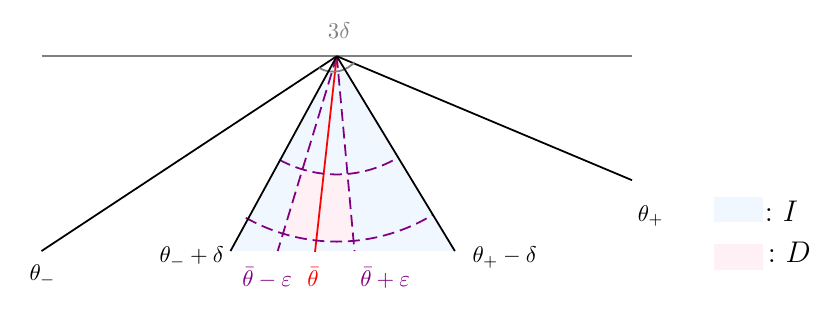}
	\caption{The construction of $D$}
\end{figure}

Consequently,
$$
\int_{  B_1\cap\{x_2<0\}}\chi_{\{u_m<0\}}dx\geqslant\int_{  D}\chi_{\{u_m<0\}}dx\geqslant C(\varepsilon)>0,
$$
where $C(\varepsilon)$ is independent of $m$. This contradicts \eqref{3.27}.
\section*{Appendix A}

 Recalling the free boundary problem \eqref{1.9}, we define the functional $E_{\text{EHD}}(u,\Omega)$ as
$$
E_{\text{EHD}}(u,\Omega)=\int_{\Omega}\left(\left|\nabla u\right|^2+\left(x_2-x_2^0\right)^+\chi_{\{u>0\}}+\left(\left(x_2-x_2^0\right)^++x_2^0-x_2\right)\chi_{\{u<0\}}\right)dx.
$$

We now introduce the definition of the weak solution of the free boundary problem \eqref{1.9}.
\begin{definitionA}
	A function $u\in W^{1,2}(\Omega)$ is said to be a weak solution of \eqref{1.9}, provided that the following conditions hold.
	
	\item[(1)] $u\in C\left(\Omega\right)\cap C^2\left(\Omega\cap\{u\ne0\}\right)$.
	
	\item[(2)] The first variation with respect to domain variations of the functional
	$$
	E_{\text{EHD}}(v,\Omega)=\int_{\Omega}\left(\left|\nabla v\right|^2+\left(x_2-x_2^0\right)^+\chi_{\{v>0\}}+\left(\left(x_2-x_2^0\right)^++x_2^0-x_2\right)\chi_{\{v<0\}}\right)dx
	$$
	vanishes at $v=u$. That is
	\begin{equation}
		\nonumber
		\begin{aligned}
			0
			=&\int_{\Omega}\left(\left|\nabla u\right|^2\mathrm{div}\phi-2\nabla uD\phi\nabla u\right)dx+\int_{\Omega\cap\{x_2<x_2^0\}}\left(\left(x_2^0-x_2\right)\chi_{\{u<0\}}\mathrm{div}\phi-\chi_{\{u<0\}}\phi_2\right)dx\\
			&+\int_{\Omega\cap\{x_2>x_2^0\}}\left(\left(x_2-x_2^0\right)\chi_{\{u>0\}}\mathrm{div}\phi+\chi_{\{u>0\}}\phi_2\right)dx
		\end{aligned}
	\end{equation}
	for each $\phi(x)=\left(\phi_1(x),\phi_2(x)\right)\in C^1_0\left(\Omega,\,\mathbb{R}^2\right)$.
	
	\item[(3)]  The topological free boundary $\partial\{u<0\}\cap\Omega$ is locally a $C^{2,\alpha}$-smooth curve.
\end{definitionA}

\begin{remark}
	Applying integration by parts in the sense of distributions, we obtain
	\begin{equation}
		\nonumber
		\begin{aligned}
			0=&2\int_{\Omega}\left(-\nabla u D^2 u\phi+\nabla u D^2 u \phi+\left(\phi,\nabla u\right)\Delta u\right)dx\\
			&+\int_{\Gamma_{\text{\rm tp}}}\left|\nabla u^-\right|^2\left(\phi,\nu^-\right)-2\left(\nabla u^-,\phi\right)\left(\nabla u^-,\nu^-\right)dS\\
			&+\int_{\Gamma_{\text{\rm tp}}}\left|\nabla u^+\right|^2\left(\phi,\nu^+\right)-2\left(\nabla u^+,\phi\right)\left(\nabla u^+,\nu^+\right)+\left(x_2^0-x_2\right)\left(\phi,\nu^-\right)dS\\
			&+\int_{\Gamma_{\text{\rm op}}^-}\left|\nabla u^-\right|^2\left(\phi,\nu^-\right)-2\left(\nabla u^-,\phi\right)\left(\nabla u^-,\nu^-\right)+\left(\left(x_2-x_2^0\right)^++x_2^0-x_2\right)\left(\phi,\nu^-\right)dS\\
			&+\int_{\Gamma_{\text{\rm op}}^+}\left|\nabla u^+\right|^2\left(\phi,\nu^+\right)-2\left(\nabla u^+,\phi\right)\left(\nabla u^+,\nu^+\right)+\left(x_2-x_2^0\right)^+\left(\phi,\nu^+\right)dS,
		\end{aligned}
	\end{equation}	
	where $\nu^-$ and $\nu^+$ denote the outer normal vectors. Substituting $\phi=\xi\nu^-$ into the integral equation, where $\xi$ is a smooth scalar test function with support in a neighborhood of $\Gamma_{\text{\rm tp}}$ and $\Gamma_{\text{\rm op}}^-$, respectively. Similarly, substituting $\phi=\xi\nu^+$, where $\xi$ has support in a neighborhood of $\Gamma_{\text{\rm op}}^+$. By the assumed smoothness, and using $\left(\phi,\nabla u\right)\Delta u=0$ in $\{u\ne0\}$, $\nu^+=-\nu^-$ on $\Gamma_{\text{\rm tp}}$ and $\nabla u^\pm=-\left|\nabla u^\pm\right|\nu^\pm$, we derive
	\begin{equation}
\nonumber
		\begin{cases}
			\begin{matrix}
				\Delta u=0,& \hspace{40.2mm}	 \text{in 
				}\,\,\{u\ne0\},\\
			\end{matrix}\\
			\begin{matrix}
				\left|\nabla u^-\right|^2-\left|\nabla u^+\right|^2=x_2^0-x_2,& \hspace{6.7mm}	 \text{on 
				}\,\,\Gamma_{\text{\rm tp}},\\
			\end{matrix}\\
			\begin{matrix}
				\left|\nabla u^-\right|^2=\left(x_2-x_2^0\right)^++x_2^0-x_2,&	 \text{on 
				}\,\,\Gamma_{\text{\rm op}}^-,\\
			\end{matrix}\\
			\begin{matrix}
				\left|\nabla u^+\right|^2=\left(x_2-x_2^0\right)^+,&	 \hspace{16.7mm}\text{on 
				}\,\,\Gamma_{\text{\rm op}}^+.\\
			\end{matrix}\\
		\end{cases}
	\end{equation}
\end{remark}

\begin{remark}
	The existence of a minimizer for the functional $E_{\text{EHD}}(u,\Omega)$ implies that the set of solutions satisfying (1) and (2) in $\text{Definition A}$ is non-empty (refer to Theorem 1.1 in \cite{ACF1}). On the other hand, based on the elegant work \cite{PSV}, the condition (3) has been verified. Therefore, the weak solution is well-defined.
\end{remark}

 To investigate stagnation point $x^0\in\partial\{u<0\}$, we observe that $u^-$ is indeed a weak solution of a one-phase free boundary problem (see $\text{Proposition 2.1}$). We now define the weak solution in $\Omega\cap\{u\leqslant0\}$.

\begin{definitionB}
	A non-positive function $w\in W^{1,2}(\Omega)$ is called a weak solution of
	\begin{equation}
		\nonumber
		\begin{cases}
			\begin{matrix}
				\Delta w=0,&	\hspace{21.4mm}	\text{\rm in}\,\,\{w<0\},\\
			\end{matrix}\\
			\begin{matrix}
				\left| \nabla w \right|^2=\left(x_{2}^{0}-x_2\right)^+,&		\text{\rm on}\,\,\partial \{w<0\},\\
			\end{matrix}\\
		\end{cases}
	\end{equation}
	if the following conditions are satisfied.
	
	\item[(1)] $w\in C\left(\Omega\right)\cap C^2\left(\Omega\cap\{w<0\}\right)$.
	
	\item[(2)] The first variation with respect to domain variations of the functional
	$$
	E_{\text{op}}(v,\Omega)=\int_{\Omega}\left(\left|\nabla v\right|^2+\left(x_2^0-x_2\right)^+\chi_{\{v<0\}}\right)dx
	$$
	vanishes at $v=w$. That is
	\begin{equation}
	\nonumber
			0=\int_{\Omega}\left(\left|\nabla w\right|^2\mathrm{div}\phi-2\nabla wD\phi\nabla w\right)dx+\int_{\Omega\cap\{x_2<x_2^0\}}\left(\left(x_2^0-x_2\right)\chi_{\{w<0\}}\mathrm{div}\phi-\chi_{\{w<0\}}\phi_2\right)dx,
	\end{equation}
	for each $\phi(x)=\left(\phi_1(x),\phi_2(x)\right)\in C^1_0\left(\Omega,\,\mathbb{R}^2\right)$.
	
	\item[(3)]  The topological free boundary $\partial\{w<0\}\cap\Omega$ is locally a $C^{2,\alpha}$-smooth curve.
\end{definitionB}

Further analysis of the free boundary relies on the study of so-called blowups. Thus, we introduce the following definitions for blowup sequences.
\begin{definitionC}
	Let $r_m>0$ converge to 0 as $m\rightarrow+\infty$ and $x^{\text{\rm st}}=x^0=(x_1^0,x_2^0)$. The blowup sequence is defined as follows.
	
	$\bullet$ In Case 1.1 and 1.2
	$$ u_m=\frac{u(x^0+r_mx)}{r_m^{3/2}},$$ 
	
	$\bullet$ In Case 1.3, as $1<\alpha<\frac{3}{2}$,
	$$ u_m=\frac{u(x^0+r_mx)}{r_m^{\alpha}}.$$ 
\end{definitionC}

To characterize the blowup limits near any stagnation point $x^0\in\partial\{u<0\}$, we introduce a Weiss type monotonicity formula for $u$ in Case 1.1 and 1.2
\begin{equation}
	\nonumber
		\begin{aligned}
M_{x^0,u}(r)=&r^{-3}\int_{ B_r(x^0)}\left(\left|\nabla u\right|^2+\left(x_2-x_2^0\right)^+\chi_{\{u>0\}}+\left(x_2^0-x_2\right)^+\chi_{\{u<0\}}\right)dx-\frac{3}{2} r^{-4}\int_{\partial B_r(x^0)}u^2dS\\
=&r^{-3}\int_{ B_r(x^0)}\left(\left|\nabla u^+\right|^2+\left(x_2-x_2^0\right)^+\chi_{\{u>0\}}\right)dx-\frac{3}{2} r^{-4}\int_{\partial B_r(x^0)}\left(u^+\right)^2dS\\
&+r^{-3}\int_{ B_r(x^0)}\left(\left|\nabla u^-\right|^2+\left(x_2^0-x_2\right)^+\chi_{\{u<0\}}\right)dx-\frac{3}{2} r^{-4}\int_{\partial B_r(x^0)}\left(u^-\right)^2dS\\
=&M_{x^0,u^+}(r)+M_{x^0,u^-}(r).
	\end{aligned}
\end{equation}
Additionally, we consider a modified Weiss type monotonicity formula for $u$ in Case 1.3 
	\begin{equation}
	\nonumber
	\begin{aligned}
M_{x^0,u,\alpha}(r)=&r^{-2\alpha}\int_{ B_r(x^0)}\left(\left|\nabla u\right|^2+\left(x_2-x_2^0\right)^+\chi_{\{u>0\}}+\left(x_2^0-x_2\right)^+\chi_{\{u<0\}}\right)dx-\alpha r^{-2\alpha-1}\int_{\partial B_r(x^0)}u^2dS\\
=&r^{-2\alpha}\int_{ B_r(x^0)}\left(\left|\nabla u^+\right|^2+\left(x_2-x_2^0\right)^+\chi_{\{u>0\}}\right)dx-\alpha r^{-2\alpha-1}\int_{\partial B_r(x^0)}\left(u^+\right)^2dS\\
&+r^{-2\alpha}\int_{ B_r(x^0)}\left(\left|\nabla u^-\right|^2+\left(x_2^0-x_2\right)^+\chi_{\{u<0\}}\right)dx-\alpha r^{-2\alpha-1}\int_{\partial B_r(x^0)}\left(u^-\right)^2dS\\
=&M_{x^0,u^+,\alpha}(r)+M_{x^0,u^-,\alpha}(r).
		\end{aligned}
\end{equation}

In the Section 2, we have established that $M_{x^0,u}(0^+)$ and $M_{x^0,u,\alpha}(0^+)$ exist as $r\rightarrow0$. Specifically,
	\begin{equation}
	\nonumber
	\begin{aligned}
M_{x^0,u}(0^+)=&\underset{m\rightarrow\infty}{\lim}r_m^{-3}\int_{ B_{r_m}(x^0)}\left(x_2-x_2^0\right)^+\chi_{\{u>0\}}+\left(x_2^0-x_2\right)^+\chi_{\{u<0\}}dx\\
=&\underset{m\rightarrow\infty}{\lim}\int_{ B_1}\left(\left(x_2\right)^+\chi_{\{u_m>0\}}+\left(-x_2\right)^+\chi_{\{u_m<0\}}\right)dx,
		\end{aligned}
\end{equation}
and
\begin{equation}
	\nonumber
	\begin{aligned}
M_{x^0,u,\alpha}(0^+)=&\underset{m\rightarrow\infty}{\lim}r_m^{-2\alpha}\int_{ B_{r_m}(x^0)}\left(x_2-x_2^0\right)^+\chi_{\{u>0\}}+\left(x_2^0-x_2\right)^+\chi_{\{u<0\}}dx\\
=&\underset{m\rightarrow\infty}{\lim}r_m^{3-2\alpha}\int_{ B_1}\left(\left(x_2\right)^+\chi_{\{u_m>0\}}+\left(-x_2\right)^+\chi_{\{u_m<0\}}\right)dx.
		\end{aligned}
\end{equation}
 In particular, we define the density of the negative phase as
$$
M_{x^0,u^-}(0^+)=\underset{m\rightarrow\infty}{\lim}\int_{ B_1}\left(-x_2\right)^+\chi_{\{u_m<0\}}dx,
$$
and
$$
M_{x^0,u^-,\alpha}(0^+)=\underset{m\rightarrow\infty}{\lim}r_m^{3-2\alpha}\int_{ B_1}\left(-x_2\right)^+\chi_{\{u_m<0\}}dx=0.
$$

	\section*{Appendix B}

Proceeding analogously to Case 2 in Subsection 1.4, we obtain that $\left|\nabla u^+(x^{\text{\rm st}})\right|>0$. Based on this, we will show that the free boundary $\partial\{u<0\}$ is a locally $C^{1,\beta}$ curve in $\Omega\cap\{x_2>x_2^0\}$, as $0<\beta<1$.

Following the arguments in the proof of $\text{Proposition 2.2 and 2.4}$, we immediately obtain two lemmas.

\begin{lemma}
	Let $r_0=\frac{1}{2}\text{dist}\left(x^{\text{\rm st}},\partial\Omega\right)$ and $u$ be a weak solution according to $\text{Definition A}$ in $\text{Appendix A}$. Consider the  Weiss-type boundary adjusted energy for $u$ in Case 2
	$$
M_{x^{\text{\rm st}},u,1}(r)=r^{-2}I(r)- r^{-3}J(r).
	$$
	Then for almost every $r\in(0,r_0)$ and $x_2^{\text{\rm st}}>x_2^0$
	\begin{equation}
		\nonumber
		\frac{d}{dr}M_{x^{\text{\rm st}},u,1}(r)=2r^{-2}\int_{\partial B_r(x^{\text{\rm st}})}{\left(\left(\nabla u^+,\nu\right)-\frac{1}{r}u^+\right)^2+\left(\left(\nabla u^-,\nu\right)-\frac{1}{r}u^-\right)^2dS}+K(r),
	\end{equation}
	where
	$$
	K(r)=r^{-3}\int_{ B_r(x^{\text{\rm st}})\cap\{x_2>x_2^0\}}{\left(x_2-x_2^{\text{\rm st}}\right)\chi_{\{u<0\}}dx}+r^{-3}\int_{ B_r(x^{\text{\rm st}})\cap\{x_2<x_2^0\}}{\left(x^{\text{\rm st}}_2-x_2\right)\chi_{\{u<0\}}dx}.
	$$
\end{lemma}

\begin{lemma}
	Let $u$ be a weak solution according to $\text{Definition A}$ in $\text{Appendix A}$. Suppose $x^{\text{\rm st}}=(x_1^{\text{\rm st}},x^{\text{\rm st}}_2)\in\partial\{u<0\}$ with $x^{\text{\rm st}}_2>x_2^0$. Then
	
	\item[(1)] The limit $M_{x^{\text{\rm st}},u,1}(0^+)=\underset{r\rightarrow0^+}{\lim}M_{x^{\text{\rm st}},u,1}(r)$ exists and is finite.
	
	\item[(2)] Let $0<r_m\rightarrow0^+$ as $m\rightarrow\infty$ be a sequence such that the blowup sequence
	$$
	u_m(x)=\frac{u(x^{\text{\rm st}}+r_mx)}{r_m}
	$$
	converges weakly in $W_{\text{\rm loc}}^{1,2}(\mathbb{R}^2)$ to a blowup limit $u_\infty$. Then $u_\infty$ is a homogeneous function of degree $1$ in $\mathbb{R}^2$, i.e., $u_\infty(\lambda x)=\lambda u_\infty(x)$ for all $x\in\mathbb{R}^2$ and all $\lambda>0$.
	
	\item[(3)] $u_m$ converges strongly in $W_{\text{\rm loc}}^{1,2}(\mathbb{R}^2)$.
	
	\item[(4)] For $x^{\text{\rm st}}_2>x^0_2$
	$$
M_{x^{\text{\rm st}},u,1}(0^+)=\lim\limits_{m\rightarrow\infty}\int_{ B_1}\left(x^{\text{\rm st}}_2-x_2^0\right)\chi_{\{u_m>0\}}dx.
	$$
\end{lemma}
\begin{pf}
	Note that 
	\begin{equation}
		\nonumber
		\begin{aligned}
			K(r)=&r^{-3}\int_{ B_r(x^{\text{\rm st}})\cap\{x_2>x_2^0\}}{\left(x_2-x_2^{\text{\rm st}}\right)\chi_{\{u<0\}}dx}+r^{-3}\int_{ B_r(x^0)\cap\{x_2<x_2^0\}}{\left(x^{\text{\rm st}}_2-x_2\right)\chi_{\{u<0\}}dx}\\
			\leqslant& Cr^{-3}r^3	\leqslant C,
		\end{aligned}
	\end{equation}
	which implies that $K(r)$ is integrable. Hence, $\text{Lemma 3.14}$ can be obtained through a proof analogous to $\text{Lemma 4.4}$ in \cite{VW2}.
\end{pf}

Furthermore, following the ideas in \cite{EG} and \cite{VW1}$-$\cite{VW3}, we can fully characterize  $u^+_\infty$ as $x^{\text{\rm st}}_2>x_2^0$.

\begin{proposition}
	{\rm (Classification of $u_\infty^+$)} For $u_m$ defined in $\text{Lemma 3.14}$. Then $u^+_\infty$ has one of the following forms.
	\item[(1)] $u_\infty^+$ has one connected component and there exists $e=\left(e_1,e_2\right)$ such that 
	$$
	u^+_\infty(x)=\sqrt{x^{\text{\rm st}}_2-x_2^0}\left(x\cdot e\right)^+.
	$$
	\item[(2)] $u_\infty^+$ has two connected components, i.e., there exist a constant $\gamma>0$ and $e=\left(e_1,e_2\right)$ such that
	$$
	u^+_\infty(x)=\gamma\left|x\cdot e\right|.
	$$
	\item[(3)] $u_\infty^+\equiv0$.
\end{proposition}
\begin{pf}
	Consider the blowup sequence $u_m=u(x^{\text{\rm st}}+r_mx)/r_m$ as in $\text{Lemma 3.14}$. Note that $B_1\cap\left\{x_2<\frac{x_2^0-x_2^{\text{\rm st}}}{r_m}\right\}=\varnothing$ for $r_m\in(0,r_0)$, where $r_0=x_2^{\text{\rm st}}-x_2^0$, which implies that 
	$$
	\int_{B_1\cap\left\{x_2<\frac{x_2^0-x_2^{\text{\rm st}}}{r_m}\right\}}(x_2^0-x_2^{\text{\rm st}}-r_mx_2)\chi_{\{u_m<0\}}dx\rightarrow0\quad\text{as}\,\,m\rightarrow\infty.
	$$
	Moreover, the strong convergence of $u_m$ to $u_\infty $ in $W^{1,2}_{\text{\rm loc}}\left(\mathbb{R}^2\right)$ and the compact embedding from BV into $L^1$ imply that $u_\infty$ satisfies
	\begin{equation}
\label{6.1}
		0=\int_{\mathbb{R}^2}\left(\left|\nabla u_\infty \right|^2\mathrm{div}\phi-2\nabla u_\infty D\phi\nabla u_\infty+\left(x^{\text{\rm st}}_2-x_2^0\right)\chi_0\mathrm{div}\phi\right)dx,
	\end{equation}
	for any $\phi\in C^1_0\left(\mathbb{R}^2,\mathbb{R}^2\right)$, where $\chi_0$ is the strong $L^1_{\text{loc}}$ limit of $\chi_{\{u_m>0\}}$. Using the facts that  
	$$
	\chi_0\in\{0,1\}\quad\text{a.e. }\quad\text{and}\quad\int_{ B_r(x^0)}u_m^+\left(1-\chi_{\{u_m>0\}}\right)dx=0,
	$$
	the locally uniform convergence of $u_m$ implies  
	\begin{equation}
		\label{6.2}
		\chi_0=1\quad\text{in}\,\,\{u_\infty >0\}.
	\end{equation}
	Now, by choosing a test function $\phi$ whose support lies inside the interior of $\{u_\infty^+ =0\}$, we obtain that 
	\begin{equation}
		\label{6.6}
		\int_{\mathbb{R}^2}\left(x^{\text{\rm st}}_2-x_2^0\right)\chi_0\mathrm{div}\phi dx=0.
	\end{equation}
	Hence, $\chi_0$ is a constant in the interior of $\{u_\infty^+ =0\}$. 
	
	Since $u_\infty $ is a homogeneous function of degree 1, we complete this proof in two steps. 
	
	{\bf Step 1.} First, suppose $\{u_\infty >0\}=\varnothing$.
	
	Combining \eqref{6.2} with \eqref{6.6} shows that $\chi_0=0$ is either $0$ or $1$ in the interior of $\left\{u_\infty^+ =0\right\}$. Consequently, the density takes one of two forms
	$$
	 M_{x^\text{\rm st},u,1}(0^+)=\begin{cases}
	 	\begin{matrix}
	 	0,&\hspace{52.77mm}	\text{if}\,\,\chi_0=0,	\\
	 	\end{matrix}\\
	 		\begin{matrix}
	 		\int_{B_{1}}\left(x^{\text{\rm st}}_2-x_2^0\right)dx=\left(x^{\text{\rm st}}_2-x_2^0\right)\omega_2,&		\text{if}\,\,\chi_0=1,	\\
	 	\end{matrix}\\
	 	\end{cases}
	$$
	where $\omega_2$ is the volume of the unit ball in $\mathbb{R}^2$
	
	{\bf Step 2.} $\{u_\infty >0\}\ne\varnothing$.
	
Since $\left|\nabla u^-\left(x^{\text{ \rm st}}\right)\right|=0$, the scaling gives $$\int_{\mathbb{R}^2}\left|\nabla u_m^-\right|^2dx=\int_{\mathbb{R}^2}\left|\nabla u^-\left(x^{\text{ \rm st}}+r_mx\right)\right|^2dx\rightarrow 0\quad\text{as}\,\, m\rightarrow\infty.$$
Together with the strong convergence of $u_m$ to $u_\infty$ in $W^{1,2}_\text{\rm loc}(\mathbb{R}^2)$, this implies
$$
\int_{\mathbb{R}^2}\left|\nabla u_\infty^-\right|^2dx=0.
$$
Hence $\left|\nabla u_\infty^-\right|=0$, and therefore $u_\infty^-=0$.
	
Now, choose a point $y\in\partial\{u_\infty >0\}\backslash\{O\}$. For sufficiently small $\delta>0$, the unit normal $\nu(y)$ to $\partial\{u_\infty >0\}$ is constant on $\partial\{u_\infty >0\}\cap B_\delta(y)$. Substituting $\phi(x)=\eta(x)\nu(y)$, with $\eta\in C^1_0\left(B_\delta(y),\mathbb{R}^2\right)$, into \eqref{6.1} and integrating by parts, we obtain
	\begin{equation}
		\label{6.4}
		\int_{B_\delta(y)\cap\partial\{u_\infty >0\}}\left|\nabla u_\infty\right|^2\left(\eta,\nu \right)dx=	\int_{B_\delta(y)\cap\partial\{u_\infty >0\}}\left(x^{\text{\rm st}}_2-x_2^0\right)\left(1-\chi_0\right)\left(\eta,\nu\right) dx,
 	\end{equation}
	Therefore, similarly to Subsection 3.2, \eqref{6.4} implies
	\begin{equation}
		\label{6.5}
		\left|\nabla u^+_\infty\right|^2=x^{\text{\rm st}}_2-x_2^0\quad\text{on} \,\,\partial\{u_\infty >0\}.
	\end{equation}
	
First, consider the case when $\{u_\infty >0\}$ has exactly one connected component. Solving $u^+_\infty$ of the ordinary differential equation on $\partial B_1$, we obtain $u_\infty^+(x)=\sqrt{x^{\text{\rm st}}_2-x_2^0}\left(x\cdot e\right)^+$ with the corresponding density is $M_{x^0,u,1}(0^+)=\frac{1}{2}\left(x^{\text{\rm st}}_2-x_2^0\right)\omega_2$.
	
	On the other hand, if $\{u_\infty >0\}$ has two connected components, the argument for \eqref{6.5} shows that $\left|\nabla u_\infty^+\right|$ takes the same constant value on each side of $\partial\{u_\infty >0\}$, and the corresponding density is $M_{x^{\text{\rm st}},u,1}(0^+)=\left(x^{\text{\rm st}}_2-x_2^0\right)\omega_2$. 
\end{pf}

To show that there are no stationary points in Case 2, it suffices to prove the following two assertions.

(1) $u_\infty^+(x)=\sqrt{x^{\text{\rm st}}_2-x_2^0}\left(x\cdot e\right)^+$.

(2) The stagnation points set $S^u$ is empty as $x^{\text{\rm st}}_2>x_2^0$.

\begin{proposition}
	For $x^{\text{\rm st}}\in\partial\{u<0\}$ with $x^{\text{\rm st}}_2>x_2^0$, let $u$ be the weak solution of \eqref{1.20}. Then  $u_\infty^+(x)=\sqrt{x^{\text{\rm st}}_2-x_2^0}\left(x\cdot e\right)^+$.
\end{proposition}
\begin{pf}
	Refer to the classification of $u_\infty^+$ in $\text{Proposition 3.15}$, it suffices to prove that $u_\infty^+\equiv0$ and $u^+_\infty(x)=\gamma\left|x\cdot e\right|$ for some $\gamma>0$ are impossible.
	
	(1) Assume that $u_\infty^+\equiv0$. Since $x^{\text{\rm st}}_2>x_2^0$, we can choose $r_0>$ small enough such that $x_2>x_2^0$ in $B_r(x^{\text{\rm st}})$ for all $0<r<r_0$. By $\text{Proposition 1.1}$, this implies $\partial\{u<0\}=\Gamma_{\text{\rm tp}}$ locally in $B_r(x^{\text{\rm st}})$.
	
	Now, it follows from $u_\infty^+\equiv0$ that $\Delta u_\infty^+=0$. Arguing as in the proof of $\text{Theorem 3.3}$, we obtain for any $\varphi\in C_0^\infty(\Omega)$ the identity
	$$
	\Delta u^+(B_r(x^{\text{\rm st}}))=\int_{\partial\{u>0\}\cap B_r(x^{\text{\rm st}})}\varphi\left|\nabla u^+\right|dS.
	$$
On $\Gamma_{\text{\rm tp}}$, we have $\left|\nabla u^+\right|^2=\left|\nabla u^-\right|^2+x_2-x_2^0\geqslant x_2-x_2^0$, while on $\Gamma_{\text{\rm op}}^+$ we have $\left|\nabla u^+\right|^2=\left(x_2-x_2^0\right)^+$. Consider
	$$
	u_m^+(x)=\frac{u^+\left(x^{\text{\rm st}}+r_mx\right)}{r_m}.
	$$
From the lower bounds above we obtain, for some constant $C_0>0$,
	$$
	\left|\nabla u^+_m\right|^2=\left|\nabla u^+(r_mx+x^{\text{\rm st}})\right|^2\geqslant x_2^{\text{\rm st}}-x_2^0+r_mx_2\geqslant C_0\quad \text{on}\,\,\partial_{\text{red}}\{u_m>0\}\cap B_1.
	$$
	Consequently,
	\begin{equation}
		\nonumber
		0\leftarrow\Delta u_m^+(B_1)\geqslant\int_{\partial_{\text{red}\{u_m>0\}\cap B_1}}\left|\nabla u_m^+\right|dS>0\quad \text{as}\,\,m\rightarrow\infty,
	\end{equation}
which is a contradiction. Hence the assumption $u_\infty^+=0$ is impossible.
	
	(2) Assume that $u^+_\infty(x)=\gamma\left|x\cdot e\right|$ for some $\gamma>0$. Observe that $u_\infty^+$ satisfies
	$$
	\begin{cases}
		\begin{matrix}
			\Delta u_\infty^+=0,&\hspace{15.8mm}	\text{in}\,\,\{u_\infty>0\},	\\
		\end{matrix}\\
		\begin{matrix}
			\left|\nabla u_\infty^+\right|^2=x^{\text{\rm st}}_2-x_2^0,&	\text{on}\,\,\partial\{u_\infty>0\},	\\
		\end{matrix}\\
	\end{cases}
	$$
	which implies that $u_\infty^+$ is a solution to a one-phase flow problem. However, a direct computation shows that the density associated with $u^+_\infty(x)=\gamma\left|x\cdot e\right|$ is strictly smaller than $\left(x^{\text{\rm st}}_2-x_2^0\right)\omega_2$, contradicting the lower bound established in $\text{Theorem 4.3}$ of \cite{AC}. Hence, such a form for $u^+_\infty(x)$ is impossible.
\end{pf}

A key observation here is that $\partial\{u<0\}=\Gamma_{\text{\rm tp}}\subset\partial\{u>0\}$. It follows that $\Gamma_{\text{\rm tp}}$ is a locally $C^{1,\beta}$ curve in $\Omega\cap\{x_2>x_2^0\}$ for any $\beta\in(0,1)$. This regularity in turn implies that the set of stagnation points $S^u$ must be empty in $\{x_2>x_2^0\}$, i.e., $x_2^{\text{\rm st}}=x_2^0$.

\begin{theorem}
	For any $x\in S^u$ with $x_2>x_2^0$, let $u$ be the weak solution of \eqref{1.20}. Then the set of stagnation points $S^u=\varnothing$.
\end{theorem}
\begin{pf}
	Suppose not, $S^u$ is not empty in $\Omega\cap\{x_2>x_2^0\}$. Then, there exists a stagnation point $x^{\text{\rm st}}\in S^u$ such that $\left|\nabla u^+(x^{\text{\rm st}})\right|=\sqrt{x^{\text{\rm st}}_2-x_2^0}>0$.
	
	The corresponding blowup limit takes the form $u_\infty^+(x)=\gamma\left(x\cdot e\right)^+$, where $\gamma=\sqrt{x^{\text{\rm st}}_2-x_2^0}$. To proceed, we invoke an $\varepsilon$-regularity result for $u^+$, which can be obtained by adapting the arguments of \cite{PSV}. Specifically, there exists $\varepsilon_0>0$ with the property that if $U_m^+(x)$ satisfies
	$$
	\gamma_1\left(x\cdot e_1-\varepsilon\right)^+\leqslant U_m^+\leqslant\gamma_1\left(x\cdot e_1+\varepsilon\right)^+\quad\text{in}\,\,B_1,
	$$
	for some $\gamma_1\geqslant\gamma>0$ and any $\varepsilon<\varepsilon_0$, then one can find constants $\gamma_2$, $C>0$, a unit vector $e_2$ with $\left|\gamma_1-\gamma_2\right|+\left|e_1-e_2\right|\leqslant C\varepsilon$ and $\rho\in(0,1)$ such that
	$$
	\gamma_2\left(x\cdot e_2-\frac{\varepsilon}{2}\right)^+\leqslant \left(U_m\right)^+_\rho(x)\leqslant\gamma_2\left(x\cdot e_2+\frac{\varepsilon}{2}\right)^+\quad\text{in}\,\,B_1,
	$$
	where $$\left(U_m\right)^+_\rho(x)=\frac{U_m^+(\rho x)}{\rho}=\frac{u^+\left(x^{\text{\rm st}}+r_m\rho x\right)}{ r_m\rho},$$
	and $x\in\Gamma_{\text{\rm tp}}=\partial\{u>0\}\cap\partial\{u<0\}$.
	
We omit the technical details of the proof and only outline the main idea. The complete argument can be found in \cite{PSV}. The crucial step is to introduce  the linearizing sequence 
	$$
	w_m=\frac{U_m^+-\gamma\left(x\cdot e\right)^+}{\varepsilon_m\gamma},
	$$
	and analyze its limiting behavior.
	\begin{itemize}
		\item {\bf Step 1. Boundary Harnack and flatness improvement.}
		
		 Using the partial boundary Harnack lemma (see \cite{PSV}), we prove that $w_m$ converges to a H$\ddot{\rm o}$lder function $w_\infty$. This demonstrates that if $u$ is $\varepsilon$-flat to $\gamma\left(x\cdot e\right)^+$ in $B_1$, it must become $(1-c)\varepsilon$-flat for some $0<c<1$ in $B_{1/2}$. Which mathematically corresponds to the following statement that if
		$$
		\gamma\left(x\cdot e-\varepsilon\right)^+\leqslant u(x)\leqslant\gamma\left(x\cdot e+\varepsilon\right)^+\quad\text{in} \,\,B_1,
		$$
		then 
		$$
		\gamma\left(x\cdot e-(1-c)\varepsilon\right)^+\leqslant u(x)\leqslant\gamma\left(x\cdot e+(1-c)\varepsilon\right)^+\quad\text{in} \,\,B_{1/2}.
		$$
		
		\item {\bf Step 2. The limiting linearized problem.} 
		
		We verify that the limit $w_\infty$ satisfies the linearized problem
		$$
		\begin{cases}
			\begin{matrix}
				\Delta w_\infty=0,&	\hspace{5.63mm}\text{in}\,\,B_1\cap\{x\cdot e>0\},	\\
			\end{matrix}\\
			\begin{matrix}
				\nabla w_\infty\cdot e=0,&	\text{on}\,\,B_1\cap\{x\cdot e=0\}.	\\
			\end{matrix}\\
		\end{cases}
		$$
		
		\item {\bf Conclusion. Regularity of the free boundary.} 
		
		The decay of flatness for $\partial\{u>0\}\cap\partial\{u<0\}$ follows directly from the regularity of $w_\infty$ obtained in Step 1 and 2.
	\end{itemize}
	
In Case 2, $\left|\nabla u^+\right|$ is strictly positive near the stagnation points. Adapting the arguments of $\text{Lemma 2.5}$ in \cite{PSV}, one readily verifies that $u$ is a viscosity solution of \eqref{1.20}. Moreover, the flatness decay of the free boundary $\partial\{u>0\}\cap\partial\{u<0\}$ implies both the uniqueness of the direction vector $e$ and the $C^{1,\beta}$-regularity of $\partial\{u>0\}\cap\partial\{u<0\}$ in a small neighborhood of each $x^{\text{\rm st}}$, for some $\beta\in(0,1)$. 
	
In Case 2, $\Gamma_{\text{\rm op}}^-\cap\{x_2>x_2^0\}=\varnothing$ implies  $\partial\{u<0\}=\Gamma_{\text{\rm tp}}$. Consequently, we obtain that the free boundary $\partial\{u<0\}$ is locally a $C^{1,\beta}$ curve in $\Omega\cap\{x_2>x_2^0\}$ for every $\beta\in(0,1)$. This regularity contradicts the standing assumption $\left|\nabla u^-(x^{\text{\rm st}})\right|=0$ at a stagnation point, thus completing the proof.
\end{pf}

\section*{Appendix C}

{\bf Proof of $\text{Proposition 2.2}$}	We divide the proof into three steps.
	
	{\bf Step 1.} We compute the derivative of $J(r)$ with respect to $r$. 
	
	Let $y=\frac{1}{r}(x-x^0)$. Then
	\begin{equation}
	\nonumber
	\begin{aligned}
	\frac{d}{dr}J(r)&=\frac{d}{dr}\left(\int_{\partial B_1}{r\left(\left(u^+(x^0+ry)\right)^2+\left(u^-(x^0+ry)\right)^2\right)}dS\right)\\
=&\int_{\partial B_1}{\left(u^+(x^0+ry)\right)^2dS}+2r\int_{\partial B_1}{u^+(x^0+ry) \left(\nabla u^+\left(x^0+ry\right),y\right)dS}\\
	&+\int_{\partial B_1}{\left(u^-(x^0+ry)\right)^2dS}+2r\int_{\partial B_1}{u^-(x^0+ry) \left(\nabla u^-\left(x^0+ry\right),y\right)dS}\\
	=&\frac{1}{r}\int_{\partial B_r(x^0)}{\left(u^+\right)^2dS}+2\int_{\partial B_r(x^0)}{u^+\left(\nabla u^+,\nu\right)dS}\\
	&+\frac{1}{r}\int_{\partial B_r(x^0)}{\left(u^-\right)^2dS}+2\int_{\partial B_r(x^0)}{u^-\left(\nabla u^-,\nu\right)dS}.
	\end{aligned}
	\end{equation}
	
	{\bf Step 2.} We next compute the derivative of $I(r)$ with respect to $r$. For small $\tau>0$, define
	\begin{equation}
	\nonumber
	\label{text}
	\eta _{\tau}\left( t \right) =\begin{cases}
	\begin{matrix}
	1,&	\hspace{10.5mm}	0\leqslant t\leqslant r-\tau ,\\
	\end{matrix}\\
	\begin{matrix}
	\frac{r-t}{\tau},&		r-\tau <t<r,\\
	\end{matrix}\\
	\begin{matrix}
	0,&	\hspace{10.5mm}	 r\leqslant t<\infty.\\
	\end{matrix}\\
	\end{cases}
	\end{equation}
	
	Let $\phi_\tau(x)=\eta_{\tau}(|x-x^0|)(x-x^0)$. Since the definition of the weak solution requires test function $\eta(x)\in C^1_0(B_r(x^0),\mathbb{R}^2)$, we replace  $\phi_\tau(x)$ by its standard mollification $\tilde{\phi}_\tau(x)\in C^1_0(B_r(x^0),\mathbb{R}^2)$ in the following calculations.
	
	By the definition of the weak solution,  
	\begin{equation}
	\nonumber
	\begin{aligned}
	0
	=&\int_{\Omega}{\left(\left|\nabla u^+\right|^2+\left(x_2-x^0_2\right)^+\chi_{\{u>0\}}\right)\left(\eta'_{\tau}(|x-x^0|)\left|x-x^0\right|+2\eta_{\tau}(|x-x^0|)\right)dx}\\
	&+\int_{\Omega}{\left(\left|\nabla u^-\right|^2+\left(x^0_2-x_2\right)^+\chi_{\{u<0\}}\right)\left(\eta'_{\tau}(|x-x^0|)\left|x-x^0\right|+2\eta_{\tau}(|x-x^0|)\right)dx}\\
	&-2\int_{\Omega}{\left|\nabla u^+\right|^2\eta_{\tau}(|x-x^0|)dx-2\left(\nabla u^+,\frac{x-x^0}{|x-x^0|}\right)^2\eta'_{\tau}(|x-x^0|)\left|x-x^0\right|dx}\\
		&-2\int_{\Omega}{\left|\nabla u^-\right|^2\eta_{\tau}(|x-x^0|)dx-2\left(\nabla u^-,\frac{x-x^0}{|x-x^0|}\right)^2\eta'_{\tau}(|x-x^0|)\left|x-x^0\right|dx}
	\\&+\int_{\Omega\cap\{x_2>x_2^0\}}\left(x_2-x_2^0\right)\chi_{\{u>0\}}\left(\eta'_{\tau}(|x-x^0|)\left|x-x^0\right|+2\eta_{\tau}(|x-x^0|)\right)dx		\\&+\int_{\Omega\cap\{x_2<x_2^0\}}\left(x_2^0-x_2\right)\chi_{\{u<0\}}\left(\eta'_{\tau}(|x-x^0|)\left|x-x^0\right|+2\eta_{\tau}(|x-x^0|)\right)dx\\
	&+\int_{\Omega\cap\{x_2>x_2^0\}}\chi_{\{u>0\}}\eta_{\tau}(|x-x^0|)\left(x_2-x^0_2\right)dx\\
	&-\int_{\Omega\cap\{x_2<x_2^0\}}\chi_{\{u<0\}}\eta_{\tau}(|x-x^0|)\left(x_2-x^0_2\right)dx.
	\end{aligned}
	\end{equation}
	
	Letting $\tau\rightarrow0$, we obtain
	\begin{equation}
	\nonumber
	\begin{aligned}
	0=&-r\int_{\partial B_r(x^0)}\left(\left|\nabla u^+\right|^2+(x_2-x_2^0)^+\chi_{\{u>0\}}\right)dS\\
	&-r\int_{\partial B_r(x^0)}\left(\left|\nabla u^-\right|^2+\left((x_2-x_2^0)^++x_2^0-x_2\right)\chi_{\{u<0\}}\right)dS\\
	&+2r\int_{\partial B_r(x^0)}{\left(\nabla u^+,\nu\right)^2+\left(\nabla u^-,\nu\right)^2dS}+\int_{B_r(x^0)\cap\{x_2>x_2^0\}}3\left(x_2-x_2^0\right)\chi_{\{u>0\}}dx\\
	&+\int_{B_r(x^0)\cap\{x_2<x_2^0\}}3\left(x_2^0-x_2\right)\chi_{\{u<0\}}dx,
	\end{aligned}
	\end{equation}
where $\nu=\frac{x-x^0}{|x-x^0|}$. Thus, 
	\begin{equation}
	\nonumber
	\begin{aligned}	
	\frac{d}{dr}I(r)=&\int_{\partial B_r(x^0)}{\left(\left|\nabla u^+\right|^2+\left(x_2-x_2^0\right)^+\chi_{\{u>0\}}\right)dS}\\
	&+\int_{\partial B_r(x^0)}{\left(\left|\nabla u^-\right|^2+\left(x_2^0-x_2\right)^+\chi_{\{u<0\}}\right)dS}
	\\=&2\int_{\partial B_r(x^0)}{\left(\nabla u^+,\nu\right)^2+\left(\nabla u^-,\nu\right)^2dS}+\frac{1}{r}\int_{B_r(x^0)\cap\{x_2>x_2^0\}}3\left(x_2-x_2^0\right)\chi_{\{u>0\}}dx\\
	&+\frac{1}{r}\int_{B_r(x^0)\cap\{x_2<x_2^0\}}3\left(x_2^0-x_2\right)\chi_{\{u<0\}}dx.
	\end{aligned}
	\end{equation}
	
	{\bf Step 3.} To obtain the derivative of $M_{x^0,u}(r)$, we combine $\frac{d}{dr}I(r)$ and $	\frac{d}{dr}J(r)$.
	\begin{equation}
	\nonumber
	\begin{aligned}	
		\frac{d}{dr}M_{x^0,u}(r)=&\frac{d}{dr}\left(r^{-3}I(r)-\frac{3}{2} r^{-4}J(r)\right)
	\\=&-3 r^{-4}I(r)+r^{-3}	\frac{d}{dr}I(r)+6r^{-5}J(r)-\frac{3}{2} r^{-4}	\frac{d}{dr}J(r)\\
	=&-3 r^{-4}\int_{B_r(x^0)}{\left(\left|\nabla u\right|^2+\left(x_2-x_2^0\right)^+\chi_{\{u>0\}}+\left(x_2^0-x_2\right)^+\chi_{\{u<0\}}\right)dx}
	\\&+2r^{-3}\int_{\partial B_r(x^0)}{\left(\nabla u^+,\nu\right)^2+\left(\nabla u^-,\nu\right)^2dS}\\
	&+3r^{-4}\left(\int_{B_r(x^0)\cap\{x_2>x_2^0\}}\left(x_2-x_2^0\right)\chi_{\{u>0\}}dx\int_{B_r(x^0)\cap\{x_2<x_2^0\}}\left(x_2^0-x_2\right)\chi_{\{u<0\}}dx\right)
	\\&+6 r^{-5}\int_{\partial B_r(x^0)}{u^2dS}-\frac{3}{2} r^{-5}\int_{\partial B_r(x^0)}{u^2dS}\\
	&-3r^{-4}\int_{\partial B_r(x^0)}{u^+\left(\nabla u^+,\nu\right)+u^-\left(\nabla u^-,\nu\right)dS}.
	\end{aligned}
	\end{equation}
	Since $\Delta u=0$ in $\{u\ne0\}$, for any $\varepsilon>0$, we have
	\begin{equation}
\nonumber
	\int_{B_r(x^0)}{\left(\nabla u^+,\nabla\left(\max\left\{u^+-\varepsilon,0\right\}\right)^{1+\varepsilon}\right)dS}=\int_{\partial B_r(x^0)}{\left(\max\left\{u^+-\varepsilon,0\right\}\right)^{1+\varepsilon}\left(\nabla u^+,\nu\right)dS},
		\end{equation}
		and
	\begin{equation}
	\nonumber
	\int_{B_r(x^0)}{\left(\nabla u^-,\nabla\left(\min\left\{u^-+\varepsilon,0\right\}\right)^{1+\varepsilon}\right)dS}=\int_{\partial B_r(x^0)}{\left(\min\left\{u^-+\varepsilon,0\right\}\right)^{1+\varepsilon}\left(\nabla u^-,\nu\right)dS}.
\end{equation}
	As $\varepsilon\rightarrow0$
	$$
	\int_{B_r(x^0)}{\left|\nabla u^+\right|^2dx}=\int_{\partial B_r(x^0)}{u^+\left(\nabla u^+,\nu\right)dS},
	$$
	and
		$$
	\int_{B_r(x^0)}{\left|\nabla u^-\right|^2dx}=\int_{\partial B_r(x^0)}{u^-\left(\nabla u^-,\nu\right)dS}.
	$$
	
	Therefore,
	\begin{equation}
	\nonumber
	\begin{aligned}	
		\frac{d}{dr}M_{x^0,u}(r)=&-6 r^{-4}\int_{\partial B_r(x^0)}{u^+\left(\nabla u^+,\nu\right)dS}-6 r^{-4}\int_{\partial B_r(x^0)}{u^-\left(\nabla u^-,\nu\right)dS}\\
		&+2r^{-3}\int_{\partial B_r(x^0)}{\left(\nabla u^+,\nu\right)^2dS}+2r^{-3}\int_{\partial B_r(x^0)}{\left(\nabla u^-,\nu\right)^2dS}\\
		&+\frac{9}{2} r^{-5}\int_{\partial B_r(x^0)}{\left(u^+\right)^2dS}+\frac{9}{2} r^{-5}\int_{\partial B_r(x^0)}{\left(u^-\right)^2dS}
	\\=&2r^{-3}\int_{\partial B_r(x^0)}{\left(\left(\nabla u^+,\nu\right)-\frac{3}{2r}u^+\right)^2+\left(\left(\nabla u^-,\nu\right)-\frac{3}{2r}u^-\right)^2dS},
	\end{aligned}
	\end{equation}
	which concludes the proof.$\hfill\Box$
	
{\bf Proof of $\text{Proposition 2.3}$} Reviewing the proof of $\text{Proposition 2.2}$, it is straightforward to derive that
		\begin{equation}
			\nonumber
				\begin{aligned}
				\frac{d}{dr}I(r)=&2\int_{\partial B_r(x^0)}\left(\nabla u^+, \nu\right)^2dS+\frac{3}{r}\int_{ B_r(x^0)}\left(x_2-x_2^0\right)^+\chi_{\{u>0\}}dx\\
				&+2\int_{\partial B_r(x^0)}\left(\nabla u^-, \nu\right)^2dS+\frac{3}{r}\int_{ B_r(x^0)}\left(x_2^0-x_2\right)^+\chi_{\{u<0\}}dx,
					\end{aligned}
		\end{equation}
		and
			\begin{equation}
			\nonumber
			\begin{aligned}	
		\frac{d}{dr}J(r)=&\frac{1}{r}\int_{\partial B_r(x^0)}\left(u^+\right)^2 dS+2\int_{\partial B_r(x^0)}u^+\left(\nabla u^+, \nu\right)dS\\
		&+\frac{1}{r}\int_{\partial B_r(x^0)}\left(u^-\right)^2 dS+2\int_{\partial B_r(x^0)}u^-\left(\nabla u^-, \nu\right)dS.
		\end{aligned}
\end{equation}
		Then for almost every $r\in(0,r_0)$ and $x_2^{\text{\rm st}}=x_2^0$
		\begin{equation}
			\nonumber
			\begin{aligned}	
				\frac{d}{dr}M_{x^0,u,\alpha}(r)=&\frac{d}{dr}\left(r^{-2\alpha}I(r)-\alpha r^{-2\alpha-1}J(r)\right)
				\\=&-2\alpha r^{-2\alpha-1}\int_{B_r(x^0)}\left(\left|\nabla u\right|^2+\left(x_2-x_2^0\right)^+\chi_{\{u>0\}}+\left(x_2^0-x_2\right)^+\chi_{\{u<0\}}\right)dx\\
				&+2r^{-2\alpha}\int_{\partial B_r(x^0)}\left(\nabla u^+, \nu\right)^2dS+\int_{\partial B_r(x^0)}\left(\nabla u^-, \nu\right)^2dS\\
				&+\frac{3}{r}\int_{ B_r(x^0)}\left(\left(x_2-x_2^0\right)^+\chi_{\{u>0\}}+\left(x_2^0-x_2\right)^+\chi_{\{u<0\}}\right)dx\\
				&+\alpha\left(2\alpha+1\right)r^{-2\alpha-2}\int_{\partial B_r(x^0)}{u^2dS}-\alpha r^{-2\alpha-2}\int_{\partial B_r(x^0)}{u^2dS}\\
				&+2\alpha r^{-2\alpha-1}\int_{\partial B_r(x^0)}{u^+\left(\nabla u^+,\nu\right)dS}+2\alpha r^{-2\alpha-1}\int_{\partial B_r(x^0)}{u^-\left(\nabla u^-,\nu\right)dS}\\
				=&2r^{-2\alpha-1}\int_{\partial B_r(x^0)}\left(\nabla u^+,\nu\right)^2dS-4\alpha r^{-2\alpha-1}\int_{\partial B_r(x^0)}u^+\left(\nabla u^+, \nu\right)dS\\
				&+2r^{-2\alpha-1}\int_{\partial B_r(x^0)}\left(\nabla u^-,\nu\right)^2dS-4\alpha r^{-2\alpha-1}\int_{\partial B_r(x^0)}u^-\left(\nabla u^-, \nu\right)dS\\
				&+2\alpha^2r^{-2\alpha-2}\int_{\partial B_r(x^0)}{\left(u^+\right)^2dS}+2\alpha^2r^{-2\alpha-2}\int_{\partial B_r(x^0)}{\left(u^-\right)^2dS}\\
				&+(3-2\alpha)r^{-2\alpha-1}\int_{ B_r(x^0)}\left(\left(x_2-x_2^0\right)^+\chi_{\{u>0\}}+\left(x_2^0-x_2\right)^+\chi_{\{u<0\}}\right)dx\\
				=&2r^{-2\alpha}\int_{\partial B_r(x^0)}{\left(\left(\nabla u^+,\nu\right)-\frac{\alpha}{r}u^+\right)^2+\left(\left(\nabla u^-,\nu\right)-\frac{\alpha}{r}u^-\right)^2dS}\\
		&+(3-2\alpha)r^{-2\alpha-1}\int_{ B_r(x^0)}\left(\left(x_2-x_2^0\right)^+\chi_{\{u>0\}}+\left(x_2^0-x_2\right)^+\chi_{\{u<0\}}\right)dx.
			\end{aligned}
		\end{equation}
$\hfill\Box$

\end{document}